\documentclass[10pt,a4paper,twoside]{amsart}
\usepackage{amsmath}
\usepackage{amsfonts}
\usepackage{amssymb}
\usepackage[applemac]{inputenc}
\usepackage[english]{babel}
\usepackage{verbatim}
\usepackage{graphicx}
\usepackage[all]{xy}
\xyoption{poly}
\usepackage{enumitem}
\usepackage{multirow}
\usepackage{longtable}
\usepackage{booktabs}
\usepackage{stmaryrd}
\usepackage[usenames,dvipsnames]{color}
	\usepackage{amsthm}
	\usepackage{aliascnt}
	\makeatletter
	\renewenvironment{proof}[1][\proofname]{\vspace{.5em}\par\noindent
  	{\scshape #1\@addpunct{.} }\normalfont\ignorespaces
	}{%
  	\qed}
	\makeatother
	\newtheoremstyle{myplain} 
    		{.5em}                    
   		{.5em}                    
    		{\itshape}                   
   		{}                           
   		{\scshape}                   
    		{.}                          
    		{.5em}                       
    		{}  
	\newtheoremstyle{mydefinition} 
    		{.5em}                    
   		{.5em}                    
    		{\normalfont}                   
   		{}                           
   		{\scshape}                   
    		{.}                          
    		{.5em}                       
    		{}  
	\newtheoremstyle{myremark} 
    		{.25em}                    
   		{.25em}                    
    		{\normalfont}                   
   		{}                           
   		{\scshape}                   
    		{.}                          
    		{.5em}                       
    		{}  
	\theoremstyle{myplain}
	\newtheorem{thm}{Theorem}[section]
	\newtheorem{main}{Theorem}
	
	\newaliascnt{main_cor}{main}
	\newtheorem{main_cor}[main_cor]{Corollary}
	\aliascntresetthe{main_cor}
	\newaliascnt{prop}{thm}
	\newtheorem{prop}[prop]{Proposition}
	\aliascntresetthe{prop}
	\newaliascnt{lem}{thm}
	\newtheorem{lem}[lem]{Lemma}
	\aliascntresetthe{lem}

	\newaliascnt{cor}{thm}
	\newtheorem{cor}[cor]{Corollary}
	\aliascntresetthe{cor}
	
	\newaliascnt{conj}{thm}
	\newtheorem{conj}[conj]{Conjecture}
	\aliascntresetthe{conj}
	\theoremstyle{mydefinition}
	\newaliascnt{defn}{thm}
	\newtheorem{defn}[defn]{Definition}
	\aliascntresetthe{defn}
	\theoremstyle{myremark}
	\newaliascnt{notation}{thm}
	\newtheorem{notation}[notation]{Notation}
	\aliascntresetthe{notation}
	\newaliascnt{rem}{thm}
	\newtheorem{rem}[rem]{Remark}
	\aliascntresetthe{rem}
	\newaliascnt{ass}{thm}
	
	\aliascntresetthe{ass}
	\newaliascnt{exmp}{thm}
	\newtheorem{exmp}[exmp]{Example}
	\aliascntresetthe{exmp}
\usepackage{etoolbox}
\AfterEndEnvironment{main}{\noindent\ignorespaces}
\AfterEndEnvironment{main_cor}{\noindent\ignorespaces}
\AfterEndEnvironment{thm}{\noindent\ignorespaces}
\AfterEndEnvironment{prop}{\noindent\ignorespaces}
\AfterEndEnvironment{lem}{\noindent\ignorespaces}
\AfterEndEnvironment{cor}{\noindent\ignorespaces}
\AfterEndEnvironment{rem}{\noindent\ignorespaces}
\AfterEndEnvironment{exmp}{\noindent\ignorespaces}
\AfterEndEnvironment{defn}{\noindent\ignorespaces}
\AfterEndEnvironment{ass}{\noindent\ignorespaces}
\AfterEndEnvironment{notation}{\noindent\ignorespaces}
\usepackage{array}
	\usepackage[bookmarks=false]{hyperref}

\usepackage[capitalise]{cleveref}
	\Crefname{equation}{}{}
	\crefname{equation}{}{}
	\Crefname{thm}{Theorem}{Theorems}
	\Crefname{main}{Theorem}{Theorems}
	\Crefname{main_cor}{Corollary}{Corollaries}
	\Crefname{prop}{Proposition}{Propositions}
	\Crefname{lem}{Lemma}{Lemmata}
	\Crefname{cor}{Corollary}{Corollaries}
	\Crefname{defn}{Definition}{Definitions}
	\Crefname{rem}{Remark}{Remarks}
	\Crefname{exmp}{Example}{Examples}
	\Crefname{notation}{Notation}{Notations}
	\Crefname{ass}{Assumption}{Assumptions}
	\usepackage[msc-links,initials]{amsrefs}
	\hypersetup{colorlinks=true,allcolors = cyan}
\numberwithin{table}{section}
\numberwithin{equation}{section}
\numberwithin{subsection}{section}
\makeatletter 
\def\@cite#1#2{[{\hbox{#1}\if@tempswa , #2\fi }]} 
\makeatother
	\newcommand{\dimalg}{d}
	\newcommand{\halfdim}{h}
	\newcommand{\permiss}{m}
	\newcommand{\element}{a}
	\newcommand{\lastindex}{\ell}
		\newcommand{\level}{r}
		 \newcommand{\otherlevel}{t}
	\newcommand{\parameterOne}{\alpha}
	\newcommand{\parameterTwo}{\beta}

\newcommand{\C}{\mathbb{C}}
\newcommand{\N}{\mathbb{N}}
\newcommand{\R}{\mathbb{R}}
\newcommand{\Z}{\mathbb{Z}}
\newcommand{\Q}{\mathbb{Q}}
\newcommand{\A}{\mathbb{A}}
\newcommand{\F}{\mathbb{F}}

\newcommand{\Mat}[2]{\mathop{\mathrm{Mat}_{#1}#2}}

\newcommand{\tuple}[1]{\mathrm{#1}}
\newcommand{\transpose}[1]{\mathop{#1^{\mathrm{t}}}}
	\DeclareMathOperator{\Hom}{Hom}
	\DeclareMathOperator{\rad}{\mathrm{Rad}}
	\newcommand{\rk}{\mathop{\mathrm{rk}}}
	\DeclareMathOperator{\linspan}{Span}
	\newcommand{\im}{\mathop{\mathrm{im}}}
	\DeclareMathOperator{\tr}{tr}
	\DeclareMathOperator{\id}{id}
	\DeclareMathOperator{\proj}{pr}
	\newcommand{\Cmatrix}[2]{{\mathop{\mathcal{R}}}_{#1}(#2)}
	\newcommand{\cmatrix}{\mathcal{R}}
	\newcommand{\redcmatrix}{\mathcal{S}}

\newcommand{\zetafunc}[1]{\mathop{\zeta_{#1}(s)}}
\newcommand{\abs}[1]{\alpha(#1)}

\newcommand{\pseries}[2]{\mathop{\mathcal{P}_{#1}#2}}
\newcommand{\numpoin}[2]{\mathop{\mathcal{F}_{#1}{(#2)}}}
\newcommand{\denpoin}[2]{\mathop{\mathcal{G}_{#1}{(#2)}}}
\newcommand{\cone}[1]{\mathop{W_{#1}(\completion)}}

\newcommand{\centrseq}{\mathfrak{S}}
\newcommand{\fcoeff}[2]{F_{#1}#2}
\newcommand{\pcoeff}[2]{N_{#1,\mathrm{r}_{#1}}^{\completion}{#2}}

\newcommand{\ccoeffiso}[2]{\mathcal{C}_{#1}{#2}}

\newcommand{\radic}{\mathfrak{r}}
\newcommand{\radiclass}{\mathbf{c}}
\newcommand{\arr}{\mathcal{A}}
\newcommand{\radseq}{\mathcal{K}}
        \newcommand{\redfin}{\sigma}
        \newcommand{\redfunction}[1]{\redfin_{I,\mathrm{r}_I,#1}}

\newcommand{\ring}{R}
\newcommand{\ringint}{\mathcal{O}}
\newcommand{\primeideal}{\mathfrak{p}}
\newcommand{\qfield}{\mathbb{F}_q}
\newcommand{\numfield}{k}
\newcommand{\completion}{\mathfrak{o}}
\newcommand{\localfrac}{\mathfrak{k}}
\newcommand{\qclosure}{\overline{\mathbb{F}}_q}
\DeclareMathOperator{\frob}{\mathrm{Frob}}

\DeclareMathOperator{\spm}{Spec}
	\newcommand{\pfaffians}[2]{P_{#1}#2}
	\newcommand{\radicalminors}[2]{R_{#1}#2}
	\newcommand{\rankvar}[2]{\mathop{\mathrm{V}_{#1}^{#2}}}
	\newcommand{\ranklocus}[3]{\mathop{\mathrm{L}_{#1}^{#2}{#3}}}
	\newcommand{\centrvar}[3]{\mathop{\mathrm{X}_{#1}^{#2}{#3}}}
	\newcommand{\sheet}{\mathrm{S}}
	\newcommand{\crosssection}{C}
	\newcommand{\cs}[2]{C_{#1}{(#2)}}

\newcommand{\choice}{\varphi}
\newcommand{\bl}[2]{\mathop{\mathrm{B}_{#1}(#2)}}
\newcommand{\killing}{\kappa}
\newcommand{\killingiso}{\beta}

\newcommand{\centr}[2]{{\mathop{\mathrm{C}}}_{#1}\left(#2\right)}
\newcommand{\heis}[1]{\mathop{\mathrm{Heis}(#1)}}

	\DeclareMathOperator{\Sh}{\mathrm{Sh}}
	\newcommand{\sh}[2]{\mathop{\Sh_{#1}\mathopen{}\left({#2}\right)\mathclose{}}}
	\newcommand{\lieshadow}{\mathfrak{s}}
        \newcommand{\orbit}{\mathcal{C}}
        
        \newcommand{\alg}[1]{\mathfrak{g}_{#1}}
        \newcommand{\algfin}{\bar{\alg{}}}
        \newcommand{\alginf}{\alg{}}
        \newcommand{\basis}{\mathcal{B}}
		\DeclareMathOperator{\alggroup}{\mathbf{G}}
		\DeclareMathOperator{\otheralggroup}{\mathbf{H}}
			\newcommand{\SL}[3]{\mathop{\mathrm{SL}_{#1}^{#3} #2}}
			\newcommand{\GL}[2]{\mathop{\mathrm{GL}_{#1} #2}}
		\DeclareMathOperator{\lie}{Lie}
		\newcommand{\liealg}{\mathcal{L}}
			\newcommand{\spl}[2]{{\mathop{}\mathfrak{sl}}_{#1}#2}
			\newcommand{\gl}[2]{{\mathop{\mathfrak{gl}}}_{#1}#2}

\newcommand{\partition}{\mathbf{d}}
\newcommand{\partel}{d}
\DeclareMathOperator{\pdual}{D}
\newcommand{\dualpartition}[1]{\pdual(#1)}

\newcommand{\reduc}[2]{\theta_{#1,#2}}

\newcommand{\reducinf}[1]{\theta_{#1}}

\newcommand{\otherbasis}{\mathcal{C}}
\newcommand{\ob}{c}
\DeclareMathOperator{\isol}{iso}
\DeclareMathOperator{\compl}{c}	

\DeclareMathOperator{\asnf}{ASNF}	
		\newcommand{\TwoOneOneSemLetter}{a}
		\newcommand{\TwoOneOneNilpLetter}{b}
		\newcommand{\TwoTwoSemDiagLetter}{c}
		\newcommand{\TwoTwoSemNonLetter}{d}
		\newcommand{\TwoTwoNilpLetter}{e}
		\newcommand{\ThreeOneThreeEVLetter}{j}
		\newcommand{\ThreeOneTwoEVLetter}{f}
		\newcommand{\ThreeOneDiagLetter}{l}
		\newcommand{\ThreeOneNonJorLetter}{i}
		\newcommand{\ThreeOneNilpLetter}{h}
		\newcommand{\regularLetter}{r}
		\newcommand{\TwoTwoNilpPlusSemLetter}{m}
		\newcommand{\TwoTwoNilpPlusNilpLetter}{n}
		\newcommand{\TwoTwoNilpPlusSemNonLetter}{o}
		\newcommand{\TwoOneOneNilpPlusSemDiagLetter}{q}
		\newcommand{\TwoOneOneNilpPlusSemNonLetter}{v}
		\newcommand{\TwoOneOneNilpPlusNilpLetter}{w}
		\newcommand{\RegLetter}{r}
		\newcommand{\SubregLetter}{u}

		\newcommand{\LieTwoOneOneSem}{\mathbf{\expandafter{\TwoOneOneSemLetter}}}
    
		\newcommand{\LieTwoOneOneNilp}{\mathbf{\expandafter{\TwoOneOneNilpLetter}}}
            	
            	\newcommand{\LieTwoTwoSemDiag}{\mathbf{\expandafter{\TwoTwoSemDiagLetter}}}
            
		\newcommand{\LieTwoTwoSemNon}{\mathbf{\expandafter{\TwoTwoSemNonLetter}}}
            
		 \newcommand{\LieTwoTwoNilp}{\mathbf{\expandafter{\TwoTwoNilpLetter}}}
            
		\newcommand{\LieThreeOneThreeEV}{\mathbf{\expandafter{\ThreeOneThreeEVLetter}}}
            
		\newcommand{\LieThreeOneTwoEV}{\mathbf{\expandafter{\ThreeOneTwoEVLetter}}}
            
		\newcommand{\LieThreeOneDiag}{\mathbf{\expandafter{\ThreeOneDiagLetter}}}
            
		\newcommand{\LieThreeOneNonJor}{\mathbf{\expandafter{\ThreeOneNonJorLetter}}}
            
            	\newcommand{\LieThreeOneNilp}{\mathbf{\expandafter{\ThreeOneNilpLetter}}}
	            
		\newcommand{\Lieregular}{\mathbf{\expandafter{\regularLetter}}}
	            
		\newcommand{\LieTwoTwoNilpPlusSem}{\mathcal{\expandafter\mathbf\expandafter{\TwoTwoNilpPlusSemLetter}}}
	            
		\newcommand{\LieTwoTwoNilpPlusNilp}{\mathcal{\expandafter\mathbf\expandafter{\TwoTwoNilpPlusNilpLetter}}}
	            
		\newcommand{\LieTwoTwoNilpPlusSemNon}{\mathcal{\expandafter\mathbf\expandafter{\TwoTwoNilpPlusSemNonLetter}}}
	            
		\newcommand{\LieTwoOneOneNilpPlusSemDiag}{\mathcal{\expandafter\mathbf\expandafter{\TwoOneOneNilpPlusSemDiagLetter}}}
	            
		\newcommand{\LieTwoOneOneNilpPlusSemNon}{\mathcal{\expandafter\mathbf\expandafter{\TwoOneOneNilpPlusSemNonLetter}}}
	            
		\newcommand{\LieTwoOneOneNilpPlusNilp}{\mathcal{\expandafter\mathbf\expandafter{\TwoOneOneNilpPlusNilpLetter}}}
		
		\newcommand{\LieReg}{\mathcal{\expandafter\mathbf\expandafter{\RegLetter}}}
		\newcommand{\LieSubreg}{\mathcal{\expandafter\mathbf\expandafter{\SubregLetter}}}
		
\usepackage[disable]{todonotes}

\title[Poincar\'e series and representation zeta functions]{Poincar\'e series of Lie lattices and representation zeta functions of arithmetic groups}

\author{Michele Zordan}
\address{KU Leuven, Departement Wiskunde, Celestijnenlaan 200B, B-3001 Leuven, Belgium.}
\email{michele.zordan@kuleuven.be}
\begin{document}
\begin{abstract}
We compute explicit formulae for Dirichlet generating functions enumerating finite-dimensional irreducible complex representations of potent and saturable principal congruence subgroups of $\SL{4}{(\completion)}{}$ for $\completion$ a compact DVR of characteristic $0$ and odd residue field characteristic. In doing so we develop a novel method for computing Poincar\'e series associated with commutator matrices of  $\completion$-Lie lattices with finite abelianization and whose rank-loci enjoy an additional smoothness property.\par
We give explicit formulae for the abscissa of convergence of the representation zeta functions of potent and saturable FAb $p$-adic analytic groups whose associated Lie lattices satisfy the hypotheses of the aforementioned method.\par
As a by-product of our computations we find that not all  $4\times 4$ traceless matrices over a  finite quotient of $\completion$ admit shadow-preserving lifts, thus disproving that smooth loci of constant centralizer dimension in $\spl{4}{(\C)}{}$ ensure presence of shadow-preserving lifts for almost all primes as  suggested in \cite[Remark 6.5]{akov2}.
\end{abstract}
\maketitle
\thispagestyle{empty}
\allowdisplaybreaks
\section{Introduction}
\label{sec:intro}
Let $G$ be a group. The main object of investigation of this paper is the asymptotic behaviour of the number  $r_i(G)$ ($i\in\N$) of isomorphism classes of $i$-dimensional irreducible complex representations of $G$.
As done by Avni, Klopsch, Onn and Voll in \cite[Section~1.1]{akov2013representation} we stipulate the following  convention: when $G$ is a topological or an algebraic group, we assume that the representations enumerated by $r_i(G)$ are continuous or rational, respectively. Throughout this work, when using the sequence $\lbrace r_i(G)\rbrace_{i\in\N}$ for a group $G$, we assume that $G$ is (representation) \emph{rigid}, i.e.\ $r_i(G)$ is finite for all $i\in\N$.\par
%
%
%
The function $r_i(G)$ as $i$ varies in $\N$ is called the representation growth function of $G$. If the sequence 
\[
R_N(G)=\sum_{i=1}^{N} r_i(G)\text{ $N\in\N$,}
\] 
is bounded by a polynomial in $N$ as $N$ tends to infinity, the group $G$ is said to have \emph{polynomial representation growth} (PRG). It is a well established tradition to study the asymptotic behaviour of such sequences of non-negative integers by encoding them in a Dirichlet generating function. We therefore define the \emph{representation zeta function} as %
 \[
 \zetafunc{G}=\sum_{i=1}^{\infty} r_i(G) i^{-s},
 \]
 where $s$ is a complex variable. The following definition is crucial in relating the analytical properties of $\zetafunc{G}$ with the asymptotics of $\lbrace r_i(G)\rbrace_{i\in\N}$.
 \begin{defn}
The infimum of all $\alpha\in\R$ such that $\zetafunc{G}$ converges on the complex half-plane $\lbrace s\in\C\mid\Re(s)>\alpha\rbrace$ is called  the \emph{abscissa of convergence} of  $\zetafunc{G}$  (or of $G$ for short) and denoted by $\abs{G}$. 
\end{defn}
In the context of representation growth, the most relevant property of the abscissa of convergence  is that it gives the degree of polynomial growth of $R_N(G)$ as $N$ tends to infinity. Indeed $\log(R_N(G))/\log (N)$ tends to $\abs{G}$ as $N$ tends to infinity.
 \subsection{Main results}
Let $\completion$ be a compact discrete valuation ring of characteristic $0$, prime ideal $\primeideal$ and residue field $\qfield$ of cardinality $q$ and characteristic $p$.
In this paper we focus on representation zeta functions of arithmetic groups and FAb $p$-adic analytic groups arising from them. In order to do so we study Poincar\'e series of certain Lie lattices giving a new method to compute them. This is used in general to prove results on the abscissae of convergence of the representation zeta functions of some FAb $p$-adic analytic groups and in particular to compute the representation zeta functions of the potent and saturable principal congruence subgroups of $\SL{4}{(\completion)}{}$. Key to this latter application is the fact that the loci of constant Lie centralizer dimension in the complexification of $\spl{4}{(\completion)}{}$ are smooth and irreducible  quasi-affine varieties. Under this point of view, the situation is the same for $\mathfrak{su}_{4}{(\completion)}{}$, $\spl{5}{(\completion)}{}$ and $\mathfrak{su}_{5}{(\completion)}{}$. The methods of this paper may therefore be used in the future to compute the representation zeta functions of potent and saturable principal congruence subgroups of $\mathrm{SU}_{4}{(\completion)}$, $\SL{5}{(\completion)}{}$ and $\mathrm{SU}_{5}{(\completion)}$. Broadly speaking, the local factors of the representation zeta function of torsion-free finitely generated nilpotent group are also given by a Poincar\'e series that is very similar to the one we employ here. The same methods of this paper may therefore be modified to work also in the latter context, giving formulae for the abscissa of convergence of those zeta functions as well.\par
Our approach compares to a number of pre-existing results. On the one hand, much like the Igusa-type zeta functions used by Avni, Klopsch, Onn and Voll in \cite{akov2013representation}, our methods only rely on intrinsic information on the Lie lattice and, when they apply, they may be regarded as a way of computing these $p$-adic integrals, avoiding the non-constructiveness caused by Hironaka's theorem on resolutions of singularities. On the other hand, similarly to the method of shadows used by Avni, Klopsch, Onn and Voll in \cite{akov2}, the approach of this paper is constructive. The methods here, however, apply to a wider class of examples  (the method of shadows is bespoke for Lie lattices arising from a reductive algebraic group of type $A_2$ defined over a number field).\par
As in \cite{zor2016adjoint} our method combines \cite{akov2013representation} and \cite{akov2} by expressing the representation zeta function as a Poincar\'e series and then computing it directly without using $p$-adic integration. In addition to this there are two main new ideas. The first one is to use a new smoothness condition (cf.\ \cref{def:g_geo_smooth_rk}) and the property of isolation of a submodule (cf.\ \cref{def:isol}) to reduce the computation of the Poincar\'e series to a problem over $\qfield$. The second new idea is to further simplify the computation breaking up the $\qfield$-Lie algebra associated with the initial $\completion$-Lie lattice into pieces, each of which has a specified behaviour (cf.\  \cref{def:kernel_class}). These concepts will be later clarified as we present the main technical result of this paper, i.e.\ \cref{main:C}.
 \subsubsection{Poincar\'e series of Lie lattices}
The first main result does not directly concern representation zeta functions but Poincar\'e series of Lie lattices with finite abelianization. Let $\alginf$ be an $\completion$-Lie lattice with finite abelianization $\alginf/[\alginf,\alginf]$ or equivalently (by \cite[Proposition~2.1]{akov2013representation}) such that $\alginf\otimes_{\Z_p} \Q_p$ is a perfect Lie algebra. Let $\dimalg$ be the $\completion$-rank of $\alginf$ as a free $\completion$-module and let $\algfin$ be the $\qfield$-Lie algebra obtained from $\alginf$ by reducing modulo $\primeideal$. \par
The Poincar\'e series of $\alginf$ encodes how the elementary divisors of a commutator matrix behave under lifting and is defined as in \cite[(3.3)]{akov2013representation} (see also \cref{sec:tech_main} and \cref{eq:N_sets} therein). The following main result characterizes the abscissa of convergence of the Poincar\'e series of a perfect Lie lattice with geometrically smooth rank-loci (see \cref{def:geo_smooth_rk}).
 \begin{main}
\label{main:abscissa_P}
For $\omega\in\Hom(\algfin,\qfield) \smallsetminus \lbrace 0 \rbrace$ define
	\begin{align*}
		\rad(\omega)		&	=\lbrace x\in\algfin\mid \omega([x,y]) = 0\,\forall y\in\algfin\rbrace\\
		\rho_\omega(s)	 	&	= 2(\dimalg - \dim_{\qfield} [ \rad (\omega), \rad (\omega)]) - (\dimalg - \dim_{\qfield} \rad (\omega)) s.
	\end{align*}
Assume that the rank-loci of $\alginf$ are geometrically smooth. Then the abscissa of convergence of the Poincar\'e series of $\alginf$ is
	\[
		\max \bigcup_{\omega\in\Hom(\algfin,\qfield) \smallsetminus \lbrace 0 \rbrace}\lbrace s\in\Q	\mid	\rho_\omega(s) = 0\rbrace.
	\]
\end{main}
The case bearing the highest interest for us is when $\alginf$ is potent and saturable as a $\Z_p$-Lie lattice. According to the version of the Kirillov orbit method in \cite[Proposition~3.1]{akov2013representation}, in this case, the Lazard correspondence  of \cite[Section~4]{san2007p-saturable} assigns a $p$-adic analytic group $G = \exp(\alginf)$, whose representation zeta function may be expressed in terms of the Poincar\'e series $\mathcal{P}_{\alginf}$ of $\alginf$:
	\[
		\zetafunc{G^\permiss} = \mathcal{P}_{\alginf}(s+2).
	\]
It is worth recalling that given an $\completion$-Lie lattice, \cite[Proposition~2.2]{akov2013representation} ensures that there is $\permiss_0\in\N_0$ such that $\alginf^\permiss = \primeideal^\permiss \alginf$ is potent and saturable as a $\Z_p$-Lie lattice for all $m\in\N$ such that $\permiss > \permiss_0$.  The actual value of $m_0\in\N_0$ depends on $\completion$, but for instance it is $0$ if $\completion$ is unramified over $\Z_p$.\par
When $\alginf$ is also {\em quadratic} i.e.\ it admits a non-degenerate symmetric bilinear form that is associative  with respect to the Lie bracket we obtain the following result:
\begin{main_cor}
\label{main:E}
For $x\in\algfin$ define $\alpha_x(s)\in \Z[s]$ as
\[
\alpha_x(s) = 2(\dim_{\qfield} \centr{\algfin}{x} - \dim_{\qfield} [\centr{\algfin}{x},\centr{\algfin}{x}]) - (\dimalg - \dim_{\qfield} \centr{\algfin}{x}) s.
\]
Let $\alginf$ be potent saturable as a $\Z_p$-Lie lattice. Assume $\alginf$ is quadratic and that the rank-loci of $\alginf$ are geometrically smooth. Then, for $G = \exp(\alginf)$,
\[
\abs{G} = \max \bigcup_{x\in \algfin  \smallsetminus \lbrace 0\rbrace}\lbrace s\in\Q	\mid	\alpha_x(s) = 0\rbrace.
\]
\end{main_cor}
\subsubsection{Special linear groups}
The second main result is an explicit formula for the representation zeta function of the principal congruence subgroups $\SL{4}{(\completion)}{\permiss}$, where $\permiss$ is permissible, i.e.\ such that $\primeideal^\permiss \spl{4}{(\completion)}{}$ is saturable and potent.
\begin{main}
\label{thmB}
Let $\completion$ have odd residue field characteristic. Then, for all permissible $\permiss$, 
\[
\zetafunc{\SL{4}{(\completion)}{\permiss}}=q^{15\permiss}\frac{\mathcal{F}(q,q^{-s})}{\mathcal{G}(q,q^{-s})}
\]
where
\begin{align*}
\mathcal{F}(q,t)=	& \,q t^{18}\\
	& - {\left(q^{7} + q^{6} + q^{5} + q^{4} - q^{3} - q^{2} - q\right)} t^{15}\\
	& + {\left(q^{8} - 2 \, q^{5} - q^{3} + q^{2}\right)} t^{14}\\
	& + {\left(q^{9} + 2 \, q^{8} + 2 \, q^{7} - 2 \, q^{5} - 4 \, q^{4} - 2 \, q^{3} - q^{2} + 2 \, q + 1\right)} t^{13} \\
	&- {\left(q^{10} + q^{9} + q^{8} - 2 \, q^{7} - 2 \, q^{6} - 2 \, q^{5} + 2 \, q^{3} + q^{2} + q\right)} t^{12}\\
	& + {\left(q^{8} + 2 \, q^{6} + q^{4} - q^{3} - q^{2} - q\right)} t^{11}\\
	&+ {\left(q^{8} + q^{7} - 2 \, q^{4} + q\right)} t^{10}\\
	&- {\left(2 \, q^{10} + q^{9} + q^{8} - q^{7} - 3 \, q^{6} - 2 \, q^{5} - 3 \, q^{4} - q^{3} + q^{2} + q + 2\right)} t^{9}\\
	& + {\left(q^{9} - 2 \, q^{6} + q^{3} + q^{2}\right)} t^{8}\\
	& - {\left(q^{9} + q^{8} + q^{7} - q^{6} - 2 \, q^{4} - q^{2}\right)} t^{7}\\
	& - {\left(q^{9} + q^{8} + 2 \, q^{7} - 2 \, q^{5} - 2 \, q^{4} - 2 \, q^{3} + q^{2} + q + 1\right)} t^{6}\\
	& + {\left(q^{10} + 2 \, q^{9} - q^{8} - 2 \, q^{7} - 4 \, q^{6} - 2 \, q^{5} + 2 \, q^{3} + 2 \, q^{2} + q\right)} t^{5}\\
	& + {\left(q^{8} - q^{7} - 2 \, q^{5} + q^{2}\right)} t^{4}\\
	&+ {\left(q^{9} + q^{8} + q^{7} - q^{6} - q^{5} - q^{4} - q^{3}\right)} t^{3}\\
	& + q^{9}\\
\mathcal{G}(q,t)=&\,q^9 {\left(1-q t^3\right)} {\left(1-q t^4\right)}  {\left( 1-q^{2}t^5 \right)} {\left(1-q^{3}t^6 \right)}.
\end{align*}
\end{main}
\begin{rem}
The palindromic symmetry of $\mathcal{F}(q,t)$ in Theorem B implies that $\zetafunc{\SL{4}{(\completion)}{\permiss}}$ satisfies the functional equation of \cite[Theorem~A]{akov2013representation}, e.g.\ when $\permiss\in\N$ is permissible for $\Z_p$:
\begin{equation*}
\zeta_{\SL{4}{(\completion)}{\permiss}}(s)_{\vert q\rightarrow q^{-1}}=q^{-15\cdot\permiss} \cdot \zeta_{\SL{4}{(\completion)}{\permiss}}(s).
\end{equation*}
Simple substitutions reveal that $\zeta_{\SL{4}{(\completion)}{\permiss}}(-2)=0$, in accordance with \cite{sanzapklo2013vanishesarxiv}; while $\mathcal{F}(1,t)=\mathcal{G}(1,t)$.
\end{rem}
In \cite{ros2015topological} T. Rossmann introduces the topological representation zeta function of a torsion-free free nilpotent group. Following his approach one may also define a topological representation zeta function attached to $\zetafunc{\SL{4}{(\completion)}{\permiss}}$. It is indeed possible to compare the properties of this function with the properties of the topological representation zeta function of nilpotent groups proved in \cite{ros2015topological}; the only caveat here is that, in order to account for the differences in the application of the Kirillov orbit method in the two cases (compare \cite[Proposition~3.1, Corollary 3.7]{akov2013representation} and \cite[Theorem~2.6]{stavol2011nilpotent}), one substitutes $s$ with $s-2$. With this in mind one computes
\begin{equation}
\zeta^{\text{top}}_{\SL{4}{}{\permiss}}(s)=\frac{8 \, {\left(15 \, s^{3} + 26 \, s^{2} + 11 \, s - 1\right)} {\left(s + 2\right)}}{{\left(5 \, s - 2\right)} {\left(4 \, s - 1\right)} {\left(3 \, s - 1\right)} {\left(2 \, s - 1\right)}},
\end{equation}
from which it follows that
\begin{equation}
\label{top}
\zeta^{\text{top}}_{\SL{4}{}{\permiss}}(s-2)=\frac{8 \, {\left(15 \, s^{3} - 64 \, s^{2} + 87 \, s - 39\right)} s}{{\left(5 \, s - 12\right)} {\left(4 \, s - 9\right)} {\left(3 \, s - 7\right)} {\left(2 \, s - 5\right)}}.
\end{equation}
One sees that, analogously to  \cite[Proposition~4.5]{ros2015topological}, its limit as $s\rightarrow\infty$ is $1$ and that, analogously to \cite[Proposition~4.8]{ros2015topological}, all its poles are rational and smaller than $15$. The substitution of $s$ with $s-2$ also makes sure that $\zeta^{\text{top}}_{\SL{4}{}{\permiss}}$ vanishes at $0$ and its zeroes have real part between $0$ and $14$ (see  \cite[Question~7.4, Question~7.5]{ros2015topological}).
\subsubsection{Non-existence of shadow-preserving lifts}
 Let $\alggroup$ be a smooth closed subgroup scheme of the $\Z$-group scheme $\GL{n}{}{}$ for some $n\in\N$. And let $\completion_\level=\completion/\primeideal^\level$. We recall the following definition from \cite{zor2016adjoint}:
\begin{defn}
\label{shadow}
Let $\level\in\N$ and $\element\in\lie(\alggroup)(\completion_{\level})$. We define the (\emph{group}) \emph{shadow} 
\[\sh{\alggroup(\completion_{\level})}{a}\allowbreak\leq\alggroup(\qfield)\] 
of $a$ to be the reduction $\bmod$ $\primeideal$ of the group stabilizer of $\element$ for the adjoint action of $\alggroup(\completion_{\level})$ on $\lie(\alggroup)(\completion_{\level})$. Analogously, the \emph{Lie shadow} 
\[\sh{\lie(\alggroup)(\completion_{\level})}{\element}\leq\lie(\alggroup)(\qfield)\]
of $a$ is the reduction $\bmod$ $\primeideal$ of the centralizer of $\element$ in $\lie(\alggroup)(\completion_{\level})$.
\end{defn}
\begin{defn}
\label{def:hereditary_shadow}
Let $\level\in\N$. We say that $b\in\lie(\alggroup)(\completion_{\level + 1})$ is a \emph{shadow-preserving lift} of $a$ when $b\equiv a\,\bmod\,\primeideal^\level$ and $\sh{\alggroup(\completion_{\level + 1})}{b}=\sh{\alggroup(\completion_{\level})}{a}$. Accordingly we say that $b\in\lie(\alggroup)(\completion_{\level + 1})$ is a \emph{Lie shadow-preserving lift} of $a$ when $b\equiv a\,\bmod\,\primeideal^\level$ and $\sh{\lie(\alggroup)(\completion_{\level + 1})}{b}=\sh{\lie(\alggroup)(\completion_{\level})}{a}$. If every element  of $\alg{\otherlevel}$ has a shadow-preserving lift for all $\otherlevel\in\N$, we say that $\alg{}$ is \emph{shadow-preserving}.
\end{defn}
\begin{main}
\label{main:G}
Let $\completion$ have odd residue field characteristic. The Lie lattice $\spl{4}{(\completion)}{}$ is not shadow-preserving.
\end{main}
This proves that, contrary to what was believed (see \cite[Remark~6.5]{akov2}), even if the residue field characteristic is allowed to be arbitrary big, being shadow-preserving does not follow from  $\spl{4}{(\C)}{}$ having smooth loci of constant centralizer dimension.
 \subsection{Background and motivation}
In order to contextualize the main results of this paper and provide motivation for them we now give a brief summary of some recent results and open problems concerning representation zeta functions. We start with an overview of some classes of groups with PRG whose representation zeta functions have been traditionally studied.\par
The most relevant examples for the purposes of this paper all arise from semisimple algebraic groups defined over a number field. Let $\otheralggroup$ be a connected, simply connected semisimple algebraic group defined over a number field $\numfield$.\par
In this context, the first example of a group with PRG is $\otheralggroup(\C)$: its representation zeta function has been studied first by Witten in \cite{wit1991wittenzeta} and its abscissa of convergence has been computed by Larsen and Lubotzky in \cite[Theorem~5.1]{larlub2008reprgrowth}.\par
Another class of groups with PRG consists of the groups $\otheralggroup(\ringint_v)$ where $\ringint_v$ is the completion of the ring of integers of $\numfield$ at a non-archimedean place $v$. These groups are known to be rigid (see \cite[Proposition~2]{baslubmag2002proalgebraic}) and the rationality results in  \cite{zap2006zeta} entail that their representation zeta function has rational abscissa of convergence.\par
 %
 %
Arithmetic groups with the congruence subgroup property (CSP) are yet another class of groups with PRG and, in a way, bring together the previous two examples. By an \emph{arithmetic group} we mean a group $\Gamma$ which is commensurable  to $\otheralggroup(\ringint_S)$, where $\ringint_S$ are the $S$-integers in $\numfield$ for a finite set $S$ of places of $\numfield$ including all the archimedean ones. %
In \cite{lubmar2004PRG} it was shown by Lubotzky and Martin that $\Gamma$ has PRG if and only if it has the CSP. %
Moreover, by a result of Larsen and Lubotzky  (see \cite[Proposition~1.3]{larlub2008reprgrowth}), when $\Gamma$ has the CSP, the representation zeta function of $\Gamma$ has an Euler product decomposition. To exemplify what this means, let for convenience $\Gamma=\otheralggroup(\ringint_S)$. Then, if $\Gamma$ has the strong CSP (i.e.\ the congruence kernel is trivial), its representation zeta function decomposes as
\begin{equation*}
\zetafunc{\Gamma}={\zetafunc{\otheralggroup(\C){}}}^{\lvert \numfield\,:\,\Q\rvert}\cdot\prod_{v\notin S} \zetafunc{\otheralggroup(\ringint_v)}
\end{equation*}
(also see  \cite[Equation (1.1)]{akov2013representation} and references therein).\par
Larsen and Lubotzky have proved in \cite[Theorem~8.1]{larlub2008reprgrowth} that the factors indexed by $v\notin S$  have abscissa of convergence bounded away from $0$ and proved that $1/15$ is a  lower bound for it. More recently Aizenbud and Avni have established upper bounds too. For instance in case $\otheralggroup = \SL{n}{}{}$ for some $n\in\N$, $22$ is an upper bound for the abscissa of convergence of $\SL{n}{(\ringint_v)}{}$ for almost all $v\notin S$ (see\cite{aizavn2016rational}, see also \cite{budzor2017representation} for an improved bound).\par
The fundamental consequence of the Euler product decomposition is that the so called \emph{global abscissa of convergence} $\abs{\Gamma}$ might be worked out if enough knowledge on the factors of the product is available. In general, Avni proves in \cite{avn2011rationalabs} that arithmetic groups with CSP have rational global abscissa of convergence. Even more strikingly, Larsen and Lubotzky conjectured that the global abscissa of convergence should depend solely on $\otheralggroup$ for higher-rank semisimple groups:
\begin{conj}[{Larsen and Lubotzky \cite[Conjecture 1.5]{larlub2008reprgrowth}}]
Let $H$ be a higher-rank semisimple group. Then, for any two irreducible lattices $\Gamma_1$ and $\Gamma_2$ in $H$, $\abs{\Gamma_1}=\abs{\Gamma_2}$.
\end{conj}%
%
We now give an overview of the state of the art in the solution of the open problem above. Fix a non-archimedean place $v$. In \cite[Theorem~1.2]{akov2011arithmetic} Avni, Klopsch, Onn and Voll prove a variant of Larsen and Lubotzky conjecture for higher-rank semisimple groups in characteristic $0$ assuming that both $\abs{\Gamma_1}$ and $\abs{\Gamma_2}$ are finite. In \cite{akov2013representation}, the same authors introduce the use of Igusa-type zeta functions in the study of representation growth. In particular they relate their representation zeta function to a generalized Igusa zeta function of the type described in \cite{voll2010functional}. In doing so, they prove functional equations for the representation zeta functions of almost all  groups in certain distinguished infinite families, e.g.\ the principal congruence subgroups $\SL{n}{(\ringint_v)}{\permiss}$ ($n\in\N$) for almost all $\permiss\in\N$ (see \cite[Theorem~A]{akov2013representation}). Using their $p$-adic integration technique, they compute explicit formulae for the representation zeta function of almost all of the principal congruence subgroups of $\SL{3}{(\ringint_v)}{}$ and $\mathrm{SU}_3(\mathfrak{D},\ringint_v)$, where $\mathfrak{D}$ is an unramified quadratic extension of $\ringint_v$. Using Clifford theory and the Euler product decomposition they are then able to deduce the abscissae of convergence of arithmetic groups of type $A_2$, thus establishing Larsen and Lubotzky's conjecture for groups of type $A_2$. %
The same authors in \cite{akov2} classify the similarity classes of $3\times 3$ matrices in $\gl{3}{(\ringint_v)}$ and $\mathfrak{gu}_3(\ringint_v)$ and obtain again the explicit formulae in \cite{akov2013representation} avoiding integrals. Using again Clifford theory, they also deduce explicit formulae for the representation zeta functions of $\SL{3}{(\ringint_v)}{}$ and of $\mathrm{SU}_3(\ringint_v)$. 
The application to the principal congruence subgroups of $\SL{4}{(\ringint_v)}{}$ (cf.\ \cref{thmB}) in the present work, starts for type $A_3$ the same line of investigation followed in
\cite{akov2013representation,akov2} for arithmetic groups of type $A_2$.
\subsection{Techniques and main technical results}
We shall now describe the techniques that allow us to directly compute the Poincar\'e series of certain Lie lattices and ultimately underpin the first three main results of this paper.  In this section we introduce these fundamental ideas and other complementary notions. At the end of the section we state the main technical result of this paper. It is worth remarking here again that this technical core does not require the Lie lattice to be quadratic.
\subsubsection{Poincar\'e series}
\label{sec:tech_main}
We start by recalling the definition of the commutator matrix associated with an $\completion$-basis of $\alginf$ and some relevant facts related to it, more details are found in \cite[Section~3.1]{akov2013representation}.\par
Let  $\mathcal{H} = \lbrace b_1,\dots,b_\dimalg\rbrace$ be an $\completion$-basis for $\alginf$. For $b_i,b_j\in \mathcal{H}$, there are $\lambda_{i,j}^1,\dots,\lambda_{i,j}^d\allowbreak\in\completion$ such that
\[[b_i,b_j]=\sum_{h=1}^d \lambda_{i,j}^h b_h.\]
The coefficients $\lambda_{i,j}^h$ for $i,j,h=1,\dots,d$ are called the \emph{structure constants} of $\alginf$ with respect to $\mathcal{H}$ (or structure constants of $\mathcal{H}$ for short). We define the \emph{commutator matrix} of $\alginf$ with respect to $\mathcal{H}$ (or the commutator matrix of $\mathcal{H}$ for short) as
	\begin{equation}
		\label{eq:comm_matrix}
		\Cmatrix{\mathcal{H}}{Y_1,\dots,Y_d}=\left( \sum_{h=1}^d \lambda_{i,j}^h Y_h\right)_{i,j}\in \Mat{d}{(\completion[Y_1,\dots,Y_d])}.
	\end{equation}
Let $\pi$ be a uniformizer of $\completion$. Fix an $\completion$-basis $\basis$. Let for convenience of notation $\cmatrix = \cmatrix_\basis$. Let now $\level\in\N$ and $\overline{\tuple{w}}\in (\completion/\primeideal^\level)^\dimalg$. Let $\tuple{w}\in \completion^\dimalg$ be a lift of $\overline{\tuple{w}}$. The matrix $\Cmatrix{}{\tuple{w}}$ is an antisymmetric $\dimalg\times\dimalg$ matrix, therefore its elementary divisors may be arranged in $\halfdim=\lfloor \dimalg/2\rfloor$ pairs $(\pi^{a_1},\pi^{a_1}),\dots,(\pi^{a_\halfdim},\pi^{a_\halfdim})$ for $0\leq a_1\leq\dots\leq a_\halfdim\in(\N_0\cup\lbrace \infty\rbrace)$ together with a single extra divisor $\pi^\infty= 0$ if $\dimalg$ is odd. We define
	\begin{align*}
		\nu_{\cmatrix, \level}(\overline{\tuple{w}})			&=(\min \lbrace a_i, \level \rbrace)_{i \in \lbrace 1,\dots, \halfdim \rbrace}.
	\end{align*}
It is easy to see that this definition does not depend on the choice of $\tuple{w}$.
\begin{defn}
\label{eq:N_sets}
Let
\begin{align*}
\cone{}		&=\completion^\dimalg \smallsetminus \primeideal\, \completion^\dimalg 	&&\\
\cone{\level}	&= (\completion/\primeideal^{\level})^\dimalg \smallsetminus \primeideal\, (\completion/\primeideal^{\level})^\dimalg	&& \level\in\N.
\end{align*}
Let $I=\lbrace i_1,\dots,i_\lastindex \rbrace_{<}$ be a (possibly empty) subset of $[\halfdim-1]_0 = \lbrace 0,\dots \halfdim - 1 \rbrace$ such that $i_1<i_2,\dots<i_\lastindex$. We set $i_0=0$ and $i_{\lastindex+1}=\halfdim$ and we write 
	\begin{align*}
		\mu_j		&=i_{j+1}-i_j 							&&\text{for } j\in \lbrace 0,\dots, \lastindex \rbrace; &
		N			&=\sum_{j=1}^\lastindex r_j 				&&\text{for } \mathbf{r}_I=(r_{1},\dots,r_{\lastindex})\in \N^{\lvert I\rvert}.
	\end{align*}
The \emph{Poincar\'e series of $\cmatrix$} is
	\begin{equation*}
		\pseries{\mathcal{R}}{(s)}=\sum_{\stackrel{I\subseteq [\halfdim-1]_0}{I=\lbrace i_1,\dots,i_\lastindex \rbrace_{<}}}\sum_{\mathbf{r}_I\in\N^{\lvert I\rvert}} \lvert \pcoeff{I}{(\cmatrix)}\rvert\,q^{-s\sum_{j=1}^\lastindex r_j(\halfdim-i_j)},
	\end{equation*}
where
	\begin{multline*}
		\pcoeff{I}{(\cmatrix)}= \lbrace \tuple{w}\in \cone{N}\mid \nu_{\cmatrix, N}(\tuple{w})=(\underbrace{0,\dots,0}_{\mu_\lastindex},\underbrace{r_{\lastindex},\dots,r_{\lastindex}}_{\mu_{\lastindex-1}},\\
			\underbrace{r_{\lastindex} + r_{\lastindex-1} ,\dots,r_{\lastindex} + r_{\lastindex-1} }_{\mu_{\lastindex-2}}%
					\dots%
						,\underbrace{N,\dots,N}_{\mu_0})\in \N_0^\halfdim\rbrace.
	\end{multline*}
If $\basis'$ is another basis for $\alginf$, it is known that $\pseries{\mathcal{R}}{(s)} = \pseries{\mathcal{R}_{\basis'}}{(s)}$, we may therefore define the {\em Poincar\'e series of $\alginf$} as
\[
\mathcal{P}_{\alginf}{(s)} = \pseries{\mathcal{R}}{(s)}.
\]
\end{defn}
\begin{rem}
\label{rem:from_N_to_reps}
The link between the Poincar\'e series and the representation zeta function of a potent and saturable FAb $p$-adic analytic group is given via the Lazard correspondence and the Kirillov orbit method. The Lazard correspondence associates with such a group $G$ a potent and saturable $\Z_p$-Lie lattice $\alginf$. There is a natural $G$-action on $\alginf^{\vee} = \Hom_{\Z}(\alginf, \C^{\times})$ given by the dual of the $G$-adjoint action on $\alginf$. The Kirillov orbit method says that the orbits of this action correspond to continuous irreducible representation of $G$ and their cardinalities are equal to the square of the degree of the representation they correspond to. In practice then the sets defining the Poincar\'e series are collecting (in coordinates) elements of $\Hom(\alginf,\completion)$ which represent elements of $\alginf^{\vee}$ whose stabilizers for the $G$-co-adjoint action have the same index; or, equivalently, whose $G$-orbits have the same size (see \cite[Section~2.2]{akov2013representation}).
\end{rem}
\subsubsection{Isolated submodules and geometrical smoothness}
We now introduce the tools we shall use to study the geometry of the sets defining the Poincar\'e series, ultimately allowing us to express it as a sum of products of geometric progressions. 
\begin{defn}
\label{def:rk_locus}
For $i \leq \dimalg$, let $\pfaffians{i}{}\subseteq\completion[X_1,\dots, X_\dimalg]$ be the ideal generated by the $i\times i$ principal minors of $\cmatrix$ and let $\pfaffians{\dimalg + 1}{} = (0)$. We write
\[
\rankvar{\cmatrix}{i}=\spm \left(\completion[X_1,\dots, X_\dimalg]/\pfaffians{i}{}\right).
\]
If the $\dimalg$-dimensional affine space over $\completion$ is defined as 
\[
\A^\dimalg = \spm\left(\completion[X_1,\dots, X_\dimalg]\right),
\]
then the canonical projection
\[
\completion[X_1,\dots, X_\dimalg]\rightarrow \completion[X_1,\dots, X_\dimalg]/\pfaffians{i}{}
\]
embeds $\rankvar{\cmatrix}{i}$ into $\A^{\dimalg}$, in what follows we shall always see $\rankvar{\cmatrix}{i}$ as an $\completion$-subscheme of $\A^{\dimalg}$.
For $2\mu\leq \dimalg$, the \emph{rank-$2\mu$ locus} $\ranklocus{\cmatrix}{2\mu}{}$ of $\cmatrix$ is the $\completion$-scheme defined as the scheme-theoretic complement of $\rankvar{\cmatrix}{2\mu}$ when seen as a closed subscheme of $\rankvar{\cmatrix}{2\mu + 1}$. If $\ring$ is an $\completion$-algebra, we write
\[
 \rankvar{\cmatrix}{2\mu}(\ring) = \Hom_\completion\left(\completion[X_1,\dots, X_\dimalg]/\pfaffians{2\mu}{},\ring\right)
\]
for the set of \emph{$\ring$-points} of $\rankvar{\cmatrix}{2\mu}$. By \cite[Theorem~1.2]{wat1979introduction} there is a natural correspondence between the set defined above and the set of solutions of $\pfaffians{2\mu}$ in $\ring^\dimalg$. For convenience, in what follows we shall identify these two sets. Accordingly we write
\[
 \ranklocus{\cmatrix}{2\mu}{(\ring)} = \rankvar{\cmatrix}{2\mu + 1}(\ring) \smallsetminus \rankvar{\cmatrix}{2\mu}{(\ring)}.
\]
\end{defn}
\begin{rem}
Before proceeding, a comment on the last definition: it might seem a little confusing that we have defined only loci consisting of points that give even rank when plugged into the commutator matrix. This is explained because antisymmetric matrices always have even rank; and so, for all $2\mu\leq \dimalg$, the zeroes of the $2\mu\times 2\mu$ principal minors coincide with the zeroes of the $(2\mu - 1)\times (2\mu - 1)$ principal minors. It follows that even extending the definition above to odd ranks one would not find new algebraic sets besides the loci of even rank.
\end{rem}
In order to compute the Poincar\'e series of $\cmatrix$ we study the set of $\completion$-points of $\ranklocus{\cmatrix}{2\mu}{}$ which by the above definition and conventions is 
\[
\ranklocus{\cmatrix}{2\mu}{(\completion)} = \ranklocus{\cmatrix}{2\mu}{(\localfrac)}\cap\,\completion^\dimalg.
\]
We recall the definition of an isolated submodule:
\begin{defn}
\label{def:isol}
Let $M$ be an $\completion$-module and $N$ be a submodule of $M$. We say that $N$ is \emph{isolated} in $M$ if $M/N$ is a torsion-free $\completion$-module. The smallest isolated $L$ containing $N$ is called the \emph{isolator} of $N$ and denoted by $\isol(N)$.
\end{defn}
Fixing a basis for $\alginf$ gives $\completion^\dimalg$ the structure of an $\completion$-Lie lattice and $\qfield^\dimalg$ the structure of an $\qfield$-Lie algebra. Here and throughout this work, the Lie bracket of two elements of $\completion^\dimalg$, or $\qfield^\dimalg$ respectively, is intended to be the one associated with the Lie structure given by $\basis$. We introduce here the following pivotal notion for this paper.
\begin{defn}
\label{def:geo_smooth_rk}
We say that $\ranklocus{\cmatrix}{2\mu}{(\completion)}$ is \emph{geometrically smooth} when 
\begin{enumerate}
\item every $\overline{\tuple{x}}\in \ranklocus{\cmatrix}{2\mu}{(\qfield)}$ has a lift $\tuple{x}\in \ranklocus{\cmatrix}{2\mu}{(\completion)}$ (i.e.\ there exists $\tuple{x}\in \ranklocus{\cmatrix}{2\mu}{(\completion)}$ such that $\tuple{x}\equiv \overline{\tuple{x}}\,\bmod\,\primeideal$).
\item For every $\tuple{x}$ of $\ranklocus{\cmatrix}{2\mu}{(\completion)}$ such that $\tuple{x}\in \ranklocus{\cmatrix}{2\mu}{(\qfield)}$ $\bmod\,\primeideal$ the following holds:
\[
[\ker \cmatrix(\tuple{x}),\ker \cmatrix(\tuple{x})]
\]
is an isolated $\completion$-submodule of $\completion^\dimalg$.
\end{enumerate}
 When $\ranklocus{\cmatrix}{2\mu}{}$ is geometrically smooth for all $2\mu\leq\dimalg$, we say that $\cmatrix$ has \emph{geometrically smooth rank-loci}.
\end{defn}
\begin{defn}
\label{def:g_geo_smooth_rk}
It is known that if another basis $\basis'$ of $\alginf$ is chosen, then the rank-loci of $\cmatrix_{\basis'}$ are isomorphic to the rank-loci of $\cmatrix$ as $\completion$-schemes (this is a straightforward consequence of \cref{lem:basis_change}). This means that having geometrically smooth rank-loci is actually a property of $\alginf$ that does not depend on the choice of the basis. So we say that $\alginf$ has {\em geometrically smooth rank-loci} when there is a basis $\otherbasis$ such that $\cmatrix_{\otherbasis}$ has geometrically smooth rank-loci.\par
\end{defn}
 \Cref{sec:rank_loci}, gives the motivation for the name ``geometrical smoothness'': indeed, we show that if the rank-loci are smooth $\completion$-schemes and their sets of $\qfield$-points are also smooth, then the geometrical smoothness condition is automatically satisfied. Although geometrical smoothness is easier to check as our computations will show, we do not know whether it is really more general than smoothness in the context of rank-loci of commutator matrices. We leave this as an interesting open question.
\subsubsection{Kernel classes and classifications}
In order to state the main technical result we still need to break up $\qfield^\dimalg$ into pieces all whose elements give a specified contribution to the Poincar\'e series. In the context of representation zeta functions of a potent and saturable FAb $p$-adic analytic group $G$, this gives us the flexibility to group together points whose Lie centralizers are $G$-conjugate (roughly replicating the approach of \cite{akov2}), or just isomorphic (thus obtaining a similar approach to the one in \cite{zor2016thesis}). 
We denote by $\overline{\cmatrix}$ the matrix of linear forms over $\qfield$ obtained reducing $\bmod\,\primeideal$ the entries of $\cmatrix$.
\begin{defn}
\label{def:kernel_class}
An $\cmatrix$-kernel class $\radiclass$ is a subset of $\qfield^\dimalg$ such that for any two $\tuple{x},\tuple{x'}\in\radiclass$
\begin{align*}
\dim_{\qfield} \ker \overline{\cmatrix}(\tuple{x}) 								&=  \dim_{\qfield} \ker \overline{\cmatrix}(\tuple{x}')\\
\dim_{\qfield} [\ker \overline{\cmatrix}(\tuple{x}), \ker \overline{\cmatrix}(\tuple{x})]	&=  \dim_{\qfield} [\ker \overline{\cmatrix}(\tuple{x}'), \ker \overline{\cmatrix}(\tuple{x}')].
\end{align*}
\end{defn}
Examples of $\cmatrix$-kernel classes are subsets all of whose elements give isomorphic kernels when plugged into $\overline{\cmatrix}$. Other examples arise when $\alginf$ is the set of $\completion$-points of the Lie algebra of a $\Z$-group scheme $\alggroup$. In this setting examples of $\cmatrix$-kernel classes are subsets all of whose elements give kernels that, for instance, are conjugate for the $\alggroup(\qfield)$-action on $\qfield^\dimalg$ induced by the  $\alggroup(\completion)$-adjoint action on $\alginf$.
\begin{defn}
A set of disjoint $\cmatrix$-kernel classes covering $\qfield^\dimalg$ is called a \emph{classification} by $\cmatrix$-kernels. Members of a classification by $\cmatrix$-kernels $\arr$ are called kernel $\arr$-classes or kernel classes when no confusion is possible.
\end{defn}
Given an $\cmatrix$-kernel class $\radiclass$ we have well defined
\begin{align*}
\dimalg_{\radiclass} 	&= \dim_{\qfield} \radic		&&\radic = \ker \overline{\cmatrix}(\tuple{x})\text{ for any } \tuple{x}\in\radiclass\\
\dimalg'_{\radiclass}	&= \dim_{\qfield} [\radic,\radic]	&&\radic = \ker \overline{\cmatrix}(\tuple{x})\text{ for any } \tuple{x}\in\radiclass.
\end{align*}
Let $\arr$ be a classification by $\cmatrix$-kernels. A \emph{sequence} of $\arr$-kernel classes is a set $\lbrace \radiclass_1,\dots,\radiclass_\lastindex\rbrace$ ($\lastindex\in \N_0$) of kernel $\arr$-classes such that $\dimalg_{\radiclass_1}> \dimalg_{\radiclass_2}>\cdots>\dimalg_{\radiclass_t}$. 
\begin{defn}
\label{def:cal_F}
Let $\mathcal{S} =  \lbrace \radiclass_1,\dots,\radiclass_\lastindex\rbrace$ be a sequence of kernel $\arr$-classes
\begin{enumerate}
\item If $\mathcal{S} = \emptyset$ we define $\mathcal{F}_{\mathcal{S}}(\cmatrix) = \lbrace \emptyset\rbrace$.
\item If $\mathcal{S} \neq \emptyset$. We define $\mathcal{F}_{\mathcal{S}}(\cmatrix)$ as the set of $\lastindex$-tuples $(\tuple{x}_1,\dots,\tuple{x}_\lastindex)$ of elements of $\qfield^\dimalg$ such that for all $j = 1,\dots,\lastindex$ for $\tuple{y}_j = \sum_{k = j}^{\lastindex} \tuple{x}_k$
\begin{enumerate}
\item $\tuple{y}_j \in \radiclass_j$
\item	$\tuple{x}_{j-1} \in [\ker \overline{\cmatrix}\left(\tuple{y}_j\right), \ker \overline{\cmatrix}\left(\tuple{y}_j\right)]$ for $j > 1$.
\end{enumerate}
\end{enumerate}
\end{defn}
Set $\halfdim = \lfloor \dimalg / 2\rfloor$. Fix  an ordered subset $I = \lbrace i_1,\dots, i_\lastindex\rbrace_{<}$ of $[\halfdim - 1]_0 = \lbrace 0,\dots, \halfdim - 1\rbrace$.
\begin{defn}
An $I$-sequence of kernel $\arr$-classes is a sequence of kernel $\arr$-classes $\lbrace \radiclass_1,\dots,\radiclass_\lastindex\rbrace$
such that $\dimalg_{\radiclass_j} = \dimalg - 2(\halfdim - i_{\lastindex + 1 - j})$ for all $j = 1,\dots,\lastindex$. The set of all $I$-sequences of kernel $\arr$-classes is denoted by $\radseq_{I}^\arr(\cmatrix)$.
\end{defn}
\begin{main}
\label{main:C}
Let $\cmatrix$ be the commutator matrix of $\alginf$ with respect to an $\completion$-basis. Let $\arr$ be a classification by $\cmatrix$-kernels. Assume that $\cmatrix$ has geometrically  smooth rank-loci. Then the Poincar\'e series of $\cmatrix$ is 
\[
\sum_{\stackrel{I\subseteq [\halfdim-1]_0}{I=\lbrace i_1,\dots,i_\lastindex \rbrace_{<}}}\sum_{\stackrel{\mathcal{S}\in\radseq^\arr_{I}(\cmatrix)}{\mathcal{S} =  \lbrace \radiclass_1,\dots,\radiclass_\lastindex\rbrace}} \left\lvert \mathcal{F}_{\mathcal{S}}(\cmatrix)\right\rvert q^{-(\dimalg - \dimalg'_{\radiclass_\lastindex})}%
\prod_{\radiclass\in\mathcal{S}} \frac{q^{\dimalg - \dimalg'_{\radiclass} - s\frac{\dimalg - \dimalg_{\radiclass}}{2}}}{1 - q^{\dimalg - \dimalg'_{\radiclass} - s\frac{\dimalg - \dimalg_{\radiclass}}{2}}}.
\]
\end{main}
The way this result is proved is an evolution of the ideas of \cite{zor2016adjoint}. Indeed the methods therein are a way of computing the Poincar\'e series of $\alginf$ by looking at the fibres of the maps induced by the projection of $\completion$ onto its finite quotients. An idea that essentially was already present in \cite{akov2}. Here, as we shall see in \cref{lem:theta}, we introduce new maps that project the sets whose cardinalities is encoded in the Poincar\'e series onto suitably defined sets involving only properties of the Lie algebra $\algfin$. The main idea is to describe how the elementary divisors $\nu_{\cmatrix,\level}(\tuple{x})$ ($\tuple{x}\in(\completion/\primeideal^\level)^\dimalg$, $\level\in\N$) behave under lifting $\tuple{x}$ to $(\completion/\primeideal^{\level + 1})^\dimalg$. This is a coarser analog of the method of shadows in \cite{zor2016adjoint,akov2}, and substantially converts the counting problem involving the elementary divisors into an easier and coarser counting problem involving ranks of matrices over $\qfield$; thus translating an enumeration over $\completion^\dimalg$ to an enumeration over $\qfield^\dimalg$.
\subsection{Organization of the paper}
The paper is divided into three parts, the first one proves \cref{main:C}, the second one applies this theorem to $\spl{4}{(\completion)}{}$, and the third one proves that this $\completion$-Lie lattice is not shadow-preserving.\par
The first part is essentially a proof of \cref{main:C} and its consequences: \cref{main:D,main:abscissa_P,main:E}. First of all, that is in \cref{sec:rank_loci}, we examine the Jacobians of the polynomials defining the rank-loci of a commutator matrix. We prove that smoothness implies geometrical smoothness of the rank-loci and give a criterion for establishing geometrical smoothness in case the rank-loci over $\C$ are smooth and irreducible (cf.\ \cref{prop:geo_sm_crit}). This will be useful in \cref{part:sl4} to determine for which primes \cref{main:D} applies to $\spl{4}{(\completion)}{}$. In \cref{sec:smooth_loci} we introduce the techniques allowing us to compute Poincar\'e series of commutator matrices with geometrically smooth rank-loci. After having proved \cref{main:C} and its corollary \cref{main:abscissa_P}, we specialize them for quadratic Lie lattices obtaining \cref{main:D,main:E}.\par
The second part applies \cref{main:D} to the Lie lattices $\spl{4}{(\completion)}{}$ when $\completion$ has odd residue field characteristic. First of all we check that the hypotheses of \cref{main:D} are satisfied. In \cref{sec:cross-sec} we elaborate an effective way of computing the cardinalities of the sets defining the Poincar\'e series in \cref{main:D}. This method presupposes arranging the elements of  $\spl{4}{(\qfield)}{}$ according to their centralizer-dimension and dimension of related derived subalgebra, which we do in \cref{sec:centr}. Finally the computation is  carried out in \cref{sec:sl4}. We decided not to overcrowd this section with details which, albeit with a different notation, are available in \cite[Chapter~4]{zor2016thesis}.\par
The last part proves that $\spl{4}{(\completion)}{}$ is not shadow-preserving, we have preferred to separate this part from the bulk of the other computations for $\spl{4}{(\completion)}{}$ as they no longer refer to shadows.
\begin{table}[h]
\caption{Some frequently used notation.}
\label{tab:notation}
\begin{tabular}{ll}
\toprule
Symbol					& Meaning\\
\midrule
$\N$						& $\lbrace 1, 2, 3,\dots\rbrace$\\
$\N_0$						& natural numbers, i.e.\ $\N\cup \lbrace 0\rbrace$\\
$[n]$						& $\lbrace 1,\dots,n\rbrace\subseteq \N$ ($n\in\N$)\\
$[n]_0$					& $[n] \cup \lbrace 0\rbrace$ ($n\in\N$)\\
$\lbrace i_1,\dots,i_\lastindex\rbrace_{<}\subseteq\N$		& ordered subset of $\N$: $I_1,\dots,i_\lastindex\N$ and $i_1<i_2<\cdots < i_\lastindex$\\
\midrule
$\Mat{m \times n}{(\ring)}$		&$m \times n$ ($m,n\in\N$) matrices with entries in the ring $\ring$\\
$\mathrm{Asym}_{n}(\ring)$	&$n \times n$ ($n\in\N$) antisymmetric matrices with entries \\
						&in the ring $\ring$\\
\midrule
$\Q_p$					& $p$-adic numbers for $p$ a prime integer\\
$\Z_p$					& $p$-adic integers for $p$ a prime integer\\
\midrule
$\completion$				& compact DVR of characteristic $0$\\
$\completion^n$			& $\completion\times\cdots\times\completion$ $n$-times ($n\in\N$)\\
$\pi$						& uniformizer of $\completion$\\
$\primeideal$				& prime ideal of $\completion$\\
$\primeideal^\level$			& $\primeideal\cdots\primeideal$ $\level$-times ($\level\in\N$)\\
$\primeideal^{(n)}$			& $\primeideal\times\cdots\times\primeideal$ $n$-times ($n\in\N$)\\
$\completion^\times$		& $\completion\smallsetminus\primeideal$\\
$\qfield$					& residue field of $\completion$\\
$\completion_\level$			& $\completion/\primeideal^\level$\\
$\localfrac$				& field of fractions of $\completion$\\
\midrule
$M^*$					& $M\smallsetminus \primeideal M$ for a non-trivial $\completion$-module $M$\\
$\lbrace 0\rbrace^*$			& $\lbrace 0\rbrace$\\
\midrule
$\alginf$					&$\completion$-Lie lattice\\
$\algfin$					& reduction of $\alginf$ modulo $\primeideal$\\
$\alg{\level}$				& reduction of $\alginf$ modulo $\primeideal^\level$ ($\level\in\N$)\\
$\dimalg$					& $\rk_\completion \alginf$\\
$\halfdim$					& $\lfloor \dimalg/2\rfloor$\\
\midrule
$\centr{\mathfrak{h}}{x}$		& Lie centralizer of $x$ in the Lie ring $\mathfrak{h}$ over a ring $\ring$\\
$[\mathfrak{h}, \mathfrak{h}]$	& $\linspan_\ring([x,y]\mid x,y\in\mathfrak{h})$, i.e.\ the derived subring  \\
						& of the $\ring$-Lie lattice $\mathfrak{h}$\\
$\mathcal{H}^\sharp$		& dual basis of the $\ring$-basis $\mathcal{H}$\\			
\midrule
$\A^n$					& $n$-dimensional affine space over $\completion$, i.e.\ $\spm(\completion[X_1,\dots,X_n])$\\
\bottomrule	
\end{tabular}		
\end{table}
\subsection{Notation}
We denote by $\N$ the set of the positive integers $\lbrace 1,2,\dots\rbrace$, while $\N_0=\lbrace 0,1,2,\dots\rbrace$ are the natural numbers. Analogously, for $n\in\N$ we set $[n]=\lbrace 1,\dots,n\rbrace$ and $[n]_0=\lbrace 0,\dots,n\rbrace$. Throughout this work $\completion$ is a compact discrete valuation ring of characteristic $0$, prime ideal $\primeideal$, uniformizer $\pi$ and residue field $\qfield$ of characteristic $p$ and  cardinality $q$. The field of fractions of $\completion$ is denoted by $\localfrac$. The field of $p$-adic numbers is denoted by $\Q_p$ and the ring of $p$-adic integers by $\Z_p$. \par
As conventional, the group of units of a ring $\ring$ is $\ring^\times$. We introduce a similar notation for non-trivial $\completion$-modules as follows. Given such a module $M$, we write $M^*=M\smallsetminus \primeideal M$. For the trivial $\completion$-module we set $\lbrace 0\rbrace^*=\lbrace 0\rbrace$. Moreover, defining $M_\ell$ to be the reduction modulo $\primeideal^\ell$ of $M$, we write
\begin{align*}
\reducinf{\level}:			& M \rightarrow M_\level			&& \level\in\N\\
\reduc{\level}{\otherlevel}:		& M_\level\rightarrow M_\otherlevel 	&& \level > \otherlevel
\end{align*}
for the maps defined by reducing modulo $\primeideal^\level$ and $\primeideal^\otherlevel$ respectively. If $\level,\otherlevel\in\N$ with $\level > \otherlevel$ and $a\in M_\otherlevel$ for some $\otherlevel\in\N$, we say that $b\in\reduc{\level}{\otherlevel}^{-1}(a)$ is a \emph{lift} of 
$a$ to $M_\level$.\par
Let $\ring$ be a ring. We denote by $\Mat{m \times n}{(\ring)}$	the $m \times n$ ($m,n\in\N$) matrices with entries in $\ring$ and by $\mathrm{Asym}_{n}(\ring)$	the $n \times n$  antisymmetric matrices with entries in $\ring$.\par
An \emph{$\ring$-Lie ring} is an $\ring$-module endowed with a Lie bracket. Let $\mathfrak{h}$ be an $\ring$-Lie ring and $x\in\mathfrak{h}$. The \emph{Lie centralizer} of $x$ in $\mathfrak{h}$ is
\[
\centr{\mathfrak{h}}{x} = \lbrace y\in\mathfrak{h}\mid [x,y] = 0\rbrace.
\]
An \emph{$\ring$-Lie lattice} is an $\ring$-Lie ring that is free and finitely generated over $\ring$. For an $\ring$-Lie lattice $\mathfrak{h}$ we define $[\mathfrak{h},\mathfrak{h}]$ to be its derived subring. If $\mathcal{H}$ is an $\ring$-basis for $\mathfrak{h}$, its dual will be denoted by $\mathcal{H}^\sharp$.\par
When $\otheralggroup$ is a  closed smooth $\Z$-subgroup scheme of $\GL{n}{}$ ($n\in\N$) and $\mathfrak{h} = \lie(\otheralggroup)(\ring)$, we write 
\[
\centr{\otheralggroup(\ring)}{x}
\]
for the stabilizer of $x$ for the adjoint action of $\otheralggroup(\ring)$ on $\mathfrak{h}$. If $\mathfrak{h}$ is stable under the adjoint action of $\GL{n}{(\ring)}$ on $\gl{n}{(\ring)}{}$ we write
\[
\centr{\GL{n}{(\ring)}}{x}
\]
for the stabilizer of $x$ for the action induced on $\mathfrak{h}$.\par
We write $R\llbracket X_1,\dots, X_t \rrbracket$ for the ring of formal power series in the $t$ variables $X_1,\dots, X_t $.\par
\Cref{tab:notation} gives an overview of the notation introduced here and also reports some frequently used symbols.
\subsection{Acknowledgements }
I am indebted to Christopher Voll and Benjamin Martin for their precious advice. I also wish to thank Tobias Rossmann, Giovanna Carnovale, Andrea Lucchini, Uri Onn, Shai Shechter, Alexander Stasinski and Wim Veys for the interesting conversations and insightful comments on this work. I wish to thank Tobias Rossmann for introducing me to the SageMath system and helping me set up the computer computations described in \cref{ex:non_shadow}.\par
\Cref{part:sl4,part:hered} of this work are part of my PhD thesis. I acknowledge financial support from the Faculty of Mathematics of the University of Bielefeld, CRC 701 during my doctoral studies. I am currently supported by the Research Project G.0939.13N of the Research Foundation - Flanders (FWO)

\part{Poincar\'e series}
\section{Rank-loci}
\label{sec:rank_loci}
In this section we investigate the geometry of the $\completion$-schemes defined in \cref{def:rk_locus}. The core is \cref{prop:rk_jacobian} which describes the image of the Jacobian matrix at a point of a rank-locus. Fix $2\mu\leq \dimalg$. We start by noticing that if a different basis $\basis'$ is used to define the commutator matrix of $\alginf$, the basis change defines an $\completion$-scheme isomorphism $\A^\dimalg\rightarrow \A^\dimalg$ that carries the rank-$2\mu$ locus of $\cmatrix$ to the rank-$2\mu$ locus of $\cmatrix_{\basis'}$.
\begin{lem}
\label{lem:basis_change}
Let $\basis'$ be another $\completion$-basis for $\alginf$, and let $S$ be the basis-change matrix from $\basis$ to $\basis'$. Then, for all $\tuple{v}\in\completion^{\dimalg}$,
\[\transpose{S}\cmatrix'(\tuple{v})S=\cmatrix(\tuple{v}S^{-\mathrm{t}}),
\]
where $\cmatrix'$ the commutator matrix of $\alginf$ with respect to $\basis'$.
\begin{proof}
Let $\tuple{v}=(v_1,\dots,v_\dimalg)\in\completion^{\dimalg}$. Let also ${\basis'}^\sharp=\lbrace {b'_1}^\sharp,\dots,{b'_\dimalg}^\sharp\rbrace$. The matrix $\cmatrix'(\tuple{v})$ is the matrix of the bilinear form $b_\omega$ on $\alginf$ defined by 
	\[	b_\omega(x,y) = \omega([x,y]),\]
where $\omega=\sum_{i=1}^{\dimalg}v_i {b'_i}^\sharp$. Since $S$ is the basis change from $\basis$ to $\basis'$, $\tuple{v}S^{-\mathrm{t}}$ expresses the coordinates of $\omega$ with respect to $\basis$. It follows that $\cmatrix(\tuple{v}S^{-\mathrm{t}})$ is the matrix of $b_\omega$ with respect to $\basis$. Hence the equality with $\transpose{S}\cmatrix'(\tuple{v})S$.
\end{proof}
\end{lem}
\begin{notation}
Let $\ring$ be a ring. Let $e\in\N$ and  $f\in\ring[X_1,\dots,X_e]$. The homomorphism of $\ring$-algebras 
\[
\xymatrix@R=3pt{\ring[t]\ar[r]&\ring[X_1,\dots,X_e]\\
t\ar@{|->}[r]&f}
\]
defines a morphism $\A^\dimalg\rightarrow\A^1$ of affine schemes over $\ring$. The preimage of $\A^1\smallsetminus\lbrace 0\rbrace$ through this morphism is an affine scheme that is open in the Zariski topology of $\A^\dimalg$. We denote it by $D(f)$. If $A$ is an $\ring$-algebra, then $D(f)(A)$ consists of the points $\tuple{a}\in\A^\dimalg(A)$ such that $f(\tuple{a})\in A^{\times}$.
\end{notation} 
\begin{prop}
\label{prop:rk_jacobian}
Let $\ring$ be an $\completion$-algebra. Let $\mathcal{J}_{2\mu}$ be the Jacobian matrix associated with a set of generators of $\pfaffians{2\mu + 1}$. Then 
\[
\im \mathcal{J}_{2\mu}(\tuple{x}) = [\ker \cmatrix(\tuple{x}), \ker \cmatrix(\tuple{x})].
\]
\begin{proof}
Assume for convenience of notation that $\ker\cmatrix(\tuple{x})$ is spanned by the last $\dimalg - 2\mu$ elements of the standard basis of $\ring^\dimalg$ (if not we may always change the basis by \cref{lem:basis_change}). Let $R_{ij}$ denote the $(i,j)$-th entry of $\cmatrix$ for $i,j\in\lbrace 1,\dots,\dimalg\rbrace$. Let also 
\begin{align*}
\mathcal{M}	&= ( R_{ij} )_{i,j\in\lbrace 1,\dots 2\mu_\lastindex\rbrace}
\end{align*}
Define $g = \det \mathcal{M}$. The set $D(g)= \spm(\completion[X_1,\dots, X_\dimalg]_g)$ is an open set in the Zariski topology of $\A^\dimalg$ and $\tuple{x}\in D(g)(\ring)$. %
Set, for convenience of notation $A = \completion[X_1,\dots, X_\dimalg]$. We embed $A$ into $A_g$ by sending each $f\in A$ to $f/g^0\in A_g$. Accordingly we define $J_\mu$ to be the ideal generated by the image of $\pfaffians{2\mu + 1}$ under this embedding. By a theorem of Kronecker, since $g$ is an invertible $2\mu\times 2\mu$ minor, $J_\mu$ may be generated by the $(2\mu + 1)\times(2\mu + 1)$ bordered minors of $\mathcal{M}$; i.e.\ by%
\begin{align*}
F_{ij} &= \det \begin{pmatrix}
		\mathcal{M}									&\begin{matrix}
																R_{1j}\\
																\vdots\\
																R_{2\mu\, j}
															\end{matrix}\\
		\begin{matrix} R_{i\,1}&\dots&R_{i\,2\mu} \end{matrix}		& R_{ij}
		   \end{pmatrix}																	&& (i,j > 2\mu).
\end{align*}
We shall now compute the Jacobian matrix of this presentation of $J_\mu$. Let 
\[
I_T 	= 	\left (\lbrace R_{ij}\mid i\leq 2\mu \; j > 2\mu \rbrace \right ).
\]
For all $i,j$ as above, we compute $F_{ij}$ by means of a Laplace expansion with respect to the last column and find that
\[
F_{ij} = g\cdot R_{ij} + G_{ij}
\]
with $G_{ij}\in I_T$.  Moreover, deriving using Leibniz's rule, one finds that for all $k\in\lbrace 1,\dots,\dimalg\rbrace$
\[
\frac{\partial F_{ij}}{\partial X_k} = g \cdot \frac{\partial R_{ij}}{\partial X_k} + G_{ij}' 
\] 
with $G_{ij}'\in I_T$. Since the last $\dimalg - 2\mu$ elements of $\basis$ span $\ker \cmatrix(\tuple{x})$, we have that $\tuple{x}$ is a common zero of the polynomials in $I_T$, so the last equality implies that  the Jacobian matrix associated with the $F_{ij}$'s has image
\[
\linspan\left(\left\lbrace\left. \begin{pmatrix}
						\frac{\partial R_{ij}}{\partial X_1}\\
						\vdots\\
						\frac{\partial R_{ij}}{\partial X_\dimalg}
						\end{pmatrix} \,\right\vert\, ij > 2\mu_\lastindex\right\rbrace\right) = [\ker \cmatrix(\tuple{x}), \ker \cmatrix(\tuple{x})].
\]
This suffices to conclude as the image of the Jacobian matrix does not change by passing to a localization or by changing the presentation of the ideal.
\end{proof}
\end{prop}
\begin{rem}
Let $\level\in\N$. If $\tuple{x}$ belongs to  $\ranklocus{\cmatrix}{2\mu}{(\completion)}$ modulo $\primeideal^\level$  and has a lift $\widehat{\tuple{x}}\in \ranklocus{\cmatrix}{2\mu}{(\completion)}$, then 
\[
\im \mathcal{J}_{2\mu}(\tuple{x}) \equiv [\ker \cmatrix(\widehat{\tuple{x}}), \ker \cmatrix(\widehat{\tuple{x}})]	\mod	\primeideal^\level.
\]
\end{rem}
With this description of the image of the Jacobian at hand, it is now easy explain how the notion of geometrical smoothness generalizes the notion of smoothness. Namely
\cref{prop:rk_jacobian} implies that if $\ranklocus{\cmatrix}{2\mu}{}$ is smooth over $\completion$ and its $\qfield$-points are also smooth, then $\ranklocus{\cmatrix}{2\mu}{(\completion)}$ is geometrically smooth. Indeed, the presence of lifts of $\ranklocus{\cmatrix}{2\mu}{(\qfield)}$ to $\ranklocus{\cmatrix}{2\mu}{(\completion)}$  is guaranteed by Hensel's lemma (cf.\ \cite[\S 4.6,~Corollary 3]{bourbaki1972commutative}). While, if $\tuple{x}\in\ranklocus{\cmatrix}{2\mu}{(\completion)}$, then there is an affine neighbourhood
\[
\mathcal{U} = D(g) \cong \spm (\completion[X_1,\dots,X_\dimalg]_g)
\]
for some $g\in\completion[X_1,\dots,X_\dimalg]$ such that $g(\tuple{x})\in\completion^\times$. Locally in $\mathcal{U}$, the rank-locus $\ranklocus{\cmatrix}{2\mu}{}$ is defined by polynomials with invertible Jacobian. In other words there is a minimal set of generators $f_1,\dots,f_c$ of $P_{2\mu + 1}$ such that 
\[
\mathcal{U}\cap \ranklocus{\cmatrix}{2\mu}{} =\spm\left( \completion[X_1,\dots,X_\dimalg]_g/(f_1,\dots,f_c)\right)
\]
and the Jacobian matrix
\[
\mathcal{J}_f^{c} = \begin{pmatrix}
\frac{\partial f_1}{\partial X_1}	&\dots	&\frac{\partial f_c}{\partial X_1}	\\
\vdots					&		&\vdots					\\
\frac{\partial f_1}{\partial X_c}	&\dots	&\frac{\partial f_c}{\partial X_c}	
\\
\end{pmatrix}
\]
is invertible in $\completion[X_1,\dots,X_\dimalg, X_{\dimalg + 1}]_g/(f_1,\dots,f_c)$. This implies that passing to the quotient by the prime ideal corresponding to $\tuple{x}$, i.e.\ evaluating $\mathcal{J}_f^c$ in $\tuple{x}$, the matrix stays invertible which in turn gives  that $\im \mathcal{J}_f^{c}(\tuple{x})$ is isolated. So by \cref{prop:rk_jacobian}, we conclude that the derived Lie sublattice of $\ker\cmatrix(\tuple{x})$ is isolated.\par
The last results of this section gives a criterion for geometrical smoothness that will be useful in the applications of \cref{part:sl4}.
\begin{prop}
\label{prop:geo_sm_crit}
Assume $\ranklocus{\cmatrix}{2\mu}{}\otimes_\completion \C$ is smooth and irreducible and that all $\tuple{x}\in\ranklocus{\cmatrix}{2\mu}{(\qfield)}$ are such that 
\[
\dim_{\qfield} [\ker \overline{\cmatrix}(\tuple{x}), \ker \overline{\cmatrix}(\tuple{x})] = \dimalg - \dim_\C (\ranklocus{\cmatrix}{2\mu}{}\otimes_\completion \C).
\]
Then $\ranklocus{\cmatrix}{2\mu}{(\completion)}$ is geometrically smooth.
\begin{proof}
Let for convenience $X = \ranklocus{\cmatrix}{2\mu}{}$. Take $\tuple{x}\in X{(\completion)}$ such that $\tuple{x}\in X(\qfield)$ modulo $\primeideal$ and let $\mathfrak{r} =[\ker \cmatrix(\completion), \ker \cmatrix(\completion)]$. Fix an embedding of $\completion$ into $\C$. The kernel of $\cmatrix(\tuple{x})$ when seen as a complex matrix is 
\[
\ker_\C \cmatrix(\tuple{x}) = \ker \cmatrix(\tuple{x})\otimes_\completion \C
\]
and consequently $[\ker_\C \cmatrix(\tuple{x}),\ker_\C \cmatrix(\tuple{x})] = \mathfrak{r}\otimes_\completion \C$. Since $\ranklocus{\cmatrix}{2\mu}{}\otimes_\completion \C$ is smooth, we must have that 
\[
{\rk}_\completion \mathfrak{r} = \dim_\C (\mathfrak{r} \otimes_\completion \C) = \dimalg - \dim_\C (X\otimes_\completion \C).
\]
Where the last equality holds by \cref{prop:rk_jacobian}. If $\dimalg - \dim_\C (X\otimes_\completion \C)$ is also the dimension over $\qfield$ of the reduction $\bmod$ $\primeideal$ of $\mathfrak{r}$ it follows that the latter $\completion$-submodule of $\completion^\dimalg$ must be isolated. Finally the presence of rank-preserving lifts follows by Hensel's lemma.
\end{proof}
\end{prop}
\section{Reduction to the finite field}
\label{sec:smooth_loci}
The aim of this section is to prove \cref{main:C} and its corollaries. Let then $\cmatrix$ have geometrically smooth rank-loci. For $\level\in\N$, the matrix obtained from $\cmatrix$ reducing its entries $\bmod\,\primeideal^\level$ is denoted by $\cmatrix^\level$. We start by proving a preliminary result involving a coarser analog of \cref{def:cal_F}.
\begin{defn}
\label{def:fcoeff}
Let  $I=\lbrace i_1,\dots,i_\lastindex \rbrace_{<}\subseteq [\halfdim-1]_0$.
\begin{enumerate}
\item If $I = \emptyset$ we define $\fcoeff{I}{(\cmatrix)} = \lbrace \emptyset\rbrace$.
\item If $I \neq \emptyset$, we define $\fcoeff{I}{(\cmatrix)}$ as the set of 
$(\tuple{x}_1,\dots,\tuple{x}_\lastindex)\in((\qfield^\dimalg)^*)^\lastindex$ with the following properties: for all $j=1,\dots,\lastindex$ and for $\tuple{y}_j = \sum_{k = j}^{\lastindex} \tuple{x}_k$
\begin{enumerate}
\item \label{eq:coeff_1}	${\rk}_{\qfield} \overline{\cmatrix}\left(\tuple{y}_j\right) = 2(\halfdim-i_j)$,
\item\label{eq:coeff_2}	$\tuple{x}_{j-1} \in [\ker \overline{\cmatrix}\left(\tuple{y}_j\right), \ker \overline{\cmatrix}\left(\tuple{y}_j\right)]$ for $j > 1$.
\end{enumerate}
\end{enumerate}
\end{defn}
We introduce the following concept.
\begin{defn}
\label{def:choice_rk}
 A \emph{choice of rank-preserving lifts} is a function $\choice:\qfield^\dimalg\rightarrow\completion^\dimalg$ such that, for all $2k\leq\dimalg$ and all $\tuple{x}\in\qfield^\dimalg$ with $\rk_{\qfield}\overline{\cmatrix}(\tuple{x}) = 2k$, $\rk_{\localfrac}\cmatrix(\choice(\tuple{x})) = 2k$. It is worth noticing here that we require $\varphi$ to be just a function.\par
\end{defn}
When $\cmatrix$ has geometrically smooth rank-loci a choice of rank-preserving lifts always exists by definition. In addition, let $\mathrm{Iso}(\alginf)$ be the set of all isolated submodules of $\alginf$. By \cite{rei2003orders}*{Theorem 4.0} every isolated submodule of $\alginf$ admits a complement.  Using the axiom of choice we may then  define a function 
\[
\compl:\mathrm{Iso}(\alginf) \rightarrow \mathrm{Iso}(\alginf)\]
assigning a complement to each isolated submodule of $\alginf$. Let us fix such a function for all the rest of this work.
The following proposition is the fundamental step in the proof of \cref{main:C}.
\begin{prop}
\label{lem:theta}
Let $\mathrm{r}_I=(r_{1},\dots,r_{\lastindex})\in \N^{\lvert I\rvert}$ and let $\pcoeff{I}{(\cmatrix)}$ be defined as in \cref{eq:N_sets}. Assume further that $\cmatrix$ has geometrically smooth rank-loci and let $\choice$ be a choice of rank-preserving lifts. Then there is a surjective function
\[
\redfunction{\choice,\compl}: \pcoeff{I}{(\cmatrix)}\rightarrow\fcoeff{I}{(\cmatrix}),
\]
\end{prop}
such that its fibre above $(\tuple{x}_1,\dots,\tuple{x}_\lastindex) \in \fcoeff{I}{(\cmatrix)}$ has cardinality
\[
q^{(\dimalg - \dimalg'_{\tuple{x}_1}) (\level_1 - 1)} \cdot \prod_{i=2}^{\lastindex} q^{(\dimalg - \dimalg'_{\tuple{x}_i}) \level_i},
\]
where, for all  $\tuple{x}\in\qfield^\dimalg$,
\begin{align*}
\dimalg'_\tuple{x}	&= \dim_{\qfield} [\ker\overline{\cmatrix}(\tuple{x}), \ker\overline{\cmatrix}(\tuple{x})].
\end{align*}
\begin{rem}
In the notation of \cref{lem:theta}, the function $\redfunction{\choice,\compl}$ depends on the non-canonical choices of $\choice$ and $\compl$ but $\fcoeff{I}{(\cmatrix})$ does not.  However, even if we shall obtain different $\redfunction{\choice,\compl}$'s for different choices of $\choice$ and $\compl$, we shall use $\redfunction{\choice,\compl}$ only to compute the cardinality of $\pcoeff{I}{(\cmatrix)}$ using the cardinality of $\fcoeff{I}{(\cmatrix)}$. These arbitrary choices will, therefore, have no consequence for our purposes.
\end{rem}
\subsection{Proof of \cref{lem:theta}} The proof  is a double induction on the cardinality of $I$ and on $\level_\lastindex$. Before starting, however, we need to fix some notation and prove some preliminary results.
\subsubsection{Notation}\label{sec:notation_proof_sigma}
For convenience we first of all introduce the following notation and definitions:
\begin{enumerate}
\item for $f\in \N$ and $e\in\N \cup \lbrace \infty \rbrace$, we write
\[
\bl{f}{e}=\begin{pmatrix}
0		&\pi^e	&0		&\dots	&\dots	&0\\
-\pi^e	&0		&		&		&		&\vdots\\
0		&		&\ddots	&		&		&\vdots\\
\vdots	&		&		&\ddots		&		&0\\
\vdots	&		&		&		&0		&\pi^e\\
0		&\dots	&\dots	&0		&-\pi^e	&0\\
\end{pmatrix}\in\Mat{2f\times 2f}{(\completion)}
\]
for the $2f\times2f$ antisymmetric matrix with $f$ blocks of the form
\[
\begin{pmatrix}
0&\pi^e\\
-\pi^e&0
\end{pmatrix}
\]
on the diagonal and zeroes elsewhere. Here $\pi^\infty = 0$ by convention.
\item Let $f\in \N$ and $M$ be an antisymmetric $2f\times 2f$ matrix over $\completion$. Then there is $S\in\GL{2f}{(\completion)}{}$, $i\in\N$ and $f_1,\dots, f_i\in \N$ and  $e_1,\dots, e_i\in\N \cup \lbrace \infty \rbrace$ such that $f = f_1 +\dots +f_i$ and 
\[
S^t M S = \begin{pmatrix} \bl{f_1}{e_1}	&0		&0\\
					0			&\ddots	&0\\
					0			&0		&\bl{f_i}{e_i}
\end{pmatrix}.
\]
If we assume that $e_1< e_2<\dots<e_i$ then the right-hand side is unique among the matrices that may be obtained from $M$ by multiplying it on the left with $T\in\GL{2f}{(\completion)}{}$ and on the right with $T^t$. We call this the \emph{antisymmetric Smith normal form} of $M$ (ASNF for short) and denote it by $\asnf (M)$. The ASNF of a $(2f + 1)\times (2f + 1)$ antisymmetric matrix over $\completion$ is defined in the obvious analogous way with the only difference that the last block on the diagonal consists only of a $0$.
%
%
%
\item
Let $\widehat{\tuple{x}}\in\completion^\dimalg$. We denote by $\omega_{\widehat{\tuple{x}}}$ the element of $\Hom_\completion(\alginf,\completion)$ whose coordinates in the dual basis $\basis^\sharp$ are given by $\widehat{\tuple{x}}$.
\item\label{def:bilinear}
Let $\mathfrak{h}$ be an $\completion$-Lie lattice and $\omega\in\Hom_\completion(\mathfrak{h},\completion)$. We define a bilinear form $b_\omega$ on $\mathfrak{h}$ as follows: for $x,y\in\mathfrak{h}$
\[
b_\omega(x,y) = \omega([x,y]).
\]
\item
In the same notation as the previous point we set
\begin{align*}
\rad(\omega) &= \lbrace x\in\mathfrak{h} \mid b_\omega(x,y) = 0 \; \forall\, y\in\mathfrak{h}\rbrace\\
\rad_\level(\omega) &= \lbrace x\in\mathfrak{h} \mid b_\omega(x,y) \in \primeideal^\level \; \forall\, y\in\mathfrak{h}\rbrace & \level\in\N.
\end{align*}
It is a straightforward application of Jacobi's identity to notice that $\rad(\omega)$ and $\rad_\level(\omega)$ ($\level\in\N$) are Lie sublattices of $\mathfrak{h}$.
\end{enumerate}
\subsubsection{Isolated radicals}
Before proceeding we also need to recall the following facts about radicals.
\begin{lem}
\label{lem:rad_isolated}
Let $\omega\in\Hom_\completion(\alginf,\completion)$. Then $\rad(\omega)$ is isolated in $\alginf$.
\begin{proof}
Let $\alpha$ be a non-zero divisor in $\completion$. Then $\alpha y\in\rad(\omega)$ if and only if 
\[
\omega([\alpha y, x]) = 0
\]
for all $x\in\alginf$. This is equivalent to $\alpha\omega([ y, x]) = 0$ for all $x\in\alginf$, which in turn is equivalent to $\omega([ y, x]) = 0$ for all $x\in\alginf$ because $\alpha$ is not a zero divisor. Hence $\alpha y\in\rad(\omega)$ if and only if $y\in\omega$, i.e.\ $\alginf/\rad(\omega)$ is torsion-free. 
\end{proof}
\end{lem}
\subsubsection{Rank-preserving lifts}
Before we start proving \cref{lem:theta} we also need an analog of \cite{zor2016adjoint}*{Theorem C} in describing how the elementary divisors of the commutator matrix behave under lifting.
\begin{lem}
\label{lem:zero_omega}
Let $\mathfrak{h}$ be an $\completion$-Lie lattice and let $\omega\in\Hom_\completion(\mathfrak{h},\completion)$. Then $b_\omega = 0$ if and only if 
\[
\omega_{\lvert \isol([\mathfrak{h},\mathfrak{h}])} = 0.
\] 
\begin{proof}
The statement is a straightforward consequence of the definition of $b_\omega$ (cf.\ \cref{sec:notation_proof_sigma}).
\end{proof}
\end{lem}
The following result constructively describes the ASNF of $\cmatrix$ when evaluated at lifts of elements of $(\completion_\level)^\dimalg$.
\begin{lem} 
\label{lem:jp_rp_lifts}
Let $\mu \leq \halfdim$. Let $\widehat{\tuple{x}}\in\completion^\dimalg$ be such that 
\[
\asnf (\cmatrix(\widehat{\tuple{x}}))	\equiv M = \begin{pmatrix}
	\bl{\mu}{0}		&0\\
	0			&0
	\end{pmatrix}								\mod \primeideal^\level.
\]
Let $V = \rad (\omega_{\widehat{\tuple{x}}})$ and $Z = \isol([V,V])$. Let $\mathcal{C}$ be an $\completion$-basis of $V$ and $\mathcal{C}'\subseteq \mathcal{C}$ be an $\completion$-basis of $Z$. In addition, let also $\widehat{\tuple{z}}\in\completion^\dimalg$ and $\widehat{\tuple{z}}_{V}$ be the coordinates with respect to $\mathcal{C}^\sharp$ of the restriction of $\omega_{\widehat{\tuple{z}}}$ to $V$. Then 
 \[
 \asnf(\cmatrix(\widehat{\tuple{x}} + \pi^\level \widehat{\tuple{z}}))	\equiv	\begin{pmatrix}
															\bl{\mu}{0}		&0\\
															0			&\pi^\level B
															\end{pmatrix}	\mod \primeideal^{\level + 1}
 \]
where $B = \asnf(\cmatrix_{\mathcal{C}}(\widehat{\tuple{z}}_{V}))$.
Moreover, if $\widehat{\tuple{z}}_{Z}$ are the coordinates with respect to $\mathcal{C}'^\sharp$ of the restriction of $\omega_{\widehat{\tuple{z}}}$ to $Z$, then 
\[
\asnf(\cmatrix(\widehat{\tuple{x}} + \pi^\level \widehat{\tuple{z}}))\equiv	M	\mod \primeideal^{\level + 1}
\]
if and only if $\widehat{\tuple{z}}_{Z}\equiv 0\,\bmod\,\primeideal$.
\begin{proof}
Let $\omega = \omega_{\widehat{\tuple{x}}}$ and $\omega' = \omega_{\widehat{\tuple{x}} + \pi^\level \widehat{\tuple{z}}}$. We may assume that $\cmatrix(\widehat{\tuple{x}})$ is already in ASNF and that $\mathcal{C} = \lbrace b_{2\mu +1}, \dots, b_\dimalg\rbrace$.
Then
\[
\cmatrix(\widehat{\tuple{x}} + \pi^\level \widehat{\tuple{z}})	=	\begin{pmatrix}
												\bl{\mu}{0} +\pi^\level A	& \pi^\level C\\
												-\pi^\level \transpose{C}	& \pi^\level B
												\end{pmatrix}
\]
for suitable $A\in\mathrm{Asym}_{2\mu}(\completion)$ and $C\in\Mat{2\mu \times (\dimalg - 2\mu)}{(\completion)}$. Now, since $\bl{\mu}{0} +\pi^\level A$ is invertible, we may use its entries to cancel $\pi^\level C$ and $-\pi^\level \transpose{C}$ with simultaneous elementary row and column operations. In practice this is equivalent to shifting $b_{2\mu + 1},\dots, b_{\dimalg}$ with linear combinations of $b_1,\dots b_\dimalg$ multiplied by $\pi^\level$. The result of these manipulations is an $\completion$-basis of $\alginf$, say $\mathcal{H}$, such that the matrix of $b_{\omega'}$ with respect to $\mathcal{H}$ is 
\[
M' = \begin{pmatrix}
\bl{\mu}{0} +\pi^\level A	& 0\\
0	& \pi^\level B + \pi^{2 \level} D
\end{pmatrix}
\]
where $D\in\mathrm{Asym}_{\dimalg - 2\mu}(\completion)$. This proves the first part of the statement, while the second part now follows by \cref{lem:zero_omega}.
\end{proof}
\end{lem}
\begin{defn}
Let $\level\in\N$. We say that $\tuple{y} \in(\completion/\primeideal^{\level + 1})^\dimalg$ is a rank-preserving lift of $\tuple{x}\in (\completion/\primeideal^{\level})^\dimalg$ when $\tuple{y}\equiv \tuple{x}\,\bmod\,\primeideal^\level$ and $\cmatrix^{\level + 1}(\tuple{y})$ and has as many maximal elementary divisors as $\cmatrix^{\level}(\tuple{x})$, i.e.\  
	\[
	\lvert \lbrace a\in \nu_{\cmatrix, \level + 1}(\tuple{y})\mid a = \level + 1\rbrace\rvert =\lvert \lbrace a\in \nu_{\cmatrix, \level}(\tuple{x})\mid a = \level \rbrace \rvert.
	\]
\end{defn}
It is now straightforward to deduce the following
\begin{cor}
\label{cor:n_lifts}
In the notation of \cref{lem:jp_rp_lifts}. There are exactly 
	\[
	q^{\dimalg - \rk_\completion Z}
	\]
distinct rank-preserving lifts of $\reducinf{\level}(\widehat{\tuple{x}})$ in $(\completion/\primeideal^{\level + 1})^\dimalg$.
\end{cor}
\subsubsection{Basis step}
The case in which $I = \emptyset$ is treated separately in the definition of $\fcoeff{I}{(\cmatrix)}$; here $\lvert \fcoeff{I}{(\cmatrix)}\rvert = \lvert \pcoeff{I}{(\cmatrix)}\rvert = 1$ and the statement of \cref{lem:theta} is readily verified setting $\redfunction{\choice,\compl}$ to be the only function between the two sets. The induction, then, ought to start by considering $\lvert I\rvert = 1$.\par
Let $i<\halfdim$. If $I = \lbrace i\rbrace$, then $\fcoeff{I}{(\cmatrix)}$ consists of elements $\tuple{x}\in\qfield^\dimalg$ such that 
\[
{\rk}_{\qfield} \overline{\cmatrix}(\tuple{x}) = 2 (\halfdim - i).
\]
When $\level_1 = 1$, the sets $\fcoeff{I}$ and $\pcoeff{I}{(\cmatrix)}$ coincide, so also in this case the statement of \cref{lem:theta} is immediately verified. Let $\level_1 = 2$. By the geometrical smoothness hypothesis, $[\rad(\omega_{\choice(\tuple{x})}),\rad(\omega_{\choice(\tuple{x})})]$ is isolated. It follows that the elements of $\pcoeff{I}{(\cmatrix)}$ lifting $\tuple{x}$ are exactly $q^{\dimalg - \dimalg'_{\tuple{x}}}$ by \cref{cor:n_lifts}. So \cref{lem:theta} follows by setting 
\[
\redfunction{\choice,\compl} = \reduc{2}{1}.
\]
More generally, when $\level_1 > 2$,  we set $\redfunction{\choice,\compl} = \reduc{\level_1}{1}$ and verify that this function has the desired properties by means of a recursive application of \cref{cor:n_lifts}.
\subsubsection{Inductive hypothesis}
We proceed now with the proof of the inductive step. Fix $I = \lbrace i_1,\dots, i_\lastindex\rbrace_{<}$ with $\lastindex > 1$ and $\tuple{\level}_I\in\N^\lastindex$.  We make two distinct inductive hypotheses according to whether $\level_\lastindex > 1$ or not:
\begin{enumerate}
\item in case $\level_\lastindex > 1$ we set $J = I$, $\tuple{r}_J = \lbrace \level_1,\dots,\level_{\lastindex - 1}, \level_\lastindex - 1\rbrace$ and assume that we have already defined
a surjective function
\[
\redfin_{J,\mathrm{r}_J,\choice,\compl} : N_{J,\mathrm{r}_{J}}^{\completion}{(\cmatrix)}\rightarrow\fcoeff{J}{(\cmatrix}),
\]
whose fibres have the cardinality specified by the expression in \cref{lem:theta}.
\item In case $\level_\lastindex = 1$ we set $J = I \smallsetminus\lbrace i_\lastindex\rbrace$, $\tuple{r}_J = \lbrace \level_1,\dots,\level_{\lastindex - 1}\rbrace$ and assume that we have already defined a surjective function
\[
\redfin_{J,\mathrm{r}_J,\choice,\compl} : \pcoeff{J}{(\cmatrix)}\rightarrow\fcoeff{J}{(\cmatrix}),
\]
whose fibres have the cardinality specified by the expression in \cref{lem:theta}.\par 
\end{enumerate}
\subsubsection{Definition of $\redfunction{\choice,\compl}$, for $\level_\lastindex > 1$}
\label{sec:r_l>1}
Let $\tuple{w}\in\pcoeff{I}{(\cmatrix)}$ and  $\widehat{\tuple{w}}\in \completion^\dimalg$ lifting $\tuple{w}$. In order to make use of the inductive hypothesis, we now show how to define a surjective function
\[
\sigma: \pcoeff{I}{(\cmatrix)}\rightarrow \pcoeff{J}{(\cmatrix)}
\]
whose fibre above each $\tuple{w}'\in\pcoeff{J}{(\cmatrix)}$ has cardinality $q^{\dimalg - e}$ where
\[
e = \dim_{\qfield} [\ker\overline{\cmatrix}(\reduc{N-1}{1}(\tuple{w}')), \ker\overline{\cmatrix}(\reduc{N-1}{1}(\tuple{w}'))].
\]
The inductive hypothesis will then make sure that $ \redfin_{I,\mathrm{r}_I,\choice,\compl} = \redfin_{J,\mathrm{r}_J,\choice,\compl}\circ \sigma$ has the required properties.\par
First of all we bring $\cmatrix(\widehat{\tuple{w}})$ in ASNF. Naturally there will be many basis changes that do so in general, what follows explains how to make sure that $\sigma$ is well defined. Let $\omega = \omega_{\widehat{\tuple{w}}}$ and $\tuple{x} =\reduc{N}{1}(\tuple{w})$. By \cref{lem:rad_isolated} we may define $U$ as $\compl(\rad(\omega_{\choice (\tuple{x})}))$.  This definition of $U$ causes the dependence of $\redfunction{\choice,\compl}$ on $\choice$ and $\compl$.\par
Clearly $b_\omega$ is non-degenerate on $U$ and hence
	\[
		V =  \lbrace x\in\alginf\mid b_\omega(x,u) = 0\; \forall u\in U\rbrace
	\]
is a complement of $U$ in $\alginf$. One checks that not only is $V\equiv \rad(\omega_{\choice (\tuple{x})})\,\bmod\,\primeideal$ but also $V$ is a Lie sublattice with isolated derived Lie sublattice. Indeed the first assertion easily follows by construction. The second instead is a consequence of the rank-loci of $\alginf$ being geometrically smooth: 
\begin{lem}
\label{lem:V_is_subalgebra}
In the notation fixed above. There is $\widehat{\tuple{x}}\in\completion^\dimalg$ such that $V = \rad (\omega_{\widehat{\tuple{x}}})$ and $[V,V]$ is an isolated $\completion$-submodule of $\alginf$.
\begin{proof}
Let us first of all fix some notation and rephrase the problem in terms of lifting approximate solutions with Hensel's lemma. Let $\basis'$ be an $\completion$-basis of $\alginf$ with the first $2\mu_\lastindex$ elements in $U$ and the remaining ones in $V$. Let  also $\cmatrix' = \cmatrix_{\basis'}$. Let $A = \completion[X_1,\dots,\ X_\dimalg]$. For $i,j\in\lbrace 1,\dots,\dimalg\rbrace$ denote by $R_{i,j}$ the $(i,j)$-th entry of $\cmatrix'$. Let
\begin{align*}
I_T &= \left (\lbrace R_{ij}\mid i\leq 2\mu_\lastindex \; j > 2\mu_\lastindex \rbrace \right ) \\
I_C &= \left (\lbrace R_{ij}\mid i,j > 2\mu_\lastindex \rbrace \right ).
\end{align*}
Let $\widehat{\tuple{w}}'$ be the coordinates of $\omega$ with respect to $\basis'^\sharp$.
We have that $\widehat{\tuple{w}}'$ is, $\bmod\,\primeideal$, a common zero of $I_T\cup I_C$. So, by Hensel's lemma, if $I_T + I_C$ is generated by linear polynomials in $A^* = A\smallsetminus \primeideal\cdot A$, we find  an exact solution $\widehat{\tuple{x}}'\in\completion^\dimalg$ such that $\widehat{\tuple{x}}'\equiv \widehat{\tuple{w}}'\,\bmod\,\primeideal$. This exact solution will be so that
\begin{align*}
\cmatrix'(\widehat{\tuple{x}}') &= \begin{pmatrix}
						M	& 0\\
						0	&0
						\end{pmatrix}		&&M\in\GL{2\mu_\lastindex}{(\completion)}{}
\end{align*}
and therefore $V = \ker \cmatrix'(\widehat{\tuple{x}}') = \rad (\omega_{\widehat{\tuple{x}}})$ for an appropriate $\widehat{\tuple{x}}\in\completion^\dimalg$, thus proving the first part of the lemma.\par
Let us prove the claim above. First of all we notice that by the construction of $U$ and $V$,
\begin{align*}
R_{ij}(\widehat{\tuple{w}}) = 0&&\forall\, i\leq 2\mu_\lastindex,\, j > 2\mu_\lastindex.
\end{align*}
So for the purposes of this proof we may assume that $I_T$ is generated by linear polynomials in $A^*$. The ideal $I_C$ may be generated in this way too, for the rank-locus is geometrically smooth which also implies that the Lie lattice $[V,V]$ is isolated. This concludes the proof.
\end{proof}
\end{lem}
%
%
%
%
%
%
%
%
%
%
What follows defines $\sigma$ and proves its surjectivity.  In order to do this we need a sequence of three lemmata. We state and prove these preliminary results for slightly more general $U$ and $V$ than the ones above. This will enable us to use them without modifications also when $\level_\lastindex = 1$. The first lemma shows that we may bring $\cmatrix(\widehat{\tuple{w}})$ in ASNF without altering the decomposition $\alginf = U \oplus V$.  This will secure a non-ambiguous definition of $\sigma$.
\begin{lem}
\label{lem:decomposition}
Let $U$ be an isolated $\completion$-submodule of $\alginf$ such that $b_\omega$ is non-degenerate when restricted to it. Let also $V$ be the orthogonal complement to $U$ with respect to $b_\omega$. Then there are isolated complementary $\completion$-submodules $U',V_1,\dots,V_\lastindex \subseteq \alginf$ such that $b_\omega$ is non degenerate when restricted to $U'$ and
\begin{align*}
U'							&\supseteq U\\
\oplus_{i=1}^{\lastindex} V_i		&\subseteq V\\
 \rad_{\level_\lastindex} (\omega)	&= \pi^{\level_\lastindex} U' + \oplus_{i=1}^{\lastindex} V_i\\
  \rad_{\level_\lastindex + \level_{\lastindex -1}} (\omega)	&= \pi^{\level_\lastindex+ \level_{\lastindex -1}} U' + \pi^{ \level_{\lastindex -1}} V_\lastindex +\oplus_{i=1}^{\lastindex -1 } V_i\\
							& \vdots\\
 \rad_{N} (\omega)				&= \pi^{N} U' + \pi^{N - \level_{\lastindex}} V_\lastindex + \pi^{N - \level_{\lastindex}- \level_{\lastindex-1}} V_{\lastindex -1} +\cdots+V_1.							
\end{align*}
Moreover $\rk_\completion U' - rk_\completion U$ is even.
\begin{proof}
The claim follows directly from the algorithm to bring the matrix of $b_\omega$ with respect to $\basis$ to its ASNF. The last statement is a consequence of the fact that a submodule of $\alginf$ needs to have even $\completion$-rank for an antisymmetric form to be non-degenerate when restricted to it.
\end{proof}
\end{lem}
The first of the following two lemmata will define $\sigma$ while the second will be used to prove its surjectivity.
\begin{lem}
\label{lem:going_down}
Let $U,V,U',V_1,\dots,V_\lastindex $ be as in \cref{lem:decomposition}. Define
\begin{align*}
I_{-}&= \begin{cases}
					\lbrace i_1,\dots,i_\lastindex\rbrace		& \level_\lastindex > 1\\
					\lbrace i_1,\dots, i_{\lastindex - 1}\rbrace	&\level_\lastindex = 1
					\end{cases}\\
\tuple{\level}_{-} &=\begin{cases}
					\lbrace \level_1,\dots,\level_{\lastindex - 1}, \level_\lastindex - 1\rbrace	& \level_\lastindex > 1\\
					\lbrace \level_1,\dots,\level_{\lastindex - 1}\rbrace &\level_\lastindex = 1.
					\end{cases}
\end{align*}
Assume that $V = \rad_{\level_\lastindex} (\omega)$ and that it has isolated derived sublattice $Z = [V,V]$. Let $W = \compl(Z)$. Then
\begin{align*}
\omega_{-} (w)	&=\omega(w)\,\forall\, w\in W\\
\omega_{-} (z)	&=\pi^{-1} 	\omega(z)\,\forall\, z\in Z
\end{align*}
defines an $\completion$-linear form on $\alginf$ whose coordinates with respect to $\basis^\sharp$ are in $N_{I_{-},r_{-}}^\completion(\cmatrix)$ when reduced $\bmod\,\primeideal^{N-1}$.\par
\begin{proof}
First of all one has to prove that it is actually possible to divide the values of $\omega$ on $Z$  by $\pi$. This is immediate because $V = \rad_{\level_\lastindex}(\omega)$ implies that $\omega( Z )\subseteq\primeideal^{\level_\lastindex}$.  We now need to show that the matrices of $b_\omega$ have the required ASNF. Fix a basis $\otherbasis = \lbrace \ob_1,\dots,\ob_\dimalg\rbrace$ bringing the matrix of $b_\omega$ in ASNF and such that
\begin{align*}
\ob_1,\dots,\ob_{2\mu_\lastindex}	&\in U\\
\ob_{2\mu_\lastindex + 1},\dots,\ob_{2\mu_{\lastindex -1}}	&\in V_\lastindex\\
\vdots&\\
\ob_{2\mu_1 + 1},\dots,\ob_{\dimalg}	&\in V_1.
\end{align*} 
Since $b_\omega$ restricted to $V_\lastindex$ coincides with a non-degenerate bilinear form on $V_\lastindex$ multiplied by $\pi^{\level_\lastindex}$ we may replace $\ob_1,\dots,\ob_{2\mu_\lastindex}$ with independent $\ob'_1,\dots,\ob'_{2\mu_\lastindex}$ in such a way that, for all $i\in\lbrace 1,\dots 2\mu_\lastindex\rbrace$  and $j\in \lbrace 2\mu_\lastindex + 1,\dots, 2\mu_{\lastindex - 1}\rbrace$
\[
b_\omega(\ob'_i,\ob_j) = \omega(w_{ij})\in\primeideal^{\level_\lastindex}  \text{ for some $w_{ij}\in W$.}
\]
Similarly, since $b_\omega$ is non-degenerate on $U'' = \linspan(\ob'_1,\dots,\ob'_{2\mu_\lastindex})$, we may further replace $\ob'_1,\dots,\ob'_{2\mu_\lastindex}$ with $\ob''_1,\dots,\ob''_{2\mu_\lastindex}\in U''$ such that, for all $i,j\in\lbrace 1,\dots 2\mu_\lastindex\rbrace$,
\[
b_\omega(\ob''_i,\ob''_j) = \omega(w_{ij}) \text{ for some $w_{ij}\in W$.}
\]\todo[color = yellow]{see 28.02 for details}
This has the following consequence: if $\otherbasis''$ is the basis obtained replacing $\ob_1,\dots,\ob_{2\mu_\lastindex}$ with $\ob''_1,\dots,\ob''_{2\mu_\lastindex}$, then there are
\begin{align*}
 C_{j}	&\in\Mat{2\mu_\lastindex \times 2\mu_{j}}{(\completion)}		&j=2,\dots,\lastindex 	- 1
\end{align*}
such that the matrix of $b_{\omega_-}$ with respect to $\otherbasis'$ and 
\[
\begin{pmatrix}
\bl{\mu_\lastindex}{0}			&\pi^{\level_\lastindex} C_{\lastindex - 1}				&\pi^{\level_{\lastindex } +\level_{\lastindex -1} - 1}C_{\lastindex-2}					&\dots	&\pi^{N-\level_1 - 1}C_{2}								&0\\
				&\bl{\mu_{\lastindex - 1}}{\level_\lastindex - 1}		&0												&\dots	&0												&\vdots\\
				&							&\bl{\mu_{\lastindex- 2}}{\level_{\lastindex } +\level_{\lastindex -1} - 1}			&\ddots	&\vdots											&\\
				&							&												&\ddots	&0												&\vdots\\
				&							&												&		&\bl{\mu_2}{N-\level_1 - 1}							&0\\
				&							&												&		&												&0				
\end{pmatrix}
\] 
coincide $\bmod\,\primeideal^{N-1}$. The conclusion now follows by cancelling the $C_i$'s with the first block on the diagonal using elementary row and column operations.
\end{proof}
\end{lem}
\begin{lem}
\label{lem:going_up}
Let $U,V,U',V_1,\dots,V_\lastindex $ be as in \cref{lem:decomposition}. Define $2\mu = \rk_\completion U' - rk_\completion U$ and
\begin{align*}
I_+&= \begin{cases}
					\lbrace i_1,\dots,i_\lastindex\rbrace		& V\subseteq\rad_{\level_\lastindex}(\omega)\\
					\lbrace i_1,\dots, i_\lastindex, i_\lastindex + \mu\rbrace	&\text{ otherwise}
					\end{cases}\\
\tuple{\level}_+ &=\begin{cases}
					\lbrace \level_1,\dots,\level_{\lastindex - 1}, \level_\lastindex + 1\rbrace	& V\subseteq\rad_{\level_\lastindex}(\omega)\\
					\lbrace \level_1,\dots,\level_{\lastindex - 1},\level_{\lastindex}, 1\rbrace &\text{ otherwise}.
					\end{cases}
\end{align*}
Assume that $V$ is a Lie-sublattice of $\alginf$ with isolated derived sublattice $Z = [V,V]$. Let $W = \compl(Z)$. Assume further there is $\omega'\in\Hom_\completion(\alginf,\completion)$ coinciding modulo $\primeideal$ with $\omega$ on $W$ and such that $V$ is equal to $\rad(\omega')$ modulo $\primeideal$. Then
\begin{align*}
\omega_+ (w)	&=\omega(w)\,\forall\, w\in W\\
\omega_+ (z)	&=\pi\omega(z)\,\forall\, z\in Z
\end{align*}
defines an $\completion$-linear form on $\alginf$ whose coordinates with respect to $\basis^\sharp$ are in $N_{I_+,r_+}^\completion(\cmatrix)$ when reduced $\bmod\,\primeideal^{N+1}$.\par
\begin{proof}
Fix a basis $\otherbasis = \lbrace \ob_1,\dots,\ob_\dimalg\rbrace$ bringing the matrix of $b_\omega$ in ASNF and such that
\begin{align*}
\ob_{1},\dots,\ob_{2\mu}							&\in U\\
\ob_{2\mu + 1},\dots,\ob_{2\mu_\lastindex}			&\in U'\\
\ob_{2\mu_\lastindex + 1},\dots,\ob_{2\mu_{\lastindex -1}}	&\in V_\lastindex\\
\vdots&\\
\ob_{2\mu_1 + 1},\dots,\ob_{\dimalg}					&\in V_1.
\end{align*} 
Since $\omega_{+}$ and $\omega'$ coincide modulo $\primeideal$, a straightforward application of \cref{lem:jp_rp_lifts} ensures that the matrix of $b_{\omega_+}$ with respect to $\otherbasis$ coincides with
\[
\begin{pmatrix}
\bl{\mu}{0}		&0					& 0\\
0		&\bl{\mu_\lastindex -\mu}{1}	& 0\\
0			&0						&0
\end{pmatrix}%
\mod \primeideal^2.
\]
Now fix $i\in\lbrace 2,\dots,\lastindex\rbrace$. If $\ob_j,\ob_k\in\otherbasis$ are such that $1\leq j\leq 2\mu$ and $2\mu_i < k \leq 2\mu_{i -1}$, then 
$b_\omega(\ob_j,\ob_k)\in\primeideal^{\level_\lastindex + \cdots +\level_i}$ and therefore also
\[
b_{\omega_+}(\ob_j,\ob_k)\in\primeideal^{\level_\lastindex + \cdots +\level_i}.
\]
Moreover, as in the proof of \cref{lem:going_down}, we may replace $\ob_1,\dots,\ob_{2\mu}$ with independent $\ob'_1,\dots,\ob'_{2\mu}$ in such a way that, for all $i\in\lbrace 1,\dots 2\mu\rbrace$  and $j\in \lbrace 2\mu_\lastindex + 1,\dots, 2\mu_{\lastindex - 1}\rbrace$
\[
b_\omega(\ob'_i,\ob_j) = \omega(w_{ij})\in\primeideal^{\level_\lastindex + 1}  \text{ for some $w_{ij}\in W$.}
\]
Let $\otherbasis'$ be the basis obtained replacing $\ob_1,\dots,\ob_{2\mu}$ with $\ob'_1,\dots,\ob'_{2\mu}$. Using the definition of $\omega_+$, it is straightforward to check that, for all $i,j\in\lbrace 2\mu, \dots, \dimalg\rbrace$, 
\[
b_{\omega_+}(\ob_i,\ob_j) = \pi b_{\omega}(\ob_i,\ob_j).
\]
This implies the following: for $j=2,\dots,\lastindex$, there are
\begin{align*}
 C_{j}	&\in\Mat{2\mu \times 2\mu_{j}}{(\completion)}	 	&j=2,\dots,\lastindex - 1\\
 C_{\lastindex}	&\in\Mat{2\mu \times 2(\mu_{\lastindex} -\mu)}{(\completion)}	 	&
\end{align*}
such that the matrix of $b_{\omega_+}$ with respect to $\otherbasis'$ coincides with 
\[
\begin{pmatrix}
\bl{\mu}{0}			&\pi^2 C_{\lastindex}				&\pi^{\level_{\lastindex} + 1}C_{\lastindex-1}					&\dots	&\pi^{N-\level_1}C_{2}								&0\\
				&\bl{\mu_{\lastindex} - \mu}{1}		&0												&\dots	&0												&\vdots\\
				&							&\bl{\mu_{\lastindex- 1}}{\level_{\lastindex } + 1}			&\ddots	&\vdots											&\\
				&							&												&\ddots	&0												&\vdots\\
				&							&												&		&\bl{\mu_2}{N-\level_1 + 1}							&0\\
				&							&												&		&												&0				
\end{pmatrix}\mod\primeideal^{N+1}.
\] 
The statement now follows immediately by cancelling the $C_i$'s with the first block on the diagonal.
\end{proof}
\end{lem}
We now reinstate $U$, $V$, $\widehat{\tuple{w}}$ and $\tuple{x}$ as defined before \cref{lem:decomposition}. We already showed that $U$ and $V$ satisfy the hypotheses of \cref{lem:decomposition,lem:going_up,lem:going_down}. Let therefore $\widehat{\tuple{w}}_{-}$ be the coordinates with respect to $\basis^\sharp$ of $\omega_{-}$ as defined in \cref{lem:going_down}. Setting $\sigma(\tuple{w}) = \reducinf{N-1}(\widehat{\tuple{w}}_{-})$ defines then by \cref{lem:going_down} a function
\[
\sigma: \pcoeff{I}{(\cmatrix)}\rightarrow \pcoeff{J}{(\cmatrix)}.
\]
Surjectivity is seen as follows: for all $\tuple{w}'\in \pcoeff{J}{(\cmatrix)}$, we produce $\tuple{w} =  \pcoeff{I}{(\cmatrix)}$ such that $\sigma(\tuple{w}) = \tuple{w}'$ by plugging in $J$ for $I$ in \cref{lem:going_up} and setting $\tuple{w}$ to be the reduction $\bmod\,\primeideal^{N}$ of the coordinates of $\omega_+$ with respect to $\basis^\sharp$. The quantitative statement on the cardinality of the fibres is now a consequence of the construction of $\omega_+$.
\subsubsection{Definition of $\redfunction{\choice,\compl}$, for $\level_\lastindex = 1$}
We finish the proof of \cref{lem:theta} by considering $\level_\lastindex = 1$. Let $\tuple{w}\in\pcoeff{I}{(\cmatrix)}$ and  $\widehat{\tuple{w}}\in \completion^\dimalg$ lifting $\tuple{w}$. Let also $U$ and $V$  be defined as at the beginning of \cref{sec:r_l>1} and set $U' = U$. We have already shown that these two $\completion$-submodules satisfy the hypotheses of \cref{lem:decomposition,lem:going_down}, we now use the inductive hypothesis to define $\redfunction{\choice,\compl}$. Namely, let $\omega = \omega_{\widehat{\tuple{w}}}$ as before and $\tuple{w}_{-}$ be the reduction $\bmod\,\primeideal^{N-1}$ of the coordinates of $\omega_-$ as defined in \cref{lem:going_down}. We define $\redfunction{\choice,\compl}(\tuple{w}) = (\tuple{x}_1,\dots,\tuple{x}_\lastindex)$ where
\begin{align*}
(\tuple{x}_1,\dots,\tuple{x}_{\lastindex - 1} + \tuple{x}_\lastindex) 	&=\redfin_{J,\mathrm{r}_J,\choice,\compl}(\tuple{w}_{-})\\
\tuple{x}_\lastindex								&=\reducinf{1}(\tuple{w}).
\end{align*}
The tuple $\redfunction{\choice,\compl}(\tuple{w})$ is easily seen to belong to $\fcoeff{I}$ because $(\tuple{x}_1,\dots,\tuple{x}_{\lastindex - 1} + \tuple{x}_\lastindex)\in\fcoeff{J}$, $\rk_{\qfield}\overline{\cmatrix}(\tuple{x}_\lastindex) = 2(\halfdim - i_\lastindex)$ and $\tuple{x}_{\lastindex - 1} \in [\ker \overline{\cmatrix}(\tuple{x}_\lastindex), \ker \overline{\cmatrix}(\tuple{x}_\lastindex)]$ by construction.\par
It remains to show surjectivity and the quantitative statement on the fibres. To see this, assume we are given $(\tuple{x}_1,\dots,\tuple{x}_{\lastindex - 1},  \tuple{x}_\lastindex)\in\fcoeff{I}$. To construct a preimage we use the inductive hypothesis and apply \cref{lem:going_up} with $J$ in place of $I$,
	\begin{align*}
		\omega	& \in \redfin_{J,\mathrm{r}_J,\choice,\compl}^{-1}((\tuple{x}_1,\dots,\tuple{x}_{\lastindex - 1} + \tuple{x}_\lastindex))\\
		\omega'	& = \omega_{\choice(\tuple{x}_\lastindex)}\\
		U		& = \compl(\rad (\omega'))\\
		V		& =  \lbrace x\in\alginf\mid b_\omega(x,u) = 0\; \forall u\in U\rbrace,
	\end{align*}
and $U'$ given by \cref{lem:decomposition}. Indeed it is easy to check that $V$ and $\rad(\omega')$ coincide modulo $\primeideal$. Moreover, since $\tuple{x}_{\lastindex - 1}\in [\ker \overline{\cmatrix}(\tuple{x}_\lastindex), \ker \overline{\cmatrix}(\tuple{x}_\lastindex)]$, it follows immediately that $\omega'$ and $\omega$ coincide on $W = \compl ([V,V])$ modulo $\primeideal$. The hypotheses of \cref{lem:going_up} are therefore satisfied and as before the reduction modulo $\primeideal^{N + 1}$ of the coordinates of $\omega_{+}$ is the sought preimage. To conclude, the construction of $\omega_-$ implies the statement on the cardinality of the fibres.
\subsection{Proof of \cref{main:C}}
Let $\arr$ be a classification by $\cmatrix$-kernels. Let $I$ be as in the previous section. It is obvious from the definition that 
\[
F_I(\cmatrix) = \bigcup_{\mathcal{S} \in \radseq^\arr_{I}(\cmatrix)} \mathcal{F}_\mathcal{S}(\cmatrix)
\]
and that the union is disjoint. The family of sets $\mathcal{F}_\mathcal{S}(\cmatrix)$ as $\mathcal{S}$ varies in $ \radseq^\arr_{I}(\cmatrix)$ is therefore a refinement of $F_I(\cmatrix)$. It also follows from \cref{def:cal_F} that, for all $\tuple{r}_I\in\N^\lastindex$,
\[
\lvert \pcoeff{I}{(\cmatrix)} \rvert = \sum_{\stackrel{\mathcal{S} \in \radseq^\arr_{I}(\cmatrix)}{\mathcal{S} = \lbrace \radiclass_1,\dots,\radiclass_\lastindex\rbrace}}%
\lvert \mathcal{F}_\mathcal{S}(\cmatrix)\rvert 
q^{(\dimalg -\dimalg'_{\radiclass_\lastindex})(\level_1 -1)}%
\prod_{i = 2}^{\lastindex} q^{(\dimalg - \dimalg'_{\radiclass_{\lastindex - i}}) r_i}
\]
It is now straightforward to compute the partial summand of the Poincar\'e series related to the above $\pcoeff{I}{(\cmatrix)}$ and obtain the formula in \cref{main:C}.\par
\subsection{Proof of \cref{main:abscissa_P}}
\Cref{main:abscissa_P} is now a corollary of \cref{main:C}. The latter gives a set of candidate poles for $\pseries{\alginf}{(s)}$, thus giving an upper bound for its abscissa of convergence. We shall now show that the abscissa cannot be smaller that the real part of the candidate pole with largest real part.\par
Let $\omega_0\in\Hom_\completion(\alginf, \completion) \smallsetminus \lbrace 0 \rbrace$ be such that $\rho_{\omega_0}$ vanishes at the maximum of
	\[
		\bigcup_{\omega\in\Hom(\algfin,\qfield) \smallsetminus \lbrace 0 \rbrace}\lbrace s\in\Q	\mid	\rho_\omega(s) = 0\rbrace.
	\]
Set
	\begin{align*}
		I 								&=\lbrace i \rbrace = \left \lbrace \frac{\dim_{\qfield} \rad(\omega_0) + 2\halfdim - \dimalg}{2}\right \rbrace\\
		M^{\completion}_{I,\level} (\cmatrix)	&= \redfunction{\choice,\compl}^{-1}(F_I(\cmatrix))	&& \level\in \N.\\
		\mathcal{P}^{abs}(s) 				&= \sum_{\level\in\N} \left \lvert M^{\completion}_{I,\level} (\cmatrix) \right\rvert \,q^{-s \level (\halfdim-i)}	&& s\in\C
	\end{align*}
It is clear that $\mathcal{P}^{abs}(s)\leq \pseries{\alginf}{(s)}$ so the latter converges at $s = s_0$ only if the former converges at $s = s_0$. A straightforward computation by means of \cref{lem:theta} gives that 
	\[
	\mathcal{P}^{abs}(s) =%
	 \lvert F_I(\cmatrix)\rvert
	 \frac{q^{ - \frac{\dimalg - \dim_{\qfield} \rad (\omega_0)}{2}s}}%
	 {1 - q^{\dimalg - \dim_{\qfield} [ \rad (\omega_0), \rad (\omega_0)] - \frac{\dimalg - \dim_{\qfield} \rad (\omega_0)}{2}s}}.
	\]
This suffices to conclude.
\section{Quadratic Lie lattices}
\label{sec:quadratic}
Assume for this section that $\alginf$ is quadratic. In this case,  there is a correspondence between kernels of the commutator matrix and centralizers in $\alginf$. The next main result takes this into account and deduces more specific versions of the formula in \cref{main:C} and a more intrinsic description of the abscissa of convergence. Set for convenience $\alginf^\sharp = \Hom_\completion(\alginf,\completion)$\par
The quadratic structure of $\alginf$ gives an isomorphism of $\completion$-modules $\lambda:\alginf\rightarrow\alginf^\sharp$. If for all $\omega_1,\omega_2\in\alginf^\sharp$, we define
\[
[\omega_1,\omega_2] =\lambda([ \lambda^{-1}(\omega_1), \lambda^{-1}(\omega_2)]),
\]
then $\alginf^\sharp$ becomes an $\completion$-Lie lattice isomorphic to $\alginf$ because the non-degenerate symmetric form defining $\lambda$ is associative with respect to the Lie bracket of $\alginf$. Fixing a basis for $\alginf$ gives an $\completion$-linear isomorphism $\iota:\completion^\dimalg\rightarrow\alginf$ and via taking the dual basis of $\basis$ this gives an $\completion$-linear isomorphism $\eta: \alginf^\sharp\rightarrow\completion^\dimalg$. The map $\eta\circ \lambda\circ \iota$ is an $\completion$-linear isomorphism of $\completion^\dimalg$ to itself and so it defines an isomorphism of $\completion$-schemes $\xi:\A^\dimalg\rightarrow \A^\dimalg$. For all $\mu\in\N_0$ such that $2\mu\leq \dimalg$ we set
\[
\centrvar{\alginf}{\dimalg - 2\mu}{} = \xi(\ranklocus{\cmatrix}{2\mu}{}).
\]
Analogously, if $\ring$ is an $\completion$-algebra, the set of $\ring$-points of $\centrvar{\alginf}{\dimalg - 2\mu}{(\ring)}$ is defined  to be the set $\xi(\ring)(\ranklocus{\cmatrix}{2\mu}{(\ring)})$.
Fix $2\mu$ as above. An argument akin to \cite[Section~5]{akov2013representation} shows that the $\completion$-points of $\centrvar{\alginf}{\dimalg - 2\mu}{}$ are exactly those whose image in $\alginf$ through $\iota$ has Lie centralizer of $\completion$-rank $\dimalg - 2\mu$. 
\begin{defn}
\label{def:geo_smooth_centr}
We call $\centrvar{\alginf}{\dimalg - 2\mu}{}$ the \emph{locus of constant centralizer dimension $\dimalg - 2\mu$} or the \emph{$(\dimalg - 2\mu)$-centralizer-locus} for short. Rephrasing \cref{def:geo_smooth_rk} in this context, we say that $\centrvar{\alginf}{\dimalg - 2\mu}{(\completion)}$ is \emph{geometrically smooth} when
\begin{enumerate}
\item for each $x\in\algfin$ with $\dim_{\qfield} \centr{\algfin}{x} = \dimalg - 2\mu$ there is $\widehat{x}\in\alginf$ such that $\rk_\completion \centr{\alginf}{\widehat{x}} = \dimalg - 2\mu$.
\item for each $\widehat{x}\in \alginf$ with  $\rk_\completion \centr{\alginf}{\widehat{x}} = \dimalg - 2\mu$ and whose reduction $x$ to $\algfin$ is such that $\dim_{\qfield} \centr{\algfin}{x} =  \dimalg - 2\mu$,
\[
[\centr{\alginf}{\widehat{x}},\centr{\alginf}{\widehat{x}}]
\]
\end{enumerate}
is an isolated submodule of $\alginf$.
\end{defn}
It is clear that by definition $\centrvar{\alginf}{\dimalg - 2\mu}{}$ is geometrically smooth if and only if $\ranklocus{\cmatrix}{2\mu}{}$ is. It is also clear that $\centrvar{\alginf}{\dimalg - 2\mu}{}$ does not depend on the choice of $\cmatrix$. In this spirit, we say that $\alginf$ has {\em geometrically smooth centralizer-loci} when $\cmatrix$ has geometrically smooth rank-loci (\cref{lem:basis_change} shows that the choice of the basis $\basis$ does not influence this definition).\par
We now definine the dual notion of kernel classes. We set
\[
\killingiso = \iota\circ\eta\circ\lambda
\]
and we denote by $\overline{\killingiso}$ the $\qfield$-linear automorphism induced on $\algfin$.
\begin{defn}
A (Lie) centralizer class $\radiclass$ is a subset of $\algfin$ such that for any two $x,x'\in\radiclass$
\begin{align*}
\dim_{\qfield} \centr{\algfin}{\overline{\killingiso}(x)}								&=  \dim_{\qfield} \centr{\algfin}{\overline{\killingiso}(x')}\\
\dim_{\qfield} [\centr{\algfin}{\overline{\killingiso}(x)},  \centr{\algfin}{\overline{\killingiso}(x)}]	&=  \dim_{\qfield} [\centr{\algfin}{\overline{\killingiso}(x')},\centr{\algfin}{\overline{\killingiso}(x')}].
\end{align*}
\end{defn}
\begin{defn}
A covering of $\algfin$ by disjoint centralizer classes is callled a \emph{classification} by centralizers. Members of a classification by centralizers $\arr$ are called centralizer $\arr$-classes or centralizer classes when no confusion is possible.
\end{defn}
As before, by definition, if $\radiclass$ is a centralizer class we have well defined
\begin{align*}
\dimalg_{\radiclass} 	&= \dim_{\qfield} \radic		&&\radic = \centr{\algfin}{\overline{\killingiso}(x)}\text{ for any } x\in\radiclass\\
\dimalg'_{\radiclass}	&= \dim_{\qfield} [\radic,\radic]	&&\radic = \centr{\algfin}{\overline{\killingiso}(x)}\text{ for any } x\in\radiclass.
\end{align*}
Let $\arr$ be a classification by centralizers. A \emph{sequence} of $\arr$-centralizer classes is a set $\lbrace \radiclass_1,\dots,\radiclass_\lastindex\rbrace$ ($\lastindex\in \N$) of centralizer $\arr$-classes such that $\dimalg_{\radiclass_1}> \dimalg_{\radiclass_2}>\cdots>\dimalg_{\radiclass_t}$. 
\begin{defn}
Let $\mathcal{S} =  \lbrace \radiclass_1,\dots,\radiclass_\lastindex\rbrace$ be a sequence of centralizer $\arr$-classes
\begin{enumerate}
\item If $\mathcal{S} = \emptyset$ we define $\ccoeffiso{\mathcal{S}}{(\algfin)} = \lbrace \emptyset\rbrace$.
\item If $\mathcal{S} \neq \emptyset$. We define $\ccoeffiso{\mathcal{S}}{(\algfin)} $ as the set of $\lastindex$-tuples $(x_1,\dots,x_\lastindex)$ of elements of $\algfin$ such that for all $j = 1,\dots,\lastindex$ for $y_j = \sum_{k = j}^{\lastindex} x_k$
\begin{enumerate}
\item $y_j \in \radiclass_j$
\item	$x_{j-1} \in [\centr{\algfin}{y_j}, \centr{\algfin}{y_j}]$ for $j > 1$.
\end{enumerate}
\end{enumerate}
\end{defn}
Fix  an ordered subset $I = \lbrace i_1,\dots, i_\lastindex\rbrace$ of $[\halfdim - 1]_0 = \lbrace 0,\dots, \halfdim - 1\rbrace$.
\begin{defn}
An $I$-sequence of centralizer $\arr$-classes is a sequence of centralizer $\arr$-classes $\lbrace \radiclass_1,\dots,\radiclass_\lastindex\rbrace$
such that $\dimalg_{\radiclass_j} = \dimalg - 2(\halfdim - i_{\lastindex + 1 - j})$ for all $j = 1,\dots,\lastindex$. The set of all $I$-sequences of centralizer $\arr$-classes is denoted by $\centrseq_{I}^\arr(\algfin)$.
\end{defn}
With this notation, the following is a straightforward consequence of \cref{main:C}
\begin{thm}
\label{main:D}
Assume $\alginf$ is quadratic and has geometrically smooth centralizer-loci.  Let $\arr$ be a classification by centralizers of $\algfin$. Then 
\[
\pseries{\alginf}{(s)} = \sum_{\stackrel{I\subseteq [\halfdim-1]_0}{I=\lbrace i_1,\dots,i_\lastindex \rbrace_{<}}}\sum_{\stackrel{\mathcal{S}\in\centrseq^\arr_{I}(\algfin)}{\mathcal{S} =  \lbrace \radiclass_1,\dots,\radiclass_\lastindex\rbrace}} \left\lvert \ccoeffiso{\mathcal{S}}{(\algfin)} \right\rvert q^{-(\dimalg - \dimalg'_{\radiclass_\lastindex})}%
\prod_{\radiclass\in\mathcal{S}} \frac{q^{\dimalg - \dimalg'_{\radiclass} - s\frac{\dimalg - \dimalg_{\radiclass}}{2}}}{1 - q^{\dimalg - \dimalg'_{\radiclass} - s\frac{\dimalg - \dimalg_{\radiclass}}{2}}}.
\]
\end{thm}
\subsection{Proof of \cref{main:E}}
\Cref{main:E} now follows easily from \cref{main:D} the same way as \cref{main:abscissa_P} followed from \cref{main:C} considering the shift in passing from the Poincar\'e series to the representation zeta function. Alternatively one may also deduce \cref{main:E} from \cref{main:abscissa_P} considering the correspondence between kernels of $\cmatrix$ and centralizers in $\alginf$.

\part{Representation zeta function of $\SL{4}{(\completion)}{\permiss}$}
\label{part:sl4}
\section{Hypotheses of \cref{main:D}}
In order to apply \cref{main:D} we first check that its hypothesis are satisfied. First of all we show that, for $2\nmid q$, $\spl{4}{(\completion)}{}$ admits a non-degenerate associative symmetric bilinear form.\par
Let $\ring$ be a ring. In this section and throughout the rest of this work $e_{ij}$ ($i,j\in\lbrace 1,2,3,4\rbrace$) denotes the element of the standard basis of $\gl{4}{(\ring)} = \Mat{4\times 4}{(\ring)}$ having $(i,j)-th$ entry equal to $1$ and zeroes elsewhere. The ring $\ring$ is not evident from this notation but it will each time be clear from the context.
\subsection{Killing form}
\label{sec:killing}
Since $\spl{4}{(\localfrac)}{}$ is a semisimple Lie algebra over a field of characteristic $0$, the restriction $\killing$ of its Killing form to $\spl{4}{(\completion)}{}$ is a natural candidate for the non-degenerate associative bilinear symmetric form required by \cref{main:D}. This form indeed inherits the last three properties form the the Killing form of $\spl{4}{(\localfrac)}{}$; however it might become degenerate, when the residue field characteristic of $\completion$ divides its determinant. We now show that $\killing$ is non-degenerate whenever $2\nmid q$. In order to do so, we compute the matrix of $\killing$ with respect to a basis of $\spl{4}{(\completion)}{}$ and show that its determinant is a power of $2$.\par
First of all, let us fix a basis: for the rest of the paper, let $\basis$ be the $\completion$-basis of $\spl{4}{(\completion)}{}$ comprising
\begin{align*}
h_{12}&=e_{11} - e_{12},&%
h_{23}&= e_{22} - e_{33},&%
h_{34}&=e_{33} - e_{44},\\
e_{12}&,&%
e_{23}&,&%
e_{34}&,\\
e_{13}&,&%
e_{24}&,&%
e_{14}&,\\
f_{21}&=e_{21},&%
f_{32}&=e_{32},&%
f_{43}&=e_{43},\\
f_{31}&=e_{31},&%
f_{42}&=e_{42},&%
f_{41}&=e_{41}.
\end{align*}
Secondly we notice that $\killing(X,Y) = 8\tr(XY)$, as it may be directly computed from the definition of the Killing form on $\spl{4}{(\localfrac)}{}$. The matrix of $\killing$ with respect to $\basis$ is computed using the well known multiplication rules among the elements of the standard basis of $\gl{4}{(\completion)}{}$ (see \cite[Section 1.2]{hum1978introlie} for instance). The result is the following $15\times15$ matrix:
\begin{equation*}
8 \begin{pmatrix}
2 & -1 &   &   &    &   &   &   &    \\
-1 & 2 & -1 &   &   &   &   &   &   \\
  & -1 & 2 &   &   &   &   &   &   \\
  &   &   &   &   &   & 1 &   &   \\
  &   &   &   &   &   &   & \ddots &   \\
  &   &   &   &   &   &   &    & 1 \\
  &   &   & 1 &   &   &   &   &   &  \\
  &   &   &   & \ddots  &   &   &   &   \\
  &   &   &   &   & 1 &   &   &   
\end{pmatrix},
\end{equation*}
which has determinant $8^{15}\cdot 4$. Hence $\killing$ is non-degenerate if $2\nmid q$, which we assume henceforth.\par
For later use, we investigate here the Lie lattice anti-automorphism of $\spl{4}{(\completion)}{}$ induced by $\killing$. Recall that having a non-degenerate symmetric bilinear form $\killing$ that is associative with respect to the Lie bracket implies that there is an isomorphism of $\completion$-Lie lattices
\[
\lambda:\xymatrix@R=3pt{\spl{4}{(\completion)}{} \ar[r]&(\spl{4}{(\completion)}{})^\sharp\\
x\ar@{|->}[r]&(y\mapsto \killing(y,x)).}
\]
Let $\iota$ be the $\completion$-linear isomorphism $\completion^{15}\rightarrow \spl{4}{(\completion)}{}$ defined by fixing the basis $\basis$ and let $\eta$ be its dual.\par
In \cref{sec:quadratic} we defined an anti-automorphism of the Lie lattice $\alginf$ as $\killingiso = \iota\circ\eta\circ \lambda$. It is clear that, for almost all primes, $\killingiso$ may be replaced with a normalization so that $\killingiso ^2 = \id$. In this case (i.e.\ for the specific choice of basis $\basis$) we just need to assume that $2\nmid q$. A straightforward computation shows that the normalized anti-automorphism is given by $\killingiso: x\mapsto \transpose{x}$.
\subsection{Rank-preserving lifts and geometrical smoothness}
We still need to show  that the centralizer-loci of $\spl{4}{(\completion)}{}$  are geometrically smooth. For convenience of notation we identify $\A^{15}(\completion)$ with $\spl{4}{(\completion)}{}$ through the coordinate system 
\[
\iota:\completion^{15}\rightarrow \spl{4}{(\completion)}{}.
\]
Accodingly, for all $2\mu\leq 15$, we identify $\centrvar{\spl{4}{(\completion)}{}}{15 - 2\mu}{(\completion)}$ with its image through $\iota$. Similarly, we identify $\centrvar{\spl{4}{(\completion)}{}}{15 - 2\mu}{(\qfield)}$ with its image through the $\qfield$-linear isomorphism $\qfield^{15}\rightarrow \spl{4}{(\qfield)}{}$ associated with the reduction modulo $\primeideal$ of $\basis$.\par
According to \cref{def:geo_smooth_centr}, showing geometrical smoothness of the centralizer-loci, amounts to showing that, for all $2\mu\leq 15$,
\begin{enumerate}
\item every $\overline{x}\in\centrvar{\spl{4}{(\completion)}{}}{15 - 2\mu}{(\qfield)}$ has a lift to $\centrvar{\spl{4}{(\completion)}{}}{15 - 2\mu}{(\completion)}$,
\item all $x\in\centrvar{\spl{4}{(\completion)}{}}{15 - 2\mu}{(\completion)}$ that reduce to an $\overline{x}\in\centrvar{\spl{4}{(\completion)}{}}{15 - 2\mu}(\qfield)$ modulo $\primeideal$, are such that $[\centr{\spl{4}{(\completion)}{}}{x},\centr{\spl{4}{(\completion)}{}}{x}]$ is an isolated $\completion$-module. 
\end{enumerate}
The first step is to determine the possible Lie centralizer dimensions and consequently the values of $2\mu$ we need to consider.
\subsubsection{Sheets of $\spl{4}{(\C)}{}$}
Fix an embedding of $\completion$ into $\C$. Let $2\mu< 15$. In order to understand the geometry of the $(15 - 2\mu)$-centralizer-locus of $\spl{4}{(\completion)}{}$ we first look at its complexification; i.e.\ we consider 
\[
X_{2\mu} = \centrvar{\spl{4}{(\completion)}{}}{\dimalg - 2\mu}{}{\otimes}_\completion \C.\]
The irreducible components of this complex scheme are the so-called \emph{sheets} of $\liealg = \mathfrak{sl}_{4}{(\C)}^{}$. Every sheet of $\liealg$ corresponds in one-to-one correspondence to a partition of $4$ (see \cite[Section 3.1]{moreau2007dimension}). So  for each partition of $4$, $\partition=[\partel_1,\dots,\partel_f]$ ($\partel_1\geq\dots\geq\partel_f$)  we have a sheet $\sheet_\partition$. The dimension of an orbit $\orbit\subseteq\sheet_\partition$ is given by equation (1) in \cite[Section 3.1]{moreau2007dimension}:
\[
\dim_\C\orbit=2\,\mathrm{m}(\partition),\text{ where }\mathrm{m}(\partition)=(4^2-\sum_{s\in \dualpartition{\partition}} s^2)/2.
\]
and $\dualpartition{\partition}=[s_i\mid i=1,\dots,f]$ ($s_i=\#\lbrace j\mid \partel_j\geq i\rbrace$) is the dual partition of $\partition$. 
A straightforward computation shows that each partition of $4$ gives a different orbit dimension. It follows that the admissible centralizer dimensions in $\spl{4}{(\C)}{}\smallsetminus\lbrace 0\rbrace$ are $9,7,5,3$. The corresponding loci are $X_{2\mu}$, ($2\mu=6,8,10,12$) which coincide with the sheets (see \cref{tab:sheets_4}) and are therefore irreducible.\par
\begin{table}[h]
\caption{The sheets of $\spl{4}{(\C)}$}
\label{tab:sheets_4}
\begin{tabular}{ccc}
\toprule
Partition $\partition$ of $4$		&Orbit dimension $2\,\mathrm{m}(\partition)$	&$\dim_\C (\centrvar{\spl{4}{(\completion)}{}}{15 - 2\,\mathrm{m}(\partition)}{}\otimes_\completion \C)$ \\
\midrule
$[1^4]$						& 0									&$0$\\
$[2,1^2]$						&6									&$7$\\
$[2,2]$						& 8									&$9$\\
$[3,1]$						& 10									&$12$\\
$[4]$							& 12									&$15$\\
\bottomrule
\end{tabular}
\end{table}
\subsubsection{Geometrical smoothness} 
\label{sec:hyp_geo_smooth}
The previous investigation shows that there are $4$ non-trivial loci of constant centralizer dimension in $\spl{4}{(\completion)}{}$ ($\centrvar{\spl{4}{(\completion)}{}}{\dimalg - 2\mu}{(\completion)}$ for $\mu\in\lbrace 3,4,5,6\rbrace$). Fix $\mu\in\lbrace 3,4,5,6\rbrace$ and let $X = \centrvar{\spl{4}{(\completion)}{}}{\dimalg - 2\mu}{}$. As $X\otimes_\completion \C$ coincides with a sheet of $\spl{4}{(\C)}{}$ which is smooth by \cite[Section 3, Korollar 2]{bong1989schichten}. Hence, in order to prove that $X(\completion)$ is geometrically smooth, we may apply the criterion in \cref{prop:geo_sm_crit}. Taking into account the dualization operated in passing from kernels of commutator matrices to centralizers, this amounts to proving that for all $x\in\spl{4}{(\completion)}{}$ reducing modulo $\primeideal$ to a point in $X(\qfield)$,
\[
\dim_{\qfield} [\centr{\spl{4}{(\qfield)}{}}{\reducinf{1}(x)}, \centr{\spl{4}{(\qfield)}{}}{\reducinf{1}(x)}]  = 15 - \dim_\C (X\otimes_\completion \C).
\]
We shall complete this step once we have obtained an overview of the Lie centralizers in $\spl{4}{(\qfield)}{}$.
\section{Classification by Lie centralizers}
\label{sec:cross-sec}
In the present and following sections we seek to define a classification by Lie centralizers that will enable us to compute the formula of \cref{main:D} in a recursive way.
More in detail we give the following definition:
\begin{defn}
\label{def:C_step}
Let $\mathcal{S}=\lbrace \mathbf{s},\mathbf{t} \rbrace$ be a sequence of $\arr$-centralizer classes. Let also $x\in\mathbf{s}$.
We define
\[
\ccoeffiso{\mathcal{S}}{(\spl{4}{(\qfield)}{})}^x=\lbrace (x_1,x_2)\in\ccoeffiso{\mathcal{S}}{(\spl{4}{(\qfield)})}\mid x_1 = x\rbrace.
\]
\end{defn}
Assume that for all $\mathcal{S}$ and $x$ as in the above definition $\lvert \ccoeffiso{\mathcal{S}}{(\spl{4}{(\qfield)}{})}^x\rvert$ does not depend on $x$ but only on $\mathbf{s}$ and $\mathbf{t}$. Let $I$ be as in \cref{main:D} and $\mathcal{T} = \lbrace \mathbf{s}_1,\dots,\mathbf{s}_\lastindex\rbrace$ be an $I$-sequence of $\arr$-centralizer classes. Then it is clear that, by definition of $\ccoeffiso{\mathcal{T}}{(\spl{4}{(\qfield)}{})}$,
\begin{equation}
\label{eq:ind_step}
\lvert\ccoeffiso{\mathcal{T}}{(\spl{4}{(\qfield)}{})}\rvert=\lvert \mathbf{s}_1\rvert\cdot
\prod_{i=1}^{\lastindex-1}
\lvert\ccoeffiso{\lbrace \mathbf{s}_i,\mathbf{s}_{i+1}\rbrace}{(\spl{4}{(\qfield)}{})}^{x_i}\rvert
\end{equation}
where $x_1,\dots,x_\lastindex\in\spl{4}{(\qfield)}{}$ such that 
$\centr{\spl{4}{(\qfield)}{}}{x_j}\in\mathbf{s}_j$ for all $j\in\lbrace 1,\dots,\lastindex\rbrace$. 
\begin{rem}
\label{rem:class_orbits}
It is also clear that, fixed $x$ as above, if $\psi$ is a Lie algebra automorphism of $\spl{4}{(\qfield)}{}$ then 
\[
\lvert\lbrace(\psi(x),\psi(y)) \mid (x,y)\in\ccoeffiso{\mathcal{S}}{(\spl{4}{(\qfield)})}\rbrace\rvert =%
\lvert \ccoeffiso{\mathcal{S}}{(\spl{4}{(\qfield)}{})}^{\psi(x)}\rvert.
\]
Thus, if $\mathbf{s}$ and $\mathbf{t}$ are two orbits for a group action by Lie algebra automorphisms of $\spl{4}{(\qfield)}{}$ then $\lvert \ccoeffiso{\mathcal{S}}{(\spl{4}{(\qfield)}{})}^x\rvert$ is the same for all $x\in\mathbf{s}$.
\end{rem}
In what follows we shall therefore see how to coarsen the classification
\[
\arr = \lbrace \orbit \mid \orbit\text{ an $\SL{4}{(\qfield)}{}$-orbit in $\spl{4}{(\completion)}{}$}\rbrace,
\] 
in order to minimize the number of its elements without compromising the possibility of using \cref{eq:ind_step}.
\subsection{Affine cross-section}
What we need is a way of parameterizing $\SL{4}{(\qfield)}{}$-orbits whose elements have Lie centralizer of the same dimension. Once we have that, coarsening $\arr$ will be just a matter of deciding for which parameters we should combine the corresponding $\SL{4}{(\qfield)}{}$-orbits in a bigger centralizer class or not. As done before we first consider the situation over $\C$.\par
As first proved by D. Peterson \cite[Chapter 3]{pet1978geometry}, every sheet $\sheet$ of $\spl{4}{(\C)}$ contains an \emph{affine cross-section}: a subset of $\sheet$ that meets each $\SL{4}{(\C)}{}$-adjoint orbit exactly once and is isomorphic to an affine space (an explicit construction is described in \cite[Section 1.4]{borho1981}). The set of $\SL{4}{(\C)}{}$-adjoint orbits and the set of $\GL{4}{(\C)}$-conjugation orbits in $\spl{4}{(\C)}{}$ coincide, because they are both in one to one correspondence with the Jordan forms in $\spl{4}{(\C)}{}$. It follows that the affine cross-section also parameterizes the $\GL{4}{(\C)}$-orbits in $\spl{4}{(\C)}$.\par
Let $\qclosure$ be an algebraic closure of $\qfield$. An analog cross-section for the $\SL{4}{(\qclosure)}{}$-orbits in $\spl{4}{(\qclosure)}{}$ has been constructed in \cite[Section 4]{bong1989schichten}. However, the $\qfield$-rational points of an $\SL{4}{(\qclosure)}{}$-orbit $\orbit\subseteq\spl{4}{(\qclosure)}$ might consist of a union of more than one $\SL{4}{(\qfield)}{}$-orbit  and the $\qfield$-rational points on the cross-section might no longer parameterize the $\SL{4}{(\qfield)}{}$-orbits. The solution is to consider the $\GL{4}{(\qfield)}{}$-conjugation instead of the $\SL{4}{(\qfield)}{}$-adjoint action. Indeed, by \cref{rem:class_orbits}, we may replace the classification consisting of $\SL{4}{(\qfield)}{}$-orbits with the classification consisting of $\GL{4}{(\qfield)}{}$-orbits without compromising the applicability of \cref{eq:ind_step}. However, now, a consequence of the Lang-Steinberg Theorem  \cite[Theorem 21.11]{maltes2011linear} guarantees that if $\orbit$ is a $\GL{4}{(\qclosure)}$-orbit in $\spl{4}{(\qclosure)}{}$, then the $\GL{4}{(\qfield)}$-action on the set of $\qfield$-points of $\orbit$ remains transitive. 
The only thing that remains to check, before we are able to employ affine cross-sections defined over $\qfield$, is that the $\qfield$-rational points of an affine cross-section defined over $\qclosure$ still parameterize $\GL{4}{(\qfield)}{}$-orbits in $\spl{4}{(\qfield)}$.
{\begin{prop}
\label{prop:affine_CS}
A $\GL{4}{(\qclosure)}$-orbit contains an $\qfield$-rational point if and only if its intersection with the \emph{affine cross-section} contains an $\qfield$-rational point.
\begin{proof}
Let $\frob$ be the Frobenius automorphism of $\qclosure$. We observe that an orbit containing an $\qfield$-rational point is $\frob$-stable while the affine cross-section is  $\frob$-stable because it is defined by equations with integer coefficients. It follows that their intersection, which consists of a single point, is $\frob$-stable and therefore $\qfield$-rational.
\end{proof}
\end{prop}}
\section{Centralizers over the finite field}
\label{sec:centr}
In this section we use the affine-cross section to set up a classification by centralizers that is fine enough to allow us to use \cref{eq:ind_step}. We operate a case distinction according to the $\qfield$-dimension of the Lie centralizer.
\subsection{Dimension $3$}
An element of $\spl{4}{(\qfield)}{}$ with $3$-dimensional Lie centralizer in $\spl{4}{(\qfield)}{}$ is called regular and its centralizer is called regular too. Since all the regular elements in $\spl{4}{(\qfield)}{}$ have abelian Lie centralizer we may as well collect them all in the same class $\LieReg$. We postpone the problem of computing $\lvert \LieReg\rvert$ to the end of the section as it is best to obtain this number by subtraction from the cardinality of $\spl{4}{(\qfield)}{}$ once the number of non-regular elements is known.
\subsection{General procedure for non-regular elements}
\label{sec:centr_method}
Non-regular elements require a little subtler classification than the one based solely on Lie centralizer dimension  operated for regular elements.\par
Fixed a Lie centralizer dimension, $e$ say, we shall first use the affine cross-section on the corresponding sheet to find all orbits with $e$-dimensional centralizer. Once this is done, we shall then decide which $\GL{4}{(\qfield)}$-orbits will form the centralizer classes. The general strategy is therefore a multi-layered case distinction based in the first place on the dimension of the Lie centralizer, and in the second place on the affine cross-section associated with the chosen dimension. We shall proceed as follows:
\begin{enumerate}
\item  \label{list:centr_1} Find the cross-section using the description in  \cite[Section 1.4]{borho1981}.
\item  \label{list:centr_2} Fix a point $\element$ on the cross-section and compute its Lie centralizer $\centr{\gl{4}{(\qfield)}{}}{\element}$.
\item \label{list:centr_3} Obtain 
\begin{align*}
\centr{\spl{4}{(\qfield)}{}}{\element}	&=\centr{\gl{4}{(\qfield)}{}}{\element}\cap \spl{4}{(\qfield)}{}\\
\centr{\SL{4}{(\qfield)}{}}{\element}	&=\centr{\gl{4}{(\qfield)}{}}{\element}\cap \SL{4}{(\qfield)}{}.
\end{align*}
The results of this step are summarized in \cref{tab:dim_5,tab:dim_7,tab:dim_9,tab:dim_9_nilp}. The first column always contains the conditions on the parameters of the affine cross-section and above them the name we give to the corresponding $\GL{4}{(\qfield)}$-orbit. For the moment we do not check for the applicability of \cref{eq:ind_step} as at a later stage we shall coarsen even further the classification by centralizers obtained at this stage.
\end{enumerate}
We operate by first considering the elements with $5$-dimensional centralizers and then proceeding in order to elements with $7$ and $9$-dimensional Lie centralizers.
\subsection{Dimension $5$}
We consider the affine cross-section on $\sheet_{[3,1]}$:
\begin{equation*}
	\crosssection_{[3,1]}(\parameterOne,\parameterTwo)=	\begin{pmatrix}
			\parameterOne&0&0&0\\
			0&\parameterOne&1&0\\
			0&0&-\parameterOne&1\\
			0&0&\parameterTwo&-\parameterOne\\
	    	\end{pmatrix}.
\end{equation*}
for $\parameterOne,\parameterTwo\in\qfield$.
\begin{defn}
\label{def:jord_form}
Every element  $\element\in\spl{4}{(\qfield)}$ admits a Jordan decomposition. When the semisimple part of $\element$ is diagonalizable over $\qfield$, we say that $\element$ \emph{admits a Jordan normal form} or that its orbit \emph{contains a Jordan normal form}.
\end{defn}
We need a case distinction between orbits that contain a Jordan normal form and orbits that do not contain such a matrix: fix $\parameterOne,\parameterTwo$ and let $\element= C_{[3,1]}(\parameterOne,\parameterTwo)$. A quick computation of the characteristic and minimal polynomials of $\element$ yields:
\begin{align*}
\chi_\element(X)&=(X-\parameterOne)^2 (X^2+2\parameterOne X+\parameterTwo+\parameterOne^2)\\
m_\element(X)&=(X-\parameterOne) (X^2+2\parameterOne X+\parameterTwo+\parameterOne^2).
\end{align*}
From this, we see that $\element$ admits a Jordan normal form if and only if $-\parameterTwo$ is a square, and it is diagonalizable if and only if 
$-\parameterTwo$ is a non-zero square and $4\parameterOne^2\neq -\parameterTwo$. \Cref{tab:dim_5} gathers the relevant information about all isomorphism classes of $5$-dimensional Lie centralizers in $\spl{4}{(\qfield)}{}$. The computations are conducted as explained at the beginning of the section, as an example we show the details of the case $\alpha \neq 0$ and $-\beta$ not a square, in which the calculations are more interesting than in the other cases. The full computations can be found in \cite[Section 5.4]{zor2016thesis}.
\subsubsection{Orbits without Jordan normal form}
We examine the orbits that do not contain the Jordan normal form of the matrix $\element=\cs{[3,1]}{\parameterOne,\parameterTwo}$ on the cross-section. In other words $\chi_\element(X)$ does not split in linear factors with coefficients in $\qfield$; this happens precisely when $-\beta\in\qfield$ is not a square. Since we are classifying elements of $\spl{4}{(\qfield)}{}$ according to their $\GL{4}{(\qfield)}{}$-conjugacy orbit, we may replace $\element$ with its Frobenius normal form
\[
\element=\begin{pmatrix}
\alpha & 0 & 0 & 0 \\
0 & \alpha & 0 & 0 \\
0 & 0 & 0 & - (\alpha^{2} +  \beta) \\
0 & 0 & 1 & -2 \alpha
\end{pmatrix}.
\]
Now let $M=(m_{ij})_{i,j}\in\Mat{4\times 4}{(R)}$. The Lie centralizer  of $\element$ is the set of solutions to the linear system defined by $[\element,M]=0$. Since $-\beta$ is not a square we deduce that $m_{ij}=0$ when $i\leq 2$, $j\geq 3$ and when $i\geq 3$, $j\leq 2$. Thus
\begin{equation*}
\centr{\gl{4}{(\qfield)}}{\element}=\left\lbrace\begin{pmatrix}
										m_{11}&m_{12}&0&0\\
										m_{21} &m_{22}&0&0\\
										0&0&2\alpha\,m_{43}+m_{44}&-(\alpha^2+\beta)\,m_{43}\\
										0&0&m_{43}&m_{44}\\
								                 \end{pmatrix}
                                                     \in\Mat{4\times 4}{(\qfield)}\right\rbrace.
\end{equation*}
Since $-\beta\in\qfield$ is not a square, the matrices 
\[
\begin{pmatrix}
	2\alpha\,m_{43}+m_{44}&-(\alpha^2+\beta)\,m_{43}\\
	m_{43}&m_{44}\\
\end{pmatrix}
\]
with $m_{43},m_{44}\in\qfield$ form a Lie algebra isomorphic to $\F_{q^2}$, and therefore $\centr{\GL{4}{(\qfield)}}{\element}$ is isomorphic to $\GL{2}{(\qfield)}\times\,{\F_{q^2}}^\times$. It follows
\begin{align*}
\vert\centr{\GL{4}{(\qfield)}}{\element}\vert		&=\vert \GL{2}{(\qfield)}\vert\cdot(q^2-1)=(q-1)^3(q+1)^2q\\
\centr{\SL{4}{(\qfield)}{}}{\element}			&\cong \SL{2}{(\qfield)}{}\times\, {\F_{q^2}}^\times.
\end{align*}
\setlength{\tabcolsep}{3pt}
\begin{table}[h]
\caption{Centralizers of $a=C_{[3,1]}(\alpha,\beta)$ with their structure}
 \label{tab:dim_5}
\begin{tabular}{llll}
\toprule
Class							&Basis of $\centr{\spl{4}{(\qfield)}{}}{a}$		&Commutators				&$\centr{\SL{4}{(\qfield)}{}}{a}\cong$\\
$\alpha, \beta$						&$\lbrace c_0,\dots, c_4\rbrace$				&$\neq 0$									&						\\
\midrule%
$\LieThreeOneNilp$					&$e_{14}$								&$[c_0,c_1] = -c_2$						&$(\heis{\qfield}\times{\qfield}^{+})\rtimes{\qfield}^{\times}$\\
$\alpha = 0$						&$ e_{21}$								&$[c_0,c_4] = -4c_0$					&\\
$\beta=0$							&$e_{24}$								&$[c_1,c_4] = -4c_1$					&\\
								&$e_{23} + e_{34}$ 							&									&\\
								&$-3e_{11} + e_{22} + e_{33} + e_{44}$			&									&\\			
\midrule

$\LieThreeOneTwoEV_{\alpha,\beta}$				&$e_{12}$								&$[c_0,c_1] = c_2$						&$\GL{2}{(\qfield)}{}\times{\qfield}^{+}$\\
 $\alpha\neq 0$			&$e_{21}$								&$[c_2,c_0] = 2 c_0$								&\\
 $\beta = 0$						&$e_{11} - e_{22}$							&$[c_2,c_1] = -2 c_1$					&\\
 								&$e_{11} + e_{22} - e_{33} - e_{44}$				&									&\\
 								&$e_{34}$								&									&\\
\midrule
$\LieThreeOneThreeEV_{\alpha,\beta}$				&$-3e_{11} + e_{22} + e_{33} + e_{44}$			&$[c_2,c_3] = c_1$						&$\heis{\qfield}\rtimes({\qfield}^\times\times{\qfield}^\times)$\\
 $\alpha\neq 0$			&$e_{23}$								&$[c_2,c_4] = c_2$									&\\
 $\beta = - 4\alpha^2$ 				&$e_{24}$								&$[c_3,c_4] = - c_3$						&\\
 								&$e_{43}$								&									&\\
 								&$e_{44}-e_{11}$							&									&\\
\midrule
$\LieThreeOneDiag_{\alpha,\beta}$					&$e_{12}$								&$[c_0,c_1] = c_2$						&$\GL{2}{(\qfield)}\times{\qfield}^\times$\\
 $\beta = -\gamma$		&$e_{21}$								&$[c_2,c_0] = 2 c_0$								&					\\
$\gamma\in\qfield^{\times}$			&$e_{11} - e_{22}$							&$[c_2,c_1] = -2 c_1$					&					\\
$\alpha\neq 0$						&$e_{11} + e_{22} - e_{33} - e_{44}$				&									&					\\
$\beta \neq -4\alpha^2$ 				&$e_{33} - e_{44}$							&									&\\

\midrule
$\LieThreeOneNonJor_{\alpha,\beta}$				&$e_{12}$						&$[c_0,c_1] = c_2$								&$\SL{2}{(\qfield)}{}\times\, {\F_{q^2}}^\times$\\
 $-\beta\neq \gamma^2$	&$e_{21}$						&$[c_2,c_0] = 2 c_0$										&\\
$\forall\gamma\in\qfield$				&$e_{11} - e_{22}$					&$[c_2,c_1] = -2 c_1$							&\\
$\alpha \neq 0$						&$e_{11} + e_{22} - e_{33} - e_{44}$		&											&\\
$\delta = (\alpha^2 +\beta)$ 			&$\alpha (e_{33} - e_{44}) - \delta e_{34} + e_{43}$	&									&\\
\bottomrule
\end{tabular}
\end{table}
\subsection{Dimension $7$}
The affine cross-section in $\sheet_{[2,2]}$ is one-dimensional. The following is a parameterization of it in terms of $\parameterOne\in\qfield$:
\begin{equation*}
	\crosssection_{[2,2]}(\parameterOne)=\begin{pmatrix}
			0&1&0&0\\
			\parameterOne&0&0&0\\
			0&0&0&1\\
			0&0&\parameterOne&0\\
	    	\end{pmatrix}.
\end{equation*}
\Cref{tab:dim_7} gives an overview of the Lie and group centralizers of elements on $\crosssection_{[2,2]}(\parameterOne)$, the computations have been conducted as explained in \cref{sec:centr_method}. We show how the second and third case are deduced, for the Lie algebras in the table might look the same at first glance. The computations for the other cases are found in \cite[Section 5.5]{zor2016thesis}.
\subsubsection{Semisimple orbits}
When $\alpha$ is a non-zero square, the corresponding point on the cross-section is a diagonalizable element of $\spl{4}{(\qfield)}{}$. Indeed, in this case $\element$ is similar to 
\begin{equation*}
	       	\begin{pmatrix}
			\gamma&0&0&0\\
			0&\gamma&0&0\\
			0&0&-\gamma&0\\
			0&0&0&-\gamma\\
	    	\end{pmatrix},
\end{equation*}
where  $\gamma\in\qfield^{\times}$ such that $\gamma^2 = \alpha$. So in reality
\begin{equation}
\label{eq:centr_22Sem}
\begin{split}
\centr{\spl{4}{(\qfield)}{}}{\element}	&\cong\spl{2}{(\qfield)}{} + \spl{2}{(\qfield)}{} + \qfield \\
\centr{\SL{4}{(\qfield)}{}}{\element}	&\cong\lbrace\GL{2}{(\qfield)}\times\GL{2}{(\qfield)}\mid \det A = \det B^{-1}\rbrace.
\end{split}
\end{equation}
When $\alpha$ is not a square in $\qfield$ however the minimal polynomial of $\element$ does not split into linear factors and hence $\element$ is not diagonalizable. Its Lie centralizer in $\gl{4}{(\qfield)}{}$ is 
\begin{equation*}
\centr{\gl{4}{(\qfield)}}{\element} = \left\lbrace \begin{pmatrix}
										m_{11}&m_{12}&m_{13}&m_{14}\\
										\parameterOne\,m_{12} &m_{11}&\parameterOne\,m_{14}&m_{13}\\
										m_{31}&m_{32}&m_{33}&m_{34}\\
										\parameterOne\,m_{32}&m_{31}&\parameterOne\,m_{34}&m_{33}\\
								                 \end{pmatrix}\in\Mat{4 \times 4}{(\qfield)}\right\rbrace.
\end{equation*}
We notice that the matrices above all consist of four blocks in 
\begin{equation*}
R=\left\lbrace\left.	\begin{pmatrix}
				x&y\\
				\parameterOne\,y &x
			\end{pmatrix}\,\right\vert\, x,y\in\qfield\right\rbrace\cong\qfield(\gamma)\cong\F_{q^2}
\end{equation*}
where $\gamma^2=\parameterOne$. Now, let $N:\F_{q^2}\rightarrow\qfield$ be the norm function defined by $(x,y)\mapsto x^2-\parameterOne y^2$, $\det_{\F_{q^2}}$ and $\det_{\qfield}$ be the determinant function on  $\Mat{2\times 2}{(\F_{q^2})}$ and $\Mat{4 \times 4}{(\F_{q})}$, respectively. If $M \in \centr{\gl{4}{(\qfield)}}{\element}$, its determinant is given by
\begin{equation*}
{\det}_{\qfield}(M)=\mathrm{N}( {\det}_{\F_{q^2}}(M)),
\end{equation*}
where $\mathrm{N}$ is the norm function on $\F_{q^2}$. Since the elements of norm $1$ in $\F_{q^2}$ form a cyclic group of order $q+1$, %
we conclude that
\begin{equation}
\begin{split}
\label{eq:centr_22Non}
\centr{\spl{4}{(\qfield)}{}}{\element}	&\cong \qfield + \spl{2}{(\F_{q^2})}{}\\
\centr{\SL{4}{(\qfield)}{}}{\element} 	&\cong C_{q+1}\times \SL{2}{(\F_{q^2})}{}.
\end{split}
\end{equation}
\begin{table}[h]
\caption{Centralizers of $a=C_{[2,2]}(\alpha)$ with their structure}
 \label{tab:dim_7}
\begin{tabular}{llll}
\toprule
Iso-class				&Basis of $\centr{\spl{4}{(\qfield)}{}}{a}$		&Commutators							&$\centr{\SL{4}{(\qfield)}{}}{a}$\\
 $\alpha$	&$\lbrace c_0,\dots, c_6\rbrace$				&$\neq 0$									&isomorphic to						\\
\midrule%
$\LieTwoTwoNilp$					&$e_{13} + e_{24}$							&$[c_0,c_1] = c_2$							&$ ({\qfield}^{+})^4\rtimes\SL{2}{(\qfield)}{}$\\
 $\alpha = 0$			&$e_{31} + e_{42}$							&$[c_0,c_2] = -2 c_0$						&						\\
								&$e_{11} + e_{22} - e_{33} -e_{44}$				&$[c_1,c_2] = 2 c_1$						&						\\
								&$e_{12} - e_{34}$							&$[c_0,c_3] = -2c_4$						&						\\
								&$e_{14}$								&$[c_0,c_5] = c_3$							&						\\
								&$e_{32}$								&$[c_1,c_3] = 2c_5$							&						\\
								&$e_{12} + e_{34}$							&$[c_1,c_4] = -c_3$							&						\\
								&										&$[c_2,c_4] = 2c_4$							&						\\
								&										&$[c_2,c_5] = -2c_5$						&						\\
\midrule
$\LieTwoTwoSemDiag_{\alpha}$				&$e_{12} + e_{34} + \alpha (e_{21} + e_{42})$		&$[c_1,c_3] = 2 c_4$						&$\lbrace\GL{2}{(\qfield)}\times\GL{2}{(\qfield)}$\\
 $\alpha = \gamma^2$		&$e_{12} - e_{34} + \alpha (e_{21} - e_{42})$		&$[c_2,c_4] = 2c_4$							&$\mid \det A = \det B^{-1}\rbrace$			\\
$\gamma \in \qfield^{\times}$			&$e_{11} + e_{22} - e_{33} -e_{44}$				&$[c_1,c_4] = 2\alpha c_3$					&						\\
								&$e_{13} + e_{24}$							&$[c_2,c_3] = 2 c_3$						&						\\
								&$e_{41} + \alpha e_{23}$					&$[c_1,c_5] = -2c_6$						&						\\
								&$e_{31} + e_{42}$							&$[c_2,c_6] = -2c_6$						&						\\
								&$e_{32} + \alpha e_{41}$					&$[c_1,c_6] = -2\alpha c_5$					&						\\
								&										&$[c_2,c_5] = -2c_5$						&						\\
								&										&$[c_3,c_5] = c_2$							&						\\
								&										&$[c_4,c_6] = \alpha c_2$						&						\\
								&										&$[c_3,c_6] = c_1$							&						\\
								&										&$[c_4,c_5] = c_1$							&						\\
\midrule
$\LieTwoTwoSemNon_{\alpha}$				&$e_{12} + e_{34} + \alpha (e_{21} + e_{42})$		&$[c_1,c_3] = 2 c_4$						&$C_{q+1}\times \SL{2}{(\F_{q^2})}{}$\\
 $\alpha \neq \gamma^2$	&$e_{12} - e_{34} + \alpha (e_{21} - e_{42})$		&$[c_2,c_4] = 2c_4$							&						\\
$\forall \gamma \in \qfield$			&$e_{11} + e_{22} - e_{33} -e_{44}$				&$[c_1,c_4] = 2\alpha c_3$					&						\\
								&$e_{13} + e_{24}$							&$[c_2,c_3] = 2 c_3$						&						\\
								&$e_{41} + \alpha e_{23}$					&$[c_1,c_5] = -2c_6$						&						\\
								&$e_{31} + e_{42}$							&$[c_2,c_6] = -2c_6$						&						\\
								&$e_{32} + \alpha e_{41}$					&$[c_1,c_6] = -2\alpha c_5$					&						\\
								&										&$[c_2,c_5] = -2c_5$						&						\\
								&										&$[c_3,c_5] = c_2$							&						\\
								&										&$[c_4,c_6] = \alpha c_2$						&						\\
								&										&$[c_3,c_6] = c_1$							&						\\
								&										&$[c_4,c_5] = c_1$							&						\\
\bottomrule
\end{tabular}
\end{table}
\subsection{Dimension 9}
\label{sec:centr_dim_9}
The affine cross-section in $\sheet_{[2,1^2]}$ is one-dimensional and, for $\parameterOne\in\qfield$, this is its parameterization in the affine space $\spl{4}{(\qfield)}$:
\begin{equation*}
	\crosssection_{[2,1^2]}(\parameterOne)=	\begin{pmatrix}
			3\parameterOne&1&0&0\\
			0&-\parameterOne&0&0\\
			0&0&-\parameterOne&0\\
			0&0&0&-\parameterOne\\
	    	\end{pmatrix}.
\end{equation*}
\subsubsection{Semisimple elements}
When $\alpha\neq 0$ the corresponding matrix in the cross-section is semisimple diagonalizable with three coincident eigenvalues as one can see by computing its characteristic and minimal polynomial. Let $\alpha \neq 0$ and $\element = 	\crosssection_{[2,1^2]}(\parameterOne)$. By what we just discussed, the $\GL{4}{(\qfield)}$-orbit of $\element$ contains
\[
\begin{pmatrix}
-3\alpha	&0		&0		&0\\
0		&\alpha	&0		&0\\
0		&0		&\alpha	&0\\
0		&0		&0		&\alpha\\
\end{pmatrix}
\]
A quick computation of the Lie centralizer in $\gl{4}{(\qfield)}$ of this last matrix yields that
\begin{align*}
\centr{\spl{4}{(\qfield)}{}}{a}	&\cong	\gl{3}{(\qfield)}\\
\centr{\SL{4}{(\qfield)}{}}{a}	&\cong	\GL{3}{(\qfield)}
\end{align*}
\setlength{\tabcolsep}{6pt}
\begin{table}[h]
\caption{Centralizers of $a=C_{[2,1,1]}(\alpha)$ $\alpha\neq 0$}
 \label{tab:dim_9}
\begin{tabular}{lll}
\toprule
Iso-class											&$\centr{\spl{4}{(\qfield)}{}}{a}$				&$\centr{\SL{4}{(\qfield)}{}}{a}$\\
$\alpha$											& isomorphic to							&isomorphic to\\	
\midrule
$\LieTwoOneOneSem_{\alpha}$						&$\gl{3}{(\qfield)}$							&$\GL{3}{(\qfield)}$\\
$\alpha \neq 0$								&										&\\
\bottomrule
\end{tabular}
\end{table}
\subsubsection{Nilpotent elements}
It only remains to determine the Lie and group centralizers of the nilpotent element
\[
	\element=\begin{pmatrix}
		0&1&0&0\\
		0&0&0&0\\
		0&0&0&0\\
		0&0&0&0\\
	     \end{pmatrix}.
\]
The Lie centralizer in $\gl{4}{(\qfield)}$ of the matrix above is the following:
\begin{equation*}
\centr{\gl{4}{(\qfield)}}{\element}=\left\lbrace \begin{pmatrix}
										m_{11}&m_{12}&m_{13}&m_{14}\\
										0&m_{11}&0&0\\
										0&m_{32}&m_{33}&m_{34}\\
										0&m_{31}&m_{43}&m_{44}\\
								                 \end{pmatrix}
                                                     \in\gl{4}{(\qfield)}\right\rbrace.
\end{equation*}
The Lie centralizer in $\spl{4}{(\qfield)}{}$ is determined by simply taking the traceless matrices in $\centr{\gl{4}{(\qfield)}}{\element}$. We fix the basis for $\centr{\spl{4}{(\qfield)}{}}{\element}$ reported in \cref{tab:dim_9_nilp} and compute its commutator relations using the well known rules to compute the Lie bracket of two elements of the standard basis of $\gl{4}{(\qfield)}$. Looking at these commutator relations we notice that the subgroup $H\leq\centr{\SL{4}{(\qfield)}{}}{a}$ defined by
\[
\left\lbrace\left. M=\begin{pmatrix}
										1&m_{12}&m_{13}&m_{14}\\
										0&1&0&0\\
										0&m_{32}&1&0\\
										0&m_{31}&0&1\\
								                 \end{pmatrix}
                                                     \,\right\vert\, M\in\GL{4}{(\qfield)}\right\rbrace
\]
is isomorphic to the direct product $\heis{\qfield}\curlyvee\heis\qfield$ of two copies of the Heisenberg group $\heis{\qfield}$ with amalgamation in the centre. Furthermore  $\centr{\SL{4}{(\qfield)}{}}{a}=HS$ where
\begin{gather*}
S=\left\lbrace\left. M=\begin{pmatrix}
										m_{11}&0&0&0\\
										0&m_{11}&0&0\\
										0&0&m_{33}&m_{34}\\
										0&0&m_{43}&m_{44}\\
								                 \end{pmatrix}
								                  \,\right\vert\, M\in\SL{4}{(\qfield)}{}\right\rbrace\cong\GL{2}{(\qfield)}.
\end{gather*}
As a consequence
\begin{equation*}
\centr{\SL{4}{(\qfield)}{}}{a}\cong (\heis{\qfield}\curlyvee\heis{\qfield})\rtimes\GL{2}{(\qfield)}{}
\end{equation*}
where the semidirect product is defined by the commutator relations in \cref{tab:dim_9_nilp}. This concludes the analysis on the structure of Lie centralizers of $\spl{4}{(\qfield)}{}\smallsetminus \lbrace 0 \rbrace$. The commutator relation in the third column of \cref{tab:dim_5,tab:dim_7,tab:dim_9,tab:dim_9_nilp} also ensure that for $2\nmid q$ the hypotheses of \cref{main:D} are fully satisfied according to the discussion in \cref{sec:hyp_geo_smooth} as the geometrical smoothness criterion in \cref{prop:geo_sm_crit} applies.
\setlength{\tabcolsep}{3pt}
\begin{table}[h]
\caption{Centralizer of $a=C_{[2,1,1]}(0)$ with its structure}
 \label{tab:dim_9_nilp}
\begin{tabular}{llll}
\toprule
Iso-class							&Basis of $\centr{\spl{4}{(\qfield)}{}}{a}$		&Commutators								&$\centr{\SL{4}{(\qfield)}{}}{a}$\\
								&$\lbrace c_0,\dots, c_8\rbrace$				&$\neq 0$									&isomorphic to				\\
\midrule%
$\LieTwoOneOneNilp$				&$e_{13}$								&$[c_0,c_2] = c_4$							&$(\heis{\qfield}\curlyvee\heis{\qfield})\rtimes\GL{2}{(\qfield)}{}$	\\
($\alpha = 0$)						&$e_{14}$								&$[c_0,c_6] = c_1$							&						\\
								&$e_{32}$								&$[c_0,c_7] = c_0$							&						\\
								&$e_{42}$								&$[c_0,c_8] = -2 c_ 0$						&						\\
								&$e_{12}$								&$[c_1,c_ 3] = c_4$							&						\\
								&$e_{43}$								&$[c_1,c_5] = c_0$							&						\\
								&$e_{34}$								&$[c_1,c_7] = -c_1$							&						\\
								&$e_{33} - e_{44}$							&$[c_1,c_ 8] = -2 c_1$						&						\\
								&$e_{11} + e_{22} - e_{33} -e_{44}$				&$[c_2,c_5] = - c_3$							&						\\
								&										&$[c_2,c_7] = - c_2$							&						\\
								&										&$[c_2,c_8] = 2 c_2$						&						\\
								&										&$[c_3,c_6] = - c_2$							&						\\
								&										&$[c_3,c_7] = c_3$							&						\\
								&										&$[c_3,c_8] = 2 c_3$						&						\\
								&										&$[c_5,c_6] = - c_7$							&						\\
								&										&$[c_5,c_ 7] = 2 c_5$						&						\\
								&										&$[c_6,c_7] = - 2 c_6$						&						\\
\bottomrule
\end{tabular}
\end{table}
\subsection{Subregular elements}
\label{sec:subreg}
We call $x\in\spl{4}{(\qfield)}{}$ \emph{subregular} when $\centr{\spl{4}{(\qfield)}{}}{x}$ has dimension $5$. Accordingly, Lie centralizers of subregular elements are called subregular too. When a centralizer class consists of subregular element, it is called a subregular centralizer class. \Cref{tab:dim_7} gives an overview of the subregular centralizer classes established so far. As it will be a crucial concept in the computations that follows, we give the following definition:
\begin{defn}
Let $\algfin$ be an $\qfield$-Lie algebra and let $\otherbasis$ be an $\qfield$-basis for $\algfin$. We may define a matrix of $\qfield$-linear forms encoding the structure constants of $\otherbasis$ in the same way as we defined commutator matrices of $\completion$-Lie lattices (cf.\ \cref{eq:comm_matrix}). We call this matrix the $\qfield$-commutator matrix of $\algfin$ with respect to $\otherbasis$.
\end{defn}
A straightforward application of \cref{lem:jp_rp_lifts,lem:theta} gives the following:  if $x$ is a subregular element with Lie centralizer $\mathfrak{h}$ and lying in the centralizer class $\mathbf{s}$, then $\lvert \ccoeffiso{\lbrace \mathbf{s},\LieReg\rbrace}{(\spl{4}{(\qfield)}{})}^x\rvert$ is equal to 
\begin{equation}
\label{eq:cardinalities}
\frac{\lvert L_2\rvert}{\lvert L_0\rvert }
\end{equation}
where $L_i$  ($i = 2, 0$) is the rank-$i$ locus of the $\qfield$-commutator matrix of $\mathfrak{h}$ with respect to one of its bases (an analogous of \cref{lem:basis_change} shows that this does not depend on the basis chosen for the Lie algebra).\par
The structure constants highlighted in the third column of \cref{tab:dim_5} give that, for  $2\nmid q$, the quantity in \cref{eq:cardinalities} are always $q^3 - 1$. For instance, this is how one shows this for $x = C_{[3,1]}(0,0)$: the basis in \cref{tab:dim_5} gives $\qfield$-commutator matrix
\[
\Cmatrix{\LieThreeOneNilp}{Y_0,\dots,Y_4}=\begin{pmatrix}
0&-Y_2&0&0&4 Y_0\\
Y_2&0&0&0&-4 Y_1\\
0&0&0&0&0\\
0&0&0&0&0\\
-4Y_0&4Y_1&0&0&0
\end{pmatrix}.
\]
Which is readily seen to have rank $2$ if and only if at least one of $Y_0,Y_1,Y_2$ is non-zero and rank $0$ otherwise.\par
The computations for the other regular elements are analogous. First one reduces the problem to considering elements on the affine cross-section. Namely, given a subregular element $x$, let $\element_x$ be the intersection between the $\GL{4}{(\qfield)}$-orbit of $x$ and the cross-section of the centralizer-locus on which $x$ is located. Indeed $\centr{\spl{4}{(\qfield)}{}}{\element_x}$ is a $\GL{4}{(\qfield)}$-conjugate of $\mathfrak{h} = \centr{\spl{4}{(\qfield)}{}}{x}$; say $g^{-1} \mathfrak{h} g = \centr{\spl{4}{(\qfield)}{}}{\element_x}$ for $g\in\GL{4}{(\qfield)}$. If $\basis_{\element_x}$ is the basis  of $\centr{\spl{4}{(\qfield)}{}}{\element_x}$ in \cref{tab:dim_5}, then the basis of $\mathfrak{h}$ defined by conjugating $\basis_{\element_x}$ by $g$ has the same structure constants as $\basis_{\element_x}$. So the commutator matrices are equal. This shows that all the possible commutator matrices of subregular Lie centralizers are found by looking at $\centr{\spl{4}{(\completion)}{}}{\element}$ where $\element$ is on the cross-section. In conclusion, considering the the structure constants in the third column of \cref{tab:dim_5}, one finds that the rank-$2$ loci of all commutator matrices of subregular Lie centralizers have the same cardinality for $2\nmid q$.\par
In addition to the observation above, one may also observe that the derived subalgebra of a subregular Lie centralizer is always $3$-dimensional. Hence, setting
\begin{equation*}
\LieSubreg 		=	 \lbrace x\in\spl{4}{(\qfield)}{}\mid \dim_{\qfield} \centr{\spl{4}{(\qfield)}{}}{x} = 5\rbrace
\end{equation*}
defines a Lie centralizer class and we have just shown that for all $x\in\LieSubreg$,
\[
\lvert \ccoeffiso{\lbrace \LieSubreg,\LieReg\rbrace}{(\spl{4}{(\qfield)}{})}^x\rvert = q^3 -1.
\] Next we consider elements with $7$-dimensional Lie centralizers.
\subsection{Dimension $7$}
\label{sec:centr_dim_7}
Let $\mathbf{s}$ be a $\GL{4}{(\qfield)}$-orbit with $7$-dimensional Lie centralizer. Let $x\in \mathbf{s}$ and let $\mathfrak{s}$ be its Lie centralizer. By a similar argument to the one in \cref{sec:subreg} one reduces the computation of $\lvert \ccoeffiso{\lbrace \mathbf{s},\LieReg\rbrace}{(\spl{4}{(\qfield)}{})}^x\rvert$ and $\lvert \ccoeffiso{\lbrace \mathbf{s},\LieSubreg\rbrace}{(\spl{4}{(\qfield)}{})}^x\rvert$ to the case of $x$ being on the affine cross section $C_{[2,2]}$. Having grouped all the regular and subregular elements into $\LieReg$ and $\LieSubreg$ it is clear that we are essentially trying to count how many $y\in\mathfrak{s}$ make $x + \transpose{y}$ regular or subregular respectively. Again this is done by looking at the rank-loci of the commutator matrix of $\mathfrak{r}$ with respect to the basis in \cref{tab:dim_5}.
\subsubsection{Nilpotent element}
Let $x\in\LieTwoTwoNilp$. The $\qfield$-commutator matrix of $\mathfrak{s}$ in this case is 
\[
\Cmatrix{\LieTwoTwoNilp}{Y_0,\dots, Y_6}=\begin{pmatrix}
0 & Y_{2} & -2 Y_{0} & -2 Y_{4} & 0 & Y_{3}&0 \\
- Y_{2} & 0 & 2 Y_{1} & 2 Y_{5} & - Y_{3} & 0&0 \\
2 Y_{0} & -2 Y_{1} & 0 & 0 & 2 Y_{4} & -2 Y_{5}&0 \\
2 Y_{4} & -2 Y_{5} & 0 & 0 & 0 & 0&0 \\
0 & Y_{3} & -2 Y_{4} & 0 & 0 & 0&0 \\
- Y_{3} & 0 & 2 Y_{5} & 0 & 0 & 0&0\\
0&0&0&0&0&0&0
\end{pmatrix}.
\]
An easy computation shows that this matrix has rank $2$ if and only if at least one of $Y_0, Y_1, Y_2$ is non-zero and  $Y_3= Y_4 = Y_5 = 0$. The rank-$0$ locus has equations  $Y_0=\dots= Y_5 = 0$, and the rank-4 locus is all the rest. By \cref{lem:jp_rp_lifts,lem:theta} it follows that for $x\in\LieTwoTwoNilp$
\begin{equation}
\label{eq:lift_22Nilp}
\begin{split}
\lvert \ccoeffiso{\lbrace \LieTwoTwoNilp,\LieReg\rbrace}{(\spl{4}{(\qfield)}{})}^x\rvert	& = q^3 (q^3 - 1)\\
\lvert \ccoeffiso{\lbrace \LieTwoTwoNilp,\LieSubreg\rbrace}{(\spl{4}{(\qfield)}{})}^x\rvert	& = q^3 -1.
\end{split}
\end{equation}
\subsubsection{Semisimple elements}
Let $\alpha$ be a square in $\qfield^\times$ and let $x\in\LieTwoTwoSemDiag_{\alpha}$. Since $\mathfrak{s}\cong\spl{2}{(\qfield)}{} + \spl{2}{(\qfield)}{} + \qfield$ as highlighted in \cref{eq:centr_22Sem}, there is a basis of $\mathfrak{s}$ such that the associated $\qfield$-commutator matrix is in block diagonal form with two $3\times 3$ blocks consisting of two commutator matrices of $\spl{2}{(\qfield)}{}$ with respect to two independent $\spl{2}{}{}$-triples and a zero $1\times 1$ block. It follows that 
\begin{equation*}
\begin{split}
\lvert \ccoeffiso{\lbrace \LieTwoTwoSemDiag_{\alpha},\LieReg\rbrace}{(\spl{4}{(\qfield)}{})}^x\rvert	& =  (q - 1)^2 (q^2 + q + 1)^2\\
\lvert \ccoeffiso{\lbrace \LieTwoTwoSemDiag_{\alpha},\LieSubreg\rbrace}{(\spl{4}{(\qfield)}{})}^x\rvert	& = 2  (q^3 -1).
\end{split}
\end{equation*}
Setting $\LieTwoTwoSemDiag = \cup_{\gamma\in\qfield^{\times}} \LieTwoTwoSemDiag_{\alpha}$ defines a Lie centralizer class and for $x\in\LieTwoTwoSemDiag$
\begin{equation}
\label{eq:lift_22Sem}
\begin{split}
\lvert \ccoeffiso{\lbrace \LieTwoTwoSemDiag,\LieReg\rbrace}{(\spl{4}{(\qfield)}{})}^x\rvert	& =  (q - 1)^2 (q^2 + q + 1)^2\\
\lvert \ccoeffiso{\lbrace \LieTwoTwoSemDiag,\LieSubreg\rbrace}{(\spl{4}{(\qfield)}{})}^x\rvert	& = 2  (q^3 -1).
\end{split}
\end{equation}
Let finally be a non-square in $\qfield$ and  let $x\in\LieTwoTwoSemNon_\alpha$. Since $\mathfrak{s}\cong \qfield + \spl{2}{(\F_{q^2})}{}$ (see \cref{eq:centr_22Non}), the commutator matrix of $\mathfrak{s}$ with respect to the basis in \cref{tab:dim_7} may only have rank $0$ or $4$. Moreover the rank-$0$ locus is isomorphic to an affine line over $\qfield$. In conclusion, setting 
\[
\LieTwoTwoSemNon = \cup_{\lbrace \alpha\in\qfield\mid \alpha \neq \gamma^2\,\forall \gamma\in\qfield\rbrace} \LieTwoTwoSemNon_\alpha,
\]
defines a Lie centralizer class and for $x\in  \LieTwoTwoSemNon$
\begin{equation}
\label{eq:lift_22SemNon}
\lvert \ccoeffiso{\lbrace \LieTwoTwoSemNon,\LieReg\rbrace}{(\spl{4}{(\qfield)}{})}^x\rvert	= q^6 - 1.
\end{equation}
\subsection{Dimension $9$}
We now consider elements with $9$-dimensional Lie centralizer. Again, by an argument akin to the one in \cref{sec:subreg} we may just consider elements on the cross section $C_{[2,1,1]}$. Furthermore, the computation in \cref{sec:centr_dim_9} shows that for $\element = C_{[2,1,1]}(\alpha)$ with $\alpha \in \qfield^{\times}$, the Lie centralizer $\centr{\spl{4}{(\qfield)}{}}{\element}$ admits a basis $\lbrace c_0,\dots,c_8\rbrace$ such that $\linspan(c_0,\dots,c_7)\cong \spl{3}{(\qfield)}$ and $c_8$ commutes with all the other basis elements. It follows that the relative $\qfield$-commutator matrix, say $\cmatrix_\LieTwoOneOneSem$, consists of a zero row and column next to an $8\times 8$ block on the diagonal that is the $\qfield$-commutator matrix $\cmatrix_{\spl{3}{}{}}$ relative to a basis of $\spl{3}{(\qfield)}{}$. As this is the case for all $\alpha\neq 0$ we may define a Lie centralizer class as
\[
\LieTwoOneOneSem = \cup_{\alpha\in\qfield^{\times}} \LieTwoOneOneSem_\alpha.
\]
The cardinalities of the rank-loci of $\cmatrix_{\spl{3}{}{}}$ were computed by Avni Klopsch Onn and Voll in \cite[Section 6.1]{akov2013representation}, in this setting (applying \cref{lem:jp_rp_lifts,lem:theta} to $\element$) their computations give that for $x\in \LieTwoOneOneSem$
\begin{equation*}
\begin{split}
\lvert \ccoeffiso{\lbrace \LieTwoOneOneSem,\LieReg\rbrace}{(\spl{4}{(\qfield)}{})}^x \rvert	& = q (q - 1) (q^6 +q^5 +q^4 -q^2 -2 q - 1)\\
\lvert \ccoeffiso{\lbrace \LieTwoOneOneSem,\LieSubreg\rbrace}{(\spl{4}{(\qfield)}{})}^x \rvert	& = q^5 + q^4 + q^3 - q^2 - q - 1.
\end{split}
\end{equation*}
Collecting all the centralizer classes defined so far we define the following classification by Lie centralizers:
\[
\arr = \lbrace \LieTwoOneOneSem, \LieTwoOneOneNilp, \LieTwoTwoSemDiag, \LieTwoTwoSemNon, \LieTwoTwoNilp, \LieSubreg, \LieReg\rbrace.
\]
\subsubsection{Nilpotent element}
Let $\element = C_{[2,1,1]}(0)$. Similar to \cref{sec:centr_dim_7} it remains to compute $\lvert\ccoeffiso{\lbrace \LieTwoOneOneNilp,\mathbf{s}\rbrace}{(\spl{4}{(\qfield)}{})}^\element\rvert$ for all $\mathbf{s}\in\arr\smallsetminus \lbrace \LieTwoOneOneSem,\LieTwoOneOneNilp\rbrace$. First of all let us look at the $\qfield$-commutator matrix of $\mathfrak{h} = \centr{\spl{4}{(\qfield)}{}}{\element}$ with respect to the basis in \cref{tab:dim_9_nilp}:
\begin{equation*}
\Cmatrix{\LieTwoOneOneNilp}{Y_1,\dots,Y_8}=\begin{pmatrix}
0 & 0 & Y_{4} & 0 & 0 & 0 & Y_{1} & Y_{0} & -2 Y_{0} \\
0 & 0 & 0 & Y_{4} & 0 & Y_{0} & 0 & - Y_{1} & -2 Y_{1} \\
- Y_{4} & 0 & 0 & 0 & 0 & - Y_{3} & 0 & - Y_{2} & 2 Y_{2} \\
0 & - Y_{4} & 0 & 0 & 0 & 0 & - Y_{2} & Y_{3} & 2 Y_{3} \\
0 & 0 & 0 & 0 & 0 & 0 & 0 & 0 & 0 \\
0 & - Y_{0} & Y_{3} & 0 & 0 & 0 & - Y_{7} & 2 Y_{5} & 0 \\
- Y_{1} & 0 & 0 & Y_{2} & 0 & Y_{7} & 0 & -2 Y_{6} & 0 \\
- Y_{0} & Y_{1} & Y_{2} & - Y_{3} & 0 & -2 Y_{5} & 2 Y_{6} & 0 & 0 \\
2 Y_{0} & 2 Y_{1} & -2 Y_{2} & -2 Y_{3} & 0 & 0 & 0 & 0 & 0
\end{pmatrix}.
\end{equation*}
As in the previous computation $\lvert\ccoeffiso{\lbrace \LieTwoOneOneNilp,\LieSubreg\rbrace}{(\spl{4}{(\qfield)}{})}^\element\rvert$ and $\lvert\ccoeffiso{\lbrace \LieTwoOneOneNilp,\LieReg\rbrace}{(\spl{4}{(\qfield)}{})}^\element\rvert$ are essentially the cardinalities of the rank-$4$ and rank-$6$ loci of the commutator matrix above. However now the centralizer $\arr$-classes with dimension $5$ are three, so it does not suffice to just apply \cref{lem:jp_rp_lifts} as we did before. We reason as follows: $\cmatrix_{\LieTwoOneOneNilp}(c_0,\dots,c_7,c_8)$ has rank $2$ if and only if $c_0=\dots=c_4=0$ and at least one of $c_5$, $c_6$ and $c_7$ is non-zero. This means that
\[
V_2 = \lbrace \tuple{x}\in\qfield^{9}\mid {\rk}_{\qfield} \cmatrix_{\LieTwoOneOneNilp}(\tuple{x}) \leq 2\rbrace
\]
is a Zariski closed set defined by the ideal $\radicalminors{2}=(Y_0,\dots,Y_4)\subseteq\qfield[Y_0,\dots,Y_8]$ and that the points giving rank $0$ in the commutator matrix constitute a Zariski closed set defined by the ideal $ (Y_0,\dots,Y_7)$. By looking at the submatrix of $\cmatrix_{\LieTwoOneOneNilp}$ corresponding to the last $4$ coordinates, we realize that $V_2$ is isomorphic to $\gl{2}{(\qfield)}$. 
As before we use \cref{lem:theta} to determine $\lvert\ccoeffiso{\lbrace \LieTwoOneOneNilp,\mathbf{s}\rbrace}{(\spl{4}{(\qfield)}{})}^\element\rvert$  for $\dimalg_\mathbf{s} = 7$. Each orbit of lifts of $\element$ to $\alg{2}$ having a $7$-dimensional centralizer contains an element of
\begin{equation}
\label{eq:param_for_211n_7dim}
C=\left\lbrace\left.\begin{pmatrix}
		0&1&0&0\\
		\pi\parameterOne&0&0&0\\
		0&0&0&\pi\\
		0&0&\pi\parameterOne&0\\
	     \end{pmatrix}\in\spl{4}{(\completion/\primeideal^2)}{}\,\right\vert\, \parameterOne\in\qfield \right\rbrace,
\end{equation}
 and vice versa distinct elements of $C$ are contained in distinct orbits. %
For each $\parameterOne\in\qfield$ let $y_\parameterOne \spl{4}{(\qfield)}{}$ be the matrix such that
\[
\element+\pi y_\alpha=\begin{pmatrix}
		0&1&0&0\\
		\pi\parameterOne&0&0&0\\
		0&0&0&\pi\\
		0&0&\pi\parameterOne&0\\
	     \end{pmatrix}\in \spl{4}{(\qfield)}{}
\]
We have that $\transpose{y_\alpha}\in\centr{\spl{4}{(\qfield)}{}}{\element}$ and that the Lie centralizer isomorphism type of
\[
\element+y_\alpha=\begin{pmatrix}
		0&1&0&0\\
		\parameterOne&0&0&0\\
		0&0&0&1\\
		0&0&\parameterOne&0\\
	     \end{pmatrix}
\] 
varies according to whether $\parameterOne$ is zero, a non-zero square or not a square in $\qfield$. Moreover each $\GL{4}{(\qfield)}$-orbit whose elements have $7$-dimensional Lie centralizer contains an element of the form $\element+y_\alpha$ for some $\alpha\in\qfield$, because $\lbrace \element+y_\alpha\mid \alpha\in\qfield\rbrace$ is the affine cross-section of $\sheet_{[2,2]}$. We thus compute
\begin{equation}
\begin{split}
\lvert\ccoeffiso{\lbrace \LieTwoOneOneNilp,\LieTwoTwoSemDiag\rbrace}{(\spl{4}{(\qfield)}{})}^\element\rvert	&=\frac{1}{2}q(q^2-1)\\
\lvert\ccoeffiso{\lbrace \LieTwoOneOneNilp,\LieTwoTwoSemNon\rbrace}{(\spl{4}{(\qfield)}{})}^\element\rvert	&=\frac{1}{2}q(q-1)^2\\
\lvert\ccoeffiso{\lbrace \LieTwoOneOneNilp,\LieTwoTwoNilp\rbrace}{(\spl{4}{(\qfield)}{})}^\element\rvert	&= q^2-1.
\end{split}
\end{equation}
The  same is valid for all the other elements of $\spl{4}{(\qfield)}{}$ that have centralizer isomorphic to $\LieTwoOneOneNilp$, for they all are $\GL{4}{(\qfield)}$-conjugates of $\element$.\par
In order to finish our investigation, we need to compute the cardinality of the rank-$4$ locus $L_4 = V_4\smallsetminus V_2$  of $\cmatrix_{\LieTwoOneOneNilp}$. We can do it by looking at its equations. The following is a generating set for the radical of the ideal generated by the $6\times 6$ principal minors of $\cmatrix_{\LieTwoOneOneNilp}$ that has been computed with SageMath \cite{sage2016}:
\begin{equation}
\label{eq:gens_rank_4}
\begin{split}
&Y_{0} Y_{3} -  Y_{4} Y_{5}\\ 
&Y_{1} Y_{2} -  Y_{4} Y_{6}\\ 
&Y_{0} Y_{2} -  Y_{1} Y_{3} -  Y_{4} Y_{7}\\ 
&Y_{2}^{2} Y_{5} -  Y_{3}^{2} Y_{6} -  Y_{2} Y_{3} Y_{7}\\ 
&Y_{1}^{2} Y_{5} -  Y_{0}^{2} Y_{6} + Y_{0} Y_{1} Y_{7}\\ 
&Y_{1} Y_{3}^{2} -  Y_{2} Y_{4} Y_{5} + Y_{3} Y_{4} Y_{7}\\ 
&Y_{1}^{2} Y_{3} -  Y_{0} Y_{4} Y_{6} + Y_{1} Y_{4} Y_{7}.
\end{split}
\end{equation}
Let $\radicalminors{4}$ be the ideal of $\qfield[Y_0,\dots,Y_8]$ generated by the polynomials in \cref{eq:gens_rank_4} and let $V_4$ be the algebraic set defined by it. The rank-$4$ locus $L_4$ is the set where all the polynomials in $\radicalminors{4}$ but not all the polynomials in $\radicalminors{2}=(Y_0,\dots,Y_4)$ vanish. Now let $\tuple{c}=(c_0,\dots,c_8)$ be a point of $L_4$. We notice that by forcing $c_0,\dots,c_4$ to $0$ we can project this point on $V_2$; this defines a function
\[\proj:\xymatrix@R=3pt{L_4	\ar[r]&	V_2\\
(c_0,\dots,c_8)\ar@{|->}[r]&(0,\dots0,c_5,\dots,c_8).}
\]
The set $V_2$ is stable under the action of $\centr{\SL{4}{(\qfield)}{}}{\element}$, so $\proj$ maps $\centr{\SL{4}{(\qfield)}{}}{\element}$-orbits to $\centr{\SL{4}{(\qfield)}{}}{\element}$-orbits and the cardinality of the fibres of $\proj$ is constant across  $\centr{\SL{4}{(\qfield)}{}}{\element}$-orbits in $V_2$. Let us identify $V_2$ with $\gl{2}{(\qfield)}$. In what follows we shall operate a case distinction according to the adjoint orbit in $\gl{2}{(\qfield)}$. The elements in the centre of $\gl{2}{(\qfield)}$ are those for which $c_5=\dots=c_7=0$.  Now we substitute the previous conditions in \cref{eq:gens_rank_4} and impose that at least one of the $c_0,\dots,c_4$ is non-zero (we want to exclude points of $V_2$ inside $V_4$). This gives that the fibre of $\proj$ above each one of these elements has cardinality
\[
2q^3-q-1.
\]
Now let us consider the elements in $\gl{2}{(\qfield)}$ that belong to a nilpotent orbit. These are the elements whose orbit contains an element defined by $c_5=1$, $c_6=c_7=0$ and $c_8$ arbitrary. Substituting these relations into the  \cref{eq:gens_rank_4} and imposing that at least one of the other variables is non-zero, we obtain that the fibre above a nilpotent point has cardinality
\[
q^2-1.
\]
The other orbits are parameterized by the following elements $c_5=1$, $c_7=0$ and $c_6=\parameterOne\in\qfield^\times$. Again by substituting we see that there is no point in $L_4$ projecting down to a point in an orbit with $\parameterOne$ a non-square in $\qfield$. It remains to compute the cardinality of the fibre above points for which $\parameterOne$ is a non-zero square in $\qfield$ (semisimple diagonalizable points). Substituting this condition in \cref{eq:gens_rank_4} and imposing that the other variables are not all zero, gives that the cardinality of the fibre of $\proj$ above each of these points is 
\[
2\cdot(q^2-1).
\]
Considered that in $\gl{2}{(\qfield)}$ there are $q$ central elements, $q\cdot(q^2-1)$ nilpotent elements and $q^2\cdot(q^2-1)/2$ semisimple diagonalizable points, we obtain that 
\begin{equation*}
\lvert L_4 \rvert=q\cdot(2q^3-q-1)+q\cdot(q^2-1)^2+q^2\cdot(q^2-1)^2=q\cdot(q^5+q^4-2q^2).
\end{equation*}
It follows that for $x\in  \LieTwoOneOneNilp$
\begin{equation}
\label{eq:211n}
\begin{split}
\lvert \ccoeffiso{\lbrace \LieTwoOneOneNilp,\LieSubreg\rbrace}{(\spl{4}{(\qfield)}{})}^x\rvert		=&q^5+q^4-2q^2\\
\lvert\ccoeffiso{\lbrace \LieTwoOneOneNilp,\LieReg\rbrace}{(\spl{4}{(\qfield)}{})}^x\rvert			=&q^8-1\\
																			&-\lvert\ccoeffiso{\lbrace \LieTwoOneOneNilp,\LieSubreg\rbrace}{(\spl{4}{(\qfield)}{})}^x\rvert\\																									&-\sum_{\stackrel{\mathbf{s}\in\arr}{\dimalg_\mathbf{s} = 7}} \lvert\ccoeffiso{\lbrace \LieTwoOneOneNilp,\mathbf{s}\rbrace}{(\spl{4}{(\qfield)}{})}^x\rvert\\
							    													=&q^8-q^3-q^2(q^3+q^2-2)\\
																				=&q^{8} - q^{5} - q^{4} - q^{3} + 2 q^{2}.
\end{split}
\end{equation}
We summarize all the previous computations in \cref{tab:jumps}.
\setlength{\tabcolsep}{10pt}
\begin{table}[h]
\caption{Data to apply \cref{eq:ind_step}.}
\label{tab:jumps}
\begin{tabular}{lcccl}
\toprule
$\mathbf{s}\in\arr$			&$\dimalg_\mathbf{s}$			&$\dimalg'_\mathbf{t}$			&$\mathbf{t}\in\arr$		&$\lvert\ccoeffiso{\lbrace \mathbf{s},\mathbf{t}\rbrace}{(\spl{4}{(\qfield)}{})}^\element\rvert$ $\forall\,\element \in \mathbf{s}$\\

\midrule
$\LieReg$					&$3$							&$0$							&n.a.\				&n.a.\ \\
\midrule
$\LieSubreg$				&$5$							&$3$							&$\LieReg$				&$q^3 -1$\\
\midrule
$\LieTwoTwoSemDiag$		&$7$							&$6$							&$\LieReg$				&$(q - 1)^2 (q^2 + q + 1)^2$\\
						&							&							&$\LieSubreg$				&$2  (q^3 -1)$\\
\cmidrule{2-5}
$\LieTwoTwoSemNon$		&$7$							&$6$							&$\LieReg$				&$q^6 - 1$\\
\cmidrule{2-5}
$\LieTwoTwoNilp$			&$7$							&$6$							&$\LieReg$				&$q^3 (q^3 - 1)$\\
						&							&							&$\LieSubreg$				&$q^3 -1$\\
\midrule
$\LieTwoOneOneSem$		&$9$							&$8$							&$\LieReg$				&$q (q - 1) (q^6 +q^5 +q^4 -q^2 -2 q - 1)$	\\
						&							&							&$\LieSubreg$				&$ q^5 + q^4 + q^3 - q^2 - q - 1$			\\
\cmidrule{2-5}
$\LieTwoOneOneNilp$		&$9$							&$8$							&$\LieReg$				&$q^{8} - q^{5} - q^{4} - q^{3} + 2 q^{2}$		\\
						&							&							&$\LieSubreg$				&$q^5+q^4-2q^2$				\\
						&							&							&$\LieTwoTwoSemDiag$		&$ \frac{1}{2}q (q^2-1)$					\\
						&							&							&$\LieTwoTwoSemNon$		&$ \frac{1}{2}q (q -1)^2$					\\
						&							&							&$\LieTwoTwoNilp$			&$q^2 - 1$							\\									
\bottomrule
\end{tabular}
\end{table}
-
\begin{table}[h]
\caption{Elements of $\spl{4}{(\qfield)}{}$ having prescribed Lie centralizer up to isomorphism.}
\label{tab:card}
\begin{tabular}{lc}
\toprule
Class $\mathbf{s}	$		&$\lvert \mathbf{s} \rvert$						\\
\midrule
$\LieTwoOneOneSem$			&$q^{3}\cdot(q-1)\cdot(q + 1) \cdot (q^{2} + 1)$						\\
$\LieTwoOneOneNilp$			&$(q - 1) \cdot (q + 1) \cdot (q^{2} + 1) \cdot (q^{2} + q + 1)$			\\
\midrule
$\LieTwoTwoSemDiag$			&$\frac{1}{2}\cdot q^{4}\cdot (q-1) \cdot (q^{2} + 1) \cdot (q^{2} + q + 1)$	\\
$\LieTwoTwoSemNon$			&$\frac{1}{2}\cdot q^{4} \cdot (q - 1)^{3} \cdot (q^{2} + q + 1)$\\
$\LieTwoTwoNilp$				&$q \cdot (q + 1) \cdot (q - 1)^{2} \cdot (q^{2} + 1) \cdot (q^{2} + q + 1)$	\\
\midrule
$\LieThreeOneTwoEV$			&$q^{4} \cdot(q - 1)^2 \cdot (q + 1) \cdot  (q^{2} + 1) \cdot (q^{2} + q + 1)$			\\
$\LieThreeOneNilp$				&$q^{2} \cdot(q - 1)^{2} \cdot  (q + 1)^{2} \cdot (q^{2} + 1) \cdot (q^{2} + q + 1)$		\\
$\LieThreeOneNonJor$			&$\frac{1}{2}\cdot q^{6} \cdot(q - 1)^2 \cdot  (q^{2} + 1) \cdot (q^{2} + q + 1)$		\\
$\LieThreeOneThreeEV$			&$q^{3} \cdot(q - 1)^2 \cdot (q + 1)^{2}  \cdot (q^{2} + 1) \cdot (q^{2} + q + 1)$		\\
$\LieThreeOneDiag$				&$\frac{1}{2}\cdot q^{5} \cdot(q-1)\cdot(q-2)\cdot (q + 1) \cdot  (q^{2} + 1) \cdot (q^{2} + q + 1)$	\\		
\midrule
$\LieReg$						&$(q - 1) \cdot (q + 1) \cdot q^{3} \cdot (q^{10} + q^{8} - q^{7} - 3 q^{5} - q^{3} + 2 q^{2} + q + 1)$\\
\bottomrule
\end{tabular}
\end{table}
\subsection{Number of elements by their Lie centralizers}
In order to use \cref{eq:ind_step} and apply \cref{main:D}, it remains to compute the cardinality of each class in $\arr$. Rather than repeating essentially the same computations for all classes of Lie centralizers we give just a quick description of how they are performed and then summarize the results in \cref{tab:card}. For each case in \cref{tab:dim_5,tab:dim_7,tab:dim_9,tab:dim_9_nilp} the size of the $\GL{4}{(\qfield)}$-orbit is the quotient of the cardinality of $\GL{4}{(\qfield)}$ and the cardinality of the $\GL{4}{(\qfield)}$-conjugation stabilizer. Alternatively, dividing both the numerator and the denominator by $q - 1$, the cardinality of the $\GL{4}{(\qfield)}{}$-orbit is the cardinality of $\SL{4}{(\qfield)}{}$ divided by the cardinality of the group centralizer. The latter is computed using the factorizations in the last column  of  \cref{tab:dim_5,tab:dim_7,tab:dim_9,tab:dim_9_nilp}, the cardinalities of $\qfield$, $\qfield^{\times}$, $C_{q+1}$ and the following well know cardinalities:
\begin{align*}
\lvert\GL{n}{(\qfield)}\rvert  									&=	\prod_{i = 0}^{n-1} (q^n - q^{i})		&& (n\in\N)\\
\lvert\SL{n}{(\qfield)}{}\rvert									&=	\frac{\lvert\GL{n}{(\qfield)}\rvert}{q-1}	&&\\
\lvert \heis{(\qfield)}\rvert										&=	q^3							&&\\
\lvert \heis{\qfield}\curlyvee\heis{\qfield}\rvert	&= 	q^5.							&&
\end{align*}
This takes care of the cardinality of each orbit. The conditions on the parameters of the cross-section that define each class in $\arr$ give the number of orbits that class consists of. We have taken care not to put  together orbits giving non-isomorphic group centralizers. This means that all the orbits in the same $\arr$-class have the same cardinality. Thus, in order to compute the cardinality of an $\arr$-class, say $\mathbf{s}$, it suffices to multiply the cardinality of an orbit in $\mathbf{s}$ by the number of distinct orbits that form $\mathbf{s}$. This latter quantity is obtained directly from the conditions on the cross-section's parameters reported in the first column of  \cref{tab:dim_5,tab:dim_7,tab:dim_9,tab:dim_9_nilp}. Finally, as anticipated, the number $\lvert \LieReg\rvert$ is deduced by subtraction.
\section{Representation zeta functions}
\label{sec:sl4}
We now have all the necessary numerical information to compute the representation zeta function of the principal congruence subgroups of $\SL{4}{(\completion)}{}$ using \cref{main:D} (\Cref{tab:jumps} gives a summary of the relevant numerical information). 
A lengthy but straightforward application of \cref{eq:ind_step} to the data in \cref{tab:card,tab:jumps} allows  us to compute all the relevant ingredients in the formula of \cref{main:D}, which gives that, for $2\nmid q$, the Poincar\'e series of $\spl{4}{(\completion)}{}$ is
\[
\pseries{\spl{4}{(\completion)}}{(s)} = \frac{\numpoin{Poin}{q,t}}{\denpoin{Poin}{q,t}}
\]
where
\allowdisplaybreaks
\begin{align*}
\label{P_series_num}
\numpoin{Poin}{q,t}=&q^{28} t^{18}\\
 		      &- {\left(q^{28} + q^{27} + q^{26} + q^{25} - q^{24} - q^{23} - q^{22}\right)} t^{15}\\
		      &+ {\left(q^{27} - 2 \, q^{24} - q^{22} + q^{21}\right)} t^{14}\\
		      & + {\left(q^{26} + 2 \, q^{25} + 2 \, q^{24} - 2 \, q^{22} - 4 \, q^{21} - 2 \, q^{20} - q^{19} + 2 \, q^{18} + q^{17}\right)} t^{13}\\
		      & - {\left(q^{25} + q^{24} + q^{23} - 2 \, q^{22} - 2 \, q^{21} - 2 \, q^{20} + 2 \, q^{18} + q^{17} + q^{16}\right)} t^{12} \\
		      &+ {\left(q^{21} + 2 \, q^{19} + q^{17} - q^{16} - q^{15} - q^{14}\right)} t^{11} \\
		      &+ {\left(q^{19} + q^{18} - 2 \, q^{15} + q^{12}\right)} t^{10} \\
		      &- \left(2 \, q^{19} + q^{18} + q^{17} - q^{16} - 3 \, q^{15} - 2 \, q^{14}\right.\\
		      &\left. - 3 \, q^{13} - q^{12} + q^{11} + q^{10} + 2 \, q^{9}\right) t^{9}\\
		      & + {\left(q^{16} - 2 \, q^{13} + q^{10} + q^{9}\right)} t^{8}\\
		      & - {\left(q^{14} + q^{13} + q^{12} - q^{11} - 2 \, q^{9} - q^{7}\right)} t^{7}\\
		      & - {\left(q^{12} + q^{11} + 2 \, q^{10} - 2 \, q^{8} - 2 \, q^{7} - 2 \, q^{6} + q^{5} + q^{4} + q^{3}\right)} t^{6}\\
		      & + {\left(q^{11} + 2 \, q^{10} - q^{9} - 2 \, q^{8} - 4 \, q^{7} - 2 \, q^{6} + 2 \, q^{4} + 2 \, q^{3} + q^{2}\right)} t^{5} \\
		      &+ {\left(q^{7} - q^{6} - 2 \, q^{4} + q\right)} t^{4} \\
		      &+ {\left(q^{6} + q^{5} + q^{4} - q^{3} - q^{2} - q - 1\right)} t^{3} \\
		      &+ 1\\
\denpoin{Poin}{q,t}=&{\left(1-q^{7}t^3\right)} {\left(1-q^{9}t^4\right)}  {\left( 1-q^{12}t^5 \right)} {\left(1-q^{15}t^6 \right)}.
\end{align*}
Operating the substitution dictated by the Kirillov orbit method in \cite[Proposition 3.1]{akov2013representation}, we deduce \cref{thmB}.
\part{$\spl{4}{(\completion)}{}$ is not shadow-preserving}
\label{part:hered}
\section{Lie shadows and centralizer-loci}
In this part we prove \cref{main:G}. Recall that $\qfield$ has characteristic $p\neq 2$. Fix $\level\in\N$. For convenience of notation we set $\completion_\level = \completion / \primeideal^\level$ and $\alg{\level} = \spl{4}{(\completion_\level)}{}$. By the correspondence between Lie and group shadows proved in \cite[Lemma 2.3]{akov2}, $x\in\spl{4}{(\completion_\level)}{}$ has a shadow-preserving lift if and only if it has a Lie shadow-preserving lift. Indeed in what follows we shall show that there are levels $\level\in\N$ and $x\in\alg{\level}$ that have no Lie shadow-preserving lift to $\alg{\level + 1}$.\par
Let $e = \dim_{\qfield}\sh{\alg{\level}}{x}$. We rephrase the problem of finding shadow-preserving lifts in terms of lifting points of $\centrvar{\spl{4}{(\completion)}{}}{e}{(\completion_{\level})}$ to points in $\centrvar{\spl{4}{(\completion)}{}}{e}{(\completion_{\level + 1})}$. Let $\iota_\level: \completion_\level^{15}\rightarrow \alg{\level}$ be the isomorphism of $\completion_\level$-modules determined by fixing the basis $\basis$ as in \cref{sec:killing}. We have the following lemma.
\begin{lem}
\label{lem:shadows_centr_locus}
Let $x\in\alg{\level}$. Then $e = \dim_{\qfield}\sh{\alg{\level}}{x}$ if and only if 
\[\iota_\level^{-1}(x)\in \centrvar{\spl{4}{(\completion)}{}}{e}{(\completion_{\level})}.\]
\begin{proof}
Recall that we denote by $\cmatrix$ the commutator matrix of $\spl{4}{(\completion)}{}$ with respect to $\basis$. 
Take $\tuple{x}\in\completion_\level^{15}$, we show that $\tuple{x} \in \centrvar{\spl{4}{(\completion)}{}}{e}{(\completion_{\level})}$ is equivalent to $\iota(x)$ having $e$-dimensional shadow.\par
Let $\widehat{\tuple{x}}\in\completion^{15}$ be a lift of $\tuple{x}$ (any lift will do here). Let $\xi$ be as defined in \cref{sec:quadratic}. Then $\widehat{\tuple{x}}\in\xi(\ranklocus{\cmatrix}{15 - e}{(\completion)})$ modulo $\primeideal^\level$. This is equivalent to the last $e$ elementary divisors of $\cmatrix(\xi^{-1}(\widehat{\tuple{x}}))$ being divisible by $\pi^\level$ and the first $15 - e$ not being divisible by $\pi^\level$. This in turn is equivalent to 
\[
\dim_{\qfield} \reduc{\level}{1}(\ker \cmatrix^\level(\xi^{-1}_\level(\tuple{x}))) = e,
\]
where $\xi_\level: \completion_\level^{15}\rightarrow \completion_\level^{15}$ is induced by $\xi$.\par
We may now conclude that the latter condition is equivalent to $\sh{\alg{\level}}{x}$ being $e$-dimensional because, by an argument analogous to the one preceding \cref{def:geo_smooth_centr}, $\ker\cmatrix^\level(\xi^{-1}_\level(\tuple{x}))$ and $\centr{\alg{\level}}{x}$ are isomorphic as $\completion_\level$-Lie lattices.
\end{proof}
\end{lem}
Now let $\level\in\N$ and  $\element\in\alg{\level}$ with Lie shadow $\lieshadow$ of $\qfield$-dimension $e$.  Let also $b\in\alginf$ be a lift of $\element$. It is easy to see that all lifts of $\element$ to $\alg{\level+1}$ have the form $\reducinf{\level+1}(b + \pi^{\level+1} c)$  for some $c\in\alginf$. \Cref{lem:jp_rp_lifts} has the following straightforward consequence:
\begin{prop}
\label{prop:shadows_lifts}
If $\element\in\alg{\level}$ has a shadow-preserving lift $b\in\alginf$, then the number of shadow-preserving lifts of $\element$ to $\alg{\level + 1}$ is 
\[
q^{\dimalg - \dim_{\qfield} [\lieshadow,\lieshadow]}.
\]
Moreover, let $\basis_\lieshadow$ be a basis of $\lieshadow$ and let $\bar{\iota}:\qfield^{\dim_{\qfield}\lieshadow}\rightarrow \lieshadow$ be the associated $\qfield$-linear isomorphism. Let also $\redcmatrix$ be the $\qfield$-commutator matrix of $\lieshadow$ with respect to $\basis_\lieshadow$ and $\basis'$ be a basis of $\alginf$ obtained by completing a lift of $\basis_\lieshadow$. Then for all $c\in\alginf$ the Lie shadow of $\reducinf{\level+1}(b + \pi^{\level+1} c)$ is $\bar{\iota}(\ker\redcmatrix(\tuple{x}_c))$ where $\tuple{x}_c$ are the first $e$ coordinates of $\transpose{c}$ with respect to $\basis'$.
\end{prop}
%
\newcommand{\LieTwoTwoNilpBasis}{\mathcal{\expandafter\uppercase\expandafter{\TwoTwoNilpLetter}}}
\section{Proof of \cref{main:G}}
We keep the notation of \cref{prop:shadows_lifts}. Let $x = e_{12} +e_{34}\in\spl{4}{(\qfield)}{}$. The first row of \cref{tab:dim_7} shows the structure of the Lie centralizer of $x$ and also gives an $\qfield$-basis for it. Let $\mathfrak{\TwoTwoNilpLetter} = \centr{\spl{4}{(\qfield)}{}}{x}$ and $\LieTwoTwoNilpBasis = \lbrace c_0,\dots,c_6\rbrace$ be its basis as reported in \cref{tab:dim_7}. According to the commutator relations (also reported in \cref{tab:dim_7}) we compute that the $\qfield$-commutator matrix of $\mathfrak{\TwoTwoNilpLetter}$ with respect to $\LieTwoTwoNilpBasis $ is 
\begin{equation}
\label{eq:22Nilp}
\redcmatrix_{\LieTwoTwoNilpBasis} (Y_0,\dots Y_6) = \begin{pmatrix}
								0 		& Y_{2} 	& -2 Y_{0} 	& -2 Y_{4} 	& 0 		& Y_{3}	&0 \\
								- Y_{2} 	& 0 		& 2 Y_{1} 		& 2 Y_{5} 		& - Y_{3} 	& 0		&0 \\
								2 Y_{0} 	& -2 Y_{1}	& 0 			& 0 			& 2 Y_{4} 	& -2 Y_{5}&0 \\
								2 Y_{4} 	& -2 Y_{5}	& 0 			& 0 			& 0 		& 0		&0 \\
								0 		& Y_{3} 	& -2 Y_{4} 	& 0 			& 0 		& 0		&0 \\
								- Y_{3} 	& 0 		& 2 Y_{5} 		& 0 			& 0 		& 0		&0\\
								0		&0		&0			&0			&0		&0		&0
\end{pmatrix}.
\end{equation}
Take now $\widehat{x}= e_{12} + e_{34}\in \spl{4}{(\completion)}{}$, it is easy to see that $\widehat{x}$ has the same shadow as $x$. Let $y\in\spl{4}{(\qfield)}{}$ be such that $\transpose{y} = \alpha_0 c_0 + \alpha_1 c_1 + \alpha_2 c_2 \in\mathfrak{\TwoTwoNilpLetter}$ with $\alpha_0,\alpha_1,\alpha_2\in\qfield$ and not all $0$ at the same time. Let also $\widehat{y}\in\spl{4}{(\completion)}{}$ be a lift of $y$ and  $\element = \reducinf{2}(\widehat{x} + \pi \widehat{y})$. Then, by \cref{prop:shadows_lifts}, $\sh{\spl{4}{(\completion_2)}{}}{\element} = \lbrace v,c_3,\dots,c_6\rbrace$ where $v$ spans the kernel of 
\begin{equation*}
\begin{pmatrix}
								0 		& \alpha_{2} 	& -2 \alpha_{0} 	&0 	& 0 	& 0	&0 \\
								- \alpha_{2} 	& 0 		& 2 \alpha_{1} 		& 0 	& 0	& 0	&0 \\
								2 \alpha_{0} 	& -2 \alpha_{1}	& 0 			& 0 	& 0 	& 0	&0 
\end{pmatrix}.
\end{equation*}
Using the commutator relations in \cref{tab:dim_7}, we compute
\[
[v,c_3]= -2\alpha_0 c_4+2\alpha_1 c_5,\,
[v,c_4]= -\alpha_1 c_3+2\alpha_2 c_4,\,
[v,c_5]= \alpha_0 c_3-2\alpha_2 c_5.
\]
One checks that the three commutators above span a $2$-dimensional space for all possible choices of $\alpha_0,\alpha_1,\alpha_2$ as above. Thus if $\element$ admits a shadow-preserving lift $b\in\spl{4}{(\completion)}{}$, \cref{prop:shadows_lifts} tells us that there are exactly $q^{15 - 2} = q^{13}$ shadow-preserving lifts of $\element$ to $\spl{4}{(\completion)}{}$.\par
We specialize now the above discussion for $\element = x + \pi (e_{13} + e_{24})\in\spl{4}{(\completion_2)}{}$. It is straightforward to see that $b = e_{12} + e_{34} +  \pi (e_{13} + e_{24}) \in \spl{4}{(\completion)}{}$ is a shadow-preserving lift of $\element$. So \cref{prop:shadows_lifts} implies that there are exactly $q^{13}$ shadow preserving lifts of $\reducinf{3}(b)$ to $\spl{4}{(\completion)}{}$.\par
On the other hand, it can be shown by direct computation that $\mathfrak{n} = \centr{\spl{4}{(\completion)}{}}{\transpose{b}}$ has $\completion$-rank $5$, thus $b\in\ranklocus{\cmatrix}{10}{(\completion)}$. The same computation also gives that $[\mathfrak{n},\mathfrak{n}] \cong \completion\oplus\completion\oplus \pi\completion$ as an $\completion$-module. Thus, if $\mathrm{b}$ are the coordinates of $b$ with respect to $\basis$, by \cref{prop:rk_jacobian},
\[
\det\mathcal{J}_{10} (\tuple{b}) = \pi.
\]
This makes it a smooth point of $\ranklocus{\cmatrix}{10}{(\localfrac)}$, which is then locally cut out by $15 - 12 = 3$ polynomials that may be taken to have coefficients in $\completion$. In other words there are $g,f_1,f_2,f_3\in\completion[X_1,\dots,X_{15}]$ and a Zariski open neighbourhood $\mathcal{U} = D(g)(\localfrac)$ of $\tuple{b}$ such that 
\[
\mathcal{U}\cap (\rankvar{\cmatrix}{10}{(\localfrac)}) \cong \spm \left( \completion[X_1,\dots,X_{15}, X']_g/(f_1,f_2,f_3, 1 - g X')\right)
\]
Also notice that by definition $D(g)(\completion)$ contains all lifts of $\reducinf{3}(\tuple{b})$ because, for all $\tuple{b}'\in\completion^{15}$ such that $\tuple{b}'\equiv\tuple{b}\,\bmod\, \primeideal^3$ we have
\[
g(\tuple{b}')	\equiv	g(\tuple{b}) \mod \primeideal^3.
\]
This means that we may safely restrict our investigation to the neighbourhood $\mathcal{U}$ without fear of missing any  lift of $\reducinf{3}(\tuple{b})$.\par
Let $\tuple{b} = (b_1,\dots,b_{15})$. By \cite[\S 4.6, Corollary 3]{bourbaki1972commutative} there are exactly $3$ formal power series without constant term, say $\phi_1,\phi_2,\phi_3\in \completion\llbracket X_1,\dots, X_{12} \rrbracket$, such that, for all $\tuple{t} = (t_1,\dots,t_{12}) \in \primeideal^{(12)}$,
\[
f_i(b_1 + \pi^2 t_1,\dots, b_{12} + \pi^2 t_{12}, b_{13} + \pi\phi_1(\tuple{t}), b_{14} + \pi\phi_2(\tuple{t}), b_{15} + \pi\phi_3(\tuple{t})) = 0
\]
for all $i=1,2,3$. By \cref{lem:shadows_centr_locus}, each 
\[z(\tuple{t}) = (b_1 + \pi^2 t_1,\dots, b_{12} + \pi^2 t_{12}, b_{13} + \pi\phi_1(\tuple{t}), b_{14} + \pi\phi_2(\tuple{t}), b_{15} + \pi\phi_3(\tuple{t})) 
\]
is a shadow-preserving lift of $\reducinf{3}(b)$ to $\spl{4}{(\completion)}{}$. Since $\lbrace \reducinf{4}(z(\tuple{t}))\mid \tuple{t}\in\primeideal^{12}\rbrace$ has cardinality $q^{12}$, we conclude that there must be $q^{13}-q^{12}$ shadow-preserving lifts of $\reducinf{3}(b)$ that do not admit shadow-preserving lifts in $\spl{4}{(\completion)}{}$ in their turn. This concludes the proof of \cref{main:G}. The following example gives an explicit element that we have checked not to have shadow-preserving lifts
\begin{exmp}
\label{ex:non_shadow}
Take the following matrix in $\spl{4}{(\Z/27\Z)}{}$
\begin{equation*}
\label{eq:nsp_lift}
z=\left(\begin{array}{rrrr}
9 & 10 & 21 & 0 \\
0 & 18 & 9 & 21 \\
0 & 9 & 0 & 10 \\
0 & 0 & 0 & 0
\end{array}\right).
\end{equation*}
By the discussion above, this is a lift of $e_{12} + e_{34} + 3 (e_{13} + e_{32})\in\spl{4}{(\Z/9\Z)}{}$. Whether a lift of $z$ has shadow-preserving lifts or not may be tested by running through its $3^{15}$ lifts to $\spl{4}{(\Z/81\Z)}{}$ (viewed in $\Z^{15}$) taking the transpose and plugging in the result in the commutator matrix $\cmatrix$ (viewed as a matrix of linear polynomials over $\Z$). If the result has $5$ elementary divisors that are divisible by $3^5$ then the lift is shadow-preserving otherwise it is not. We have done this using SageMath \cite{sage2016} which provides an algorithm to compute the Smith normal form of a matrix. Due to a memory-leak in the said algorithm we could not run though all the $3^{15}$ lifts at once. We have however run $3^{5}$ times through randomly chosen sample sets of the lifts, each of size of size $3^7$. The code is available upon request to the author.
\end{exmp}



\end{document}